\documentclass{amsart}
\usepackage{amssymb}

\usepackage{amscd}
\usepackage{thmdefs}

\input{tcilatex}
\theoremstyle{definition}
\theoremstyle{remark}
\numberwithin{equation}{section}

\input tcilatex

\begin{document}
\title[Generating Operator of $L(F_{2})*_{L(F_{1})}L(F_{2})$]{The Moments of the Generating Operator of $L(F_{2})*_{L(F_{1})}L(F_{2})$}
\author{Ilwoo Cho}
\address{Dep. of Math, Univ. of Iowa, Iowa City, IA, U. S. A}
\email{ilcho@math.uiowa.edu}
\maketitle

\begin{abstract}
In this paper, we will consider an example of a (scalar-valued) moment
series, under the compatibility. Suppose that we have an amalgamated free
product of free group algebras, $%
L(F_{2})*_{L(F_{1})}L(F_{2})=L(<a,b>)*_{L(<h>)}L(<c,d>)$

We will provide the method how to find the moment series of $%
a+b+a^{-1}+b^{-1}+c+d+c^{-1}+d^{-1}.$ Amalgamated freeness of $%
a+b+a^{-1}+b^{-1}$ and $c+d+c^{-1}+d^{-1}$ over $L(F_{1})$ is used and some
combinatorial functions (to explain the recurrence relations) are used to
figure out the $n$-th moment of this element.
\end{abstract}

\strut

\strut

Voiculescu developed Free Probability Theory. Here, the classical concept of
Independence in Probability theory is replaced by a noncommutative analogue
called Freeness (See [9]). There are two approaches to study Free
Probability Theory. One of them is the original analytic approach of
Voiculescu and the other one is the combinatorial approach of Speicher and
Nica (See [23], [1] and [24]).\medskip

Speicher defined the free cumulants which are the main objects in
Combinatorial approach of Free Probability Theory. And he developed free
probability theory by using Combinatorics and Lattice theory on collections
of noncrossing partitions (See [24]). Also, Speicher considered the
operator-valued free probability theory, which is also defined and observed
analytically by Voiculescu, when $\Bbb{C}$ is replaced to an arbitrary
algebra $B$ (See [23]). Nica defined R-transforms of several random
variables (See [1]). He defined these R-transforms as multivariable formal
series in noncommutative several indeterminants. To observe the R-transform,
the M\"{o}bius Inversion under the embedding of lattices plays a key
role.\strut

\strut

In [16], we observed the amalgamated R-transform calculus. Actually,
amalgamated R-transforms are defined originally by Voiculescu and are
characterized combinatorially by Speicher (See [23]). In [16], we defined
amalgamated R-transforms slightly differently from those defined in [23] and
[13]. We defined them as $B$-formal series and tried to characterize, like
in [1] and [24].

\strut

In [15], we observed the compatibility of a noncommutative probability space
and an amalgamated noncommutative probability space over an unital algebra.
In this paper, we have a nice compatibility of $\left(
L(F_{2})*_{L(F_{1})}L(F_{2}),\varphi \right) $ and $\left(
L(F_{2})*_{L(F_{1})}L(F_{2}),E\right) ,$ where \strut $%
tr:L(F_{2})*_{L(F_{1})}L(F_{2})\rightarrow \Bbb{C}$ is the canonical trace
and $F$ is the free product of conditional expectations $E:L(F_{2})%
\rightarrow L(F_{1}).$

\strut

In this paper, we will compute the $n$-th (scalar-valued) moment

\strut

\begin{center}
$a+b+a^{-1}+b^{-1}+c+d+c^{-1}+d^{-1}\in L(F_{2})*_{L(F_{1})}L(F_{2}),$
\end{center}

\strut

where $<a,b>\,=F_{2}=\,<c,d>.$

\strut

\strut

\strut

\section{Preliminaries}

\strut

\strut

\subsection{Amalgamated Free Probability}

\strut

In this section, we will summarize and introduced the basic results from
[23] and [16]. Throughout this section, let $B$ be a unital algebra. The
algebraic pair $(A,\varphi )$ is said to be a noncommutative probability
space over $B$ (shortly, NCPSpace over $B$) if $A$ is an algebra over $B$
(i.e $1_{B}=1_{A}\in B\subset A$) and $\varphi :A\rightarrow B$ is a $B$%
-functional (or a conditional expectation) ; $\varphi $ satisfies

\begin{center}
$\varphi (b)=b,$ for all $b\in B$
\end{center}

and

\begin{center}
$\varphi (bxb^{\prime })=b\varphi (x)b^{\prime },$ for all $b,b^{\prime }\in
B$ and $x\in A.$
\end{center}

\strut

\strut Let $(A,\varphi )$ be a NCPSpace over $B.$ Then, for the given $B$%
-functional, we can determine a moment multiplicative function $\widehat{%
\varphi }=(\varphi ^{(n)})_{n=1}^{\infty }\in I(A,B),$ where

\strut

\begin{center}
$\varphi ^{(n)}(a_{1}\otimes ...\otimes a_{n})=\varphi (a_{1}....a_{n}),$
\end{center}

\strut \strut

for all $a_{1}\otimes ...\otimes a_{n}\in A^{\otimes _{B}n},$ $\forall n\in 
\Bbb{N}.$

\strut

\strut We will denote noncrossing partitions over $\{1,...,n\}$ ($n\in \Bbb{N%
}$) by $NC(n).$ Define an ordering on $NC(n)$ ;

\strut

$\theta =\{V_{1},...,V_{k}\}\leq \pi =\{W_{1},...,W_{l}\}$\strut $\overset{%
def}{\Leftrightarrow }$ For each block $V_{j}\in \theta $, there exists only
one block $W_{p}\in \pi $ such that $V_{j}\subset W_{p},$ for $j=1,...,k$
and $p=1,...,l.$

\strut

Then $(NC(n),\leq )$ is a complete lattice with its minimal element $%
0_{n}=\{(1),...,(n)\}$ and its maximal element $1_{n}=\{(1,...,n)\}$. We
define the incidence algebra $I_{2}$ by a set of all complex-valued\
functions $\eta $ on $\cup _{n=1}^{\infty }\left( NC(n)\times NC(n)\right) $
satisfying $\eta (\theta ,\pi )=0,$ whenever $\theta \nleq \pi .$ Then,
under the convolution

\begin{center}
$*:I_{2}\times I_{2}\rightarrow \Bbb{C}$
\end{center}

defined by

\begin{center}
$\eta _{1}*\eta _{2}(\theta ,\pi )=\underset{\theta \leq \sigma \leq \pi }{%
\sum }\eta _{1}(\theta ,\sigma )\cdot \eta _{2}(\sigma ,\pi ),$
\end{center}

$I_{2}$ is indeed an algebra of complex-valued functions. Denote zeta,
M\"{o}bius and delta functions in the incidence algebra $I_{2}$ by $\zeta ,$ 
$\mu $ and $\delta ,$ respectively. i.e

\strut

\begin{center}
$\zeta (\theta ,\pi )=\left\{ 
\begin{array}{lll}
1 &  & \theta \leq \pi \\ 
0 &  & otherwise,
\end{array}
\right. $
\end{center}

\strut

\begin{center}
$\delta (\theta ,\pi )=\left\{ 
\begin{array}{lll}
1 &  & \theta =\pi \\ 
0 &  & otherwise,
\end{array}
\right. $
\end{center}

\strut

and $\mu $ is the ($*$)-inverse of $\zeta .$ Notice that $\delta $ is the ($%
* $)-identity of $I_{2}.$ By using the same notation ($*$), we can define a
convolution between $I(A,B)$ and $I_{2}$ by

\strut

\begin{center}
$\widehat{f}\,*\,\eta \left( a_{1},...,a_{n}\,;\,\pi \right) =\underset{\pi
\in NC(n)}{\sum }\widehat{f}(\pi )(a_{1}\otimes ...\otimes a_{n})\eta (\pi
,1_{n}),$
\end{center}

\strut

where $\widehat{f}\in I(A,B)$, $\eta \in I_{1},$ $\pi \in NC(n)$ and $%
a_{j}\in A$ ($j=1,...,n$), for all $n\in \Bbb{N}.$ Notice that $\widehat{f}%
*\eta \in I(A,B),$ too. Let $\widehat{\varphi }$ be a moment multiplicative
function in $I(A,B)$ which we determined before. Then we can naturally
define a cumulant multiplicative function $\widehat{c}=(c^{(n)})_{n=1}^{%
\infty }\in I(A,B)$ by

\begin{center}
$\widehat{c}=\widehat{\varphi }*\mu $ \ \ \ or \ \ $\widehat{\varphi }=%
\widehat{c}*\zeta .$
\end{center}

This says that if we have a moment\ multiplicative function, then we always
get a cumulant multiplicative function and vice versa, by $(*).$ This
relation is so-called ''M\"{o}bius Inversion''. More precisely, we have

\strut

\begin{center}
$
\begin{array}{ll}
\varphi (a_{1}...a_{n}) & =\varphi ^{(n)}(a_{1}\otimes ...\otimes a_{n}) \\ 
& =\underset{\pi \in NC(n)}{\sum }\widehat{c}(\pi )(a_{1}\otimes ...\otimes
a_{n})\zeta (\pi ,1_{n}) \\ 
& =\underset{\pi \in NC(n)}{\sum }\widehat{c}(\pi )(a_{1}\otimes ...\otimes
a_{n}),
\end{array}
$
\end{center}

\strut

for all $a_{j}\in A$ and $n\in \Bbb{N}.$ Or equivalently,

\strut

\begin{center}
$
\begin{array}{ll}
c^{(n)}(a_{1}\otimes ...\otimes a_{n}) & =\underset{\pi \in NC(n)}{\sum }%
\widehat{\varphi }(\pi )(a_{1}\otimes ...\otimes a_{n})\mu (\pi ,1_{n}).
\end{array}
$
\end{center}

\strut \strut

Now, let $(A_{i},\varphi _{i})$ be NCPSpaces over $B,$ for all $i\in I.$
Then we can define a amalgamated free product of $A_{i}$ 's and amalgamated
free product of $\varphi _{i}$'s by

\begin{center}
$A\equiv *_{B}A_{i}$ \ \ and \ $\varphi \equiv *_{i}\varphi _{i},$
\end{center}

respectively. Then, by Voiculescu, $(A,\varphi )$ is again a NCPSpace over $%
B $ and, as a vector space, $A$ can be represented by

\begin{center}
\strut $A=B\oplus \left( \oplus _{n=1}^{\infty }\left( \underset{i_{1}\neq
...\neq i_{n}}{\oplus }(A_{i_{1}}\ominus B)\otimes ...\otimes
(A_{i_{n}}\ominus B)\right) \right) ,$
\end{center}

where $A_{i_{j}}\ominus B=\ker \varphi _{i_{j}}.$ We will use Speicher's
combinatorial definition of amalgamated free product of $B$-functionals ;

\strut

\begin{definition}
Let $(A_{i},\varphi _{i})$ be NCPSpaces over $B,$ for all $i\in I.$ Then $%
\varphi =*_{i}\varphi _{i}$ is the amalgamated free product of $B$%
-functionals $\varphi _{i}$'s on $A=*_{B}A_{i}$ if the cumulant
multiplicative function $\widehat{c}=\widehat{\varphi }*\mu \in I(A,B)$ has
its restriction to $\underset{i\in I}{\cup }A_{i},$ $\underset{i\in I}{%
\oplus }\widehat{c_{i}},$ where $\widehat{c_{i}}$ is the cumulant
multiplicative function induced by $\varphi _{i},$ for all $i\in I$ and, for
each $n\in \Bbb{N},$

\strut 

\begin{center}
$c^{(n)}(a_{1}\otimes ...\otimes a_{n})=\left\{ 
\begin{array}{lll}
c_{i}^{(n)}(a_{1}\otimes ...\otimes a_{n}) &  & \text{if }\forall a_{j}\in
A_{i} \\ 
0_{B} &  & otherwise.
\end{array}
\right. $
\end{center}
\end{definition}

\strut

Now, we will observe the freeness over $B$ ;

\strut

\begin{definition}
Let $(A,\varphi )$\strut be a NCPSpace over $B.$

\strut 

(1) Subalgebras containing $B,$ $A_{i}\subset A$ ($i\in I$) are free (over $B
$) if we let $\varphi _{i}=\varphi \mid _{A_{i}},$ for all $i\in I,$ then $%
*_{i}\varphi _{i}$ has its cumulant multiplicative function $\widehat{c}$
such that its restriction to $\underset{i\in I}{\cup }A_{i}$ is $\underset{%
i\in I}{\oplus }\widehat{c_{i}},$ where $\widehat{c_{i}}$ is the cumulant
multiplicative function induced by each $\varphi _{i},$ for all $i\in I.$

\strut 

(2) Sebsets $X_{i}$ ($i\in I$) are free (over $B$) if subalgebras $A_{i}$'s
generated by $B$ and $X_{i}$'s are free in the sense of (1). i.e If we let $%
A_{i}=A\lg \left( X_{i},B\right) ,$ for all $i\in I,$ then $A_{i}$'s are
free over $B.$
\end{definition}

\strut

In [23], Speicher showed that the above combinatorial freeness with
amalgamation can be used alternatively with respect to Voiculescu's original
freeness with amalgamation.

\strut

Let $(A,\varphi )$ be a NCPSpace over $B$ and let $x_{1},...,x_{s}$ be $B$%
-valued random variables ($s\in \Bbb{N}$). Define $(i_{1},...,i_{n})$-th
moment of $x_{1},...,x_{s}$ by

\strut

\begin{center}
$\varphi (x_{i_{1}}b_{i_{2}}x_{i_{2}}...b_{i_{n}}x_{i_{n}}),$
\end{center}

\strut

for arbitrary $b_{i_{2}},...,b_{i_{n}}\in B,$ where $(i_{1},...,i_{n})\in
\{1,...,s\}^{n},$ $\forall n\in \Bbb{N}.$ Similarly, define a symmetric $%
(i_{1},...,i_{n})$-th moment by the fixed $b_{0}\in B$ by

\strut

\begin{center}
$\varphi (x_{i_{1}}b_{0}x_{i_{2}}...b_{0}x_{i_{n}}).$
\end{center}

\strut

If $b_{0}=1_{B},$ then we call this symmetric moments, trivial moments.

\strut

Cumulants defined below are main tool of combinatorial free probability
theory ; in [16], we defined the $(i_{1},...,i_{n})$-th cumulant of $%
x_{1},...,x_{s}$ by

\strut

\begin{center}
$k_{n}(x_{i_{1}},...,x_{i_{n}})=c^{(n)}(x_{i_{1}}\otimes
b_{i_{2}}x_{i_{2}}\otimes ...\otimes b_{i_{n}}x_{i_{n}}),$
\end{center}

\strut

for $b_{i_{2}},...,b_{i_{n}}\in B,$ arbitrary, and $(i_{1},...,i_{n})\in
\{1,...,s\}^{n},$ $\forall n\in \Bbb{N},$ where $\widehat{c}%
=(c^{(n)})_{n=1}^{\infty }$ is the cumulant multiplicative function induced
by $\varphi .$ Notice that, by M\"{o}bius inversion, we can always take such 
$B$-value whenever we have $(i_{1},...,i_{n})$-th moment of $%
x_{1},...,x_{s}. $ And, vice versa, if we have cumulants, then we can always
take moments. Hence we can define a symmetric $(i_{1},...,i_{n})$-th\
cumulant by $b_{0}\in B$ of $x_{1},...,x_{s}$ by

\strut

\begin{center}
$k_{n}^{symm(b_{0})}(x_{i_{1}},...,x_{i_{n}})=c^{(n)}(x_{i_{1}}\otimes
b_{0}x_{i_{2}}\otimes ...\otimes b_{0}x_{i_{n}}).$
\end{center}

\strut

If $b_{0}=1_{B},$ then it is said to be trivial cumulants of $%
x_{1},...,x_{s} $.

\strut

By Speicher, it is shown that subalgebras $A_{i}$ ($i\in I$) are free over $%
B $ if and only if all mixed cumulants vanish.

\strut

\begin{proposition}
(See [23] and [16]) Let $(A,\varphi )$ be a NCPSpace over $B$ and let $%
x_{1},...,x_{s}\in (A,\varphi )$ be $B$-valued random variables ($s\in \Bbb{N%
}$). Then $x_{1},...,x_{s}$ are free if and only if all their mixed
cumulants vanish. $\square $
\end{proposition}

\strut

\begin{remark}
The above noncommutative probability space with amalgamation can be replaced
by $W^{*}$-probability space with amalgamation and later, we will use the $%
W^{*}$-probability framework.
\end{remark}

\strut \strut

\strut

\subsection{Amalgamated R-transform Theory}

\strut

\strut

In this section, we will define an R-transform of several $B$-valued random
variables. Note that to study R-transforms is to study operator-valued
distributions. R-transforms with single variable is defined by Voiculescu
(over $B,$ in particular, $B=\Bbb{C}$. See [9] and [13]). Over $\Bbb{C},$
Nica defined multi-variable R-transforms in [1]. In [16], we extended his
concepts, over $B.$ R-transforms of $B$-valued random variables can be
defined as $B$-formal series with its $(i_{1},...,i_{n})$-th coefficients, $%
(i_{1},...,i_{n})$-th cumulants of $B$-valued random variables, where $%
(i_{1},...,i_{n})\in \{1,...,s\}^{n},$ $\forall n\in \Bbb{N}.$

\strut

\begin{definition}
Let $(A,\varphi )$ be a NCPSpace over $B$ and let $x_{1},...,x_{s}\in
(A,\varphi )$ be $B$-valued random variables ($s\in \Bbb{N}$). Let $%
z_{1},...,z_{s}$ be noncommutative indeterminants. Define a moment series of 
$x_{1},...,x_{s}$, as a $B$-formal series, by

\strut 

\begin{center}
$M_{x_{1},...,x_{s}}(z_{1},...,z_{s})=\sum_{n=1}^{\infty }\underset{%
i_{1},..,i_{n}\in \{1,...,s\}}{\sum }\varphi
(x_{i_{1}}b_{i_{2}}x_{i_{2}}...b_{i_{n}}x_{i_{n}})\,z_{i_{1}}...z_{i_{n}},$
\end{center}

\strut 

where $b_{i_{2}},...,b_{i_{n}}\in B$ are arbitrary for all $%
(i_{2},...,i_{n})\in \{1,...,s\}^{n-1},$ $\forall n\in \Bbb{N}.$

\strut 

Define an R-transform of $x_{1},...,x_{s}$, as a $B$-formal series, by

\strut 

\begin{center}
$R_{x_{1},...,x_{s}}(z_{1},...,z_{s})=\sum_{n=1}^{\infty }\underset{%
i_{1},...,i_{n}\in \{1,...,s\}}{\sum }k_{n}(x_{i_{1}},...,x_{i_{n}})%
\,z_{i_{1}}...z_{i_{n}},$
\end{center}

\strut with

\begin{center}
$k_{n}(x_{i_{1}},...,x_{i_{n}})=c^{(n)}(x_{i_{1}}\otimes
b_{i_{2}}x_{i_{2}}\otimes ...\otimes b_{i_{n}}x_{i_{n}}),$
\end{center}

\strut 

where $b_{i_{2}},...,b_{i_{n}}\in B$ are arbitrary for all $%
(i_{2},...,i_{n})\in \{1,...,s\}^{n-1},$ $\forall n\in \Bbb{N}.$ Here, $%
\widehat{c}=(c^{(n)})_{n=1}^{\infty }$ is a cumulant multiplicative function
induced by $\varphi $ in $I(A,B).$
\end{definition}

\strut

Denote a set of all $B$-formal series with $s$-noncommutative indeterminants
($s\in \Bbb{N}$), by $\Theta _{B}^{s}$. i.e if $g\in \Theta _{B}^{s},$ then

\begin{center}
$g(z_{1},...,z_{s})=\sum_{n=1}^{\infty }\underset{i_{1},...,i_{n}\in
\{1,...,s\}}{\sum }b_{i_{1},...,i_{n}}\,z_{i_{1}}...z_{i_{n}},$
\end{center}

where $b_{i_{1},...,i_{n}}\in B,$ for all $(i_{1},...,i_{n})\in
\{1,...,s\}^{n},$ $\forall n\in \Bbb{N}.$ Trivially, by definition, $%
M_{x_{1},...,x_{s}},$ $R_{x_{1},...,x_{s}}\in \Theta _{B}^{s}.$ By $\mathcal{%
R}_{B}^{s},$\ we denote a set of all R-transforms of $s$-$B$-valued random
variables. Recall that, set-theoratically,

\begin{center}
$\Theta _{B}^{s}=\mathcal{R}_{B}^{s},$ sor all $s\in \Bbb{N}.$
\end{center}

\strut

We can also define symmetric moment series and symmetric R-transform by $%
b_{0}\in B,$ by

\strut

\begin{center}
$M_{x_{1},...,x_{s}}^{symm(b_{0})}(z_{1},...,z_{s})=\sum_{n=1}^{\infty }%
\underset{i_{1},...,i_{n}\in \{1,...,s\}}{\sum }\varphi
(x_{i_{1}}b_{0}x_{i_{2}}...b_{0}x_{i_{n}})\,z_{i_{1}}...z_{i_{n}}$
\end{center}

and

\begin{center}
$R_{x_{1},...,x_{s}}^{symm(b_{0})}(z_{1},...,z_{s})=\sum_{n=1}^{\infty }%
\underset{i_{1},..,i_{n}\in \{1,...,s\}}{\sum }%
k_{n}^{symm(b_{0})}(x_{i_{1}},...,x_{i_{n}})\,z_{i_{1}}...z_{i_{n}},$
\end{center}

with

\begin{center}
$k_{n}^{symm(b_{0})}(x_{i_{1}},...,x_{i_{n}})=c^{(n)}(x_{i_{1}}\otimes
b_{0}x_{i_{2}}\otimes ...\otimes b_{0}x_{i_{n}}),$
\end{center}

for all $(i_{1},...,i_{n})\in \{1,...,s\}^{n},$ $\forall n\in \Bbb{N}.$

\strut

If $b_{0}=1_{B},$ then we have trivial moment series and trivial R-transform
of $x_{1},...,x_{s}$ denoted by $M_{x_{1},...,x_{s}}^{t}$ and $%
R_{x_{1},...,x_{s}}^{t},$ respectively.  By definition, for the fixed random
variables $x_{1},...,x_{s}\in (A,\varphi ),$ there are infinitely many
R-transforms of them (resp. moment series of them). Symmetric and trivial
R-transforms of them are special examples. Let

\strut 

\begin{center}
$C=\underset{(i_{1},...,i_{n})\in \Bbb{N}^{n}}{\cup }%
\{(1_{B},b_{i_{2}},...,b_{i_{n}}):b_{i_{j}}\in B\}.$
\end{center}

\strut 

Suppose that we have

\strut 

\begin{center}
$coef_{i_{1},...,i_{n}}\left( R_{x_{1},...,x_{s}}\right) =c^{(n)}\left(
x_{i_{1}}\otimes b_{i_{2}}x_{i_{2}}\otimes ...\otimes
b_{i_{n}}x_{i_{n}}\right) ,$
\end{center}

\strut 

where $(1_{B},b_{i_{2}},...,b_{i_{n}})\in C,$ for all $(i_{1},...,i_{n})\in 
\Bbb{N}^{n}.$ Then we can rewite the R-transform of $x_{1},...,x_{s},$ $%
R_{x_{1},...,x_{s}}$ by $R_{x_{1},...,x_{s}}^{C}.$ If $C_{1}$ and $C_{2}$
are such collections, then in general $R_{x_{1},...,x_{s}}^{C_{1}}\neq
R_{x_{1},...,x_{s}}^{C_{2}}$ (resp. $M_{x_{1},...,x_{s}}^{C_{1}}\neq
M_{x_{1},...,x_{s}}^{C_{2}}$). From now, for the random variables $%
x_{1},...,x_{s},y_{1},...,y_{s},$ if we write $R_{x_{1},...,x_{s}}$ and $%
R_{y_{1},...,y_{s}},$ then it means that $%
R_{x_{1},...,x_{s}}^{C}=R_{y_{1},....,y_{s}}^{C},$ for the same collection $%
C.$ If there's no confusion, we will omit to write such collection.

\strut

The followings are known in [23] and [16] ;

\strut

\begin{proposition}
Let $(A,\varphi )$ be a NCPSpace over $B$ and let $%
x_{1},...,x_{s},y_{1},...,y_{p}\in (A,\varphi )$ be $B$-valued random
variables, where $s,p\in \Bbb{N}.$ Suppose that $\{x_{1},...,x_{s}\}$ and $%
\{y_{1},...,y_{p}\}$ are free in $(A,\varphi ).$ Then

\strut 

(1) $
R_{x_{1},...,x_{s},y_{1},...,y_{p}}(z_{1},...,z_{s+p})=R_{x_{1},...,x_{s}}(z_{1},...,z_{s})+R_{y_{1},...,y_{p}}(z_{1},...,z_{p}).
$

\strut 

(2) If $s=p,$ then $R_{x_{1}+y_{1},...,x_{s}+y_{s}}(z_{1},...,z_{s})=\left(
R_{x_{1},...,x_{s}}+R_{y_{1},...,y_{s}}\right) (z_{1},...,z_{s}).$

$\square $
\end{proposition}

\strut

The above proposition is proved by the characterization of freeness with
respect to cumulants. i.e $\{x_{1},...,x_{s}\}$ and $\{y_{1},...,y_{p}\}$
are free in $(A,\varphi )$ if and only if their mixed cumulants vanish. Thus
we have

\strut

$k_{n}(p_{i_{1}},...,p_{i_{n}})=c^{(n)}(p_{i_{1}}\otimes
b_{i_{2}}p_{i_{2}}\otimes ...\otimes b_{i_{n}}p_{i_{n}})$

$\ \ \ \ =\left( \widehat{c_{x}}\oplus \widehat{c_{y}}\right)
^{(n)}(p_{i_{1}}\otimes b_{i_{2}}p_{i_{2}}\otimes ...\otimes
b_{i_{n}}p_{i_{n}})$

$\ \ \ \ =\left\{ 
\begin{array}{lll}
k_{n}(x_{i_{1}},...,x_{i_{n}}) &  & or \\ 
k_{n}(y_{i_{1}},...,y_{i_{n}}) &  & 
\end{array}
\right. $

\strut

and if $s=p,$ then

$k_{n}(x_{i_{1}}+y_{i_{1}},...,x_{i_{n}}+y_{i_{n}})$

$\ =c^{(n)}\left( (x_{i_{1}}+y_{i_{1}})\otimes
b_{i_{2}}(x_{i_{2}}+y_{i_{2}})\otimes ...\otimes
b_{i_{n}}(x_{i_{n}}+y_{i_{n}})\right) $

$\ \ =c^{(n)}(x_{i_{1}}\otimes b_{i_{2}}x_{i_{2}}\otimes ...\otimes
b_{i_{n}}x_{i_{n}})+c^{(n)}(y_{i_{1}}\otimes b_{i_{2}}y_{i_{2}}\otimes
...\otimes b_{i_{n}}y_{i_{n}})+[Mixed]$

\strut

where $[Mixed]$ is the sum of mixed cumulants of $x_{j}$'s and $y_{i}$'s, by
the bimodule map property of $c^{(n)}$

\strut

$\ =k_{n}(x_{i_{1}},...,x_{i_{n}})+k_{n}(y_{i_{1}},...,y_{i_{n}})+0_{B}.$

\strut

Note that if $f,g\in \Theta _{B}^{s},$ then we can always choose free $%
\{x_{1},...,x_{s}\}$ and $\{y_{1},...,y_{s}\}$ in (some) NCPSpace over $B,$ $%
(A,\varphi ),$ such that

\begin{center}
$f=R_{x_{1},...,x_{s}}$ \ \ and \ \ $g=R_{y_{1},...,y_{s}}.$
\end{center}

\strut

\begin{definition}
(1) Let $s\in \Bbb{N}.$ Let $(f,g)\in \Theta _{B}^{s}\times \Theta _{B}^{s}.$
Define \frame{*}\thinspace \thinspace $:\Theta _{B}^{s}\times \Theta
_{B}^{s}\rightarrow \Theta _{B}^{s}$ by

\strut 

\begin{center}
$\left( f,g\right) =\left(
R_{x_{1},...,x_{s}}^{C_{1}},\,R_{y_{1},...,y_{s}}^{C_{2}}\right) \longmapsto
R_{x_{1},...,x_{s}}^{C_{1}}\,\,\frame{*}\,\,R_{y_{1},...,y_{s}}^{C_{2}}.$
\end{center}

\strut 

Here, $\{x_{1},...,x_{s}\}$ and $\{y_{1},...,y_{s}\}$ are free in $%
(A,\varphi )$. Suppose that

\strut 

\begin{center}
$coef_{i_{1},..,i_{n}}\left( R_{x_{1},...,x_{s}}^{C_{1}}\right)
=c^{(n)}(x_{i_{1}}\otimes b_{i_{2}}x_{i_{2}}\otimes ...\otimes
b_{i_{n}}x_{i_{n}})$
\end{center}

and

\begin{center}
$coef_{i_{1},...,i_{n}}(R_{y_{1},...,y_{s}}^{C_{2}})=c^{(n)}(y_{i_{1}}%
\otimes b_{i_{2}}^{\prime }y_{i_{2}}\otimes ...\otimes b_{i_{n}}^{\prime
}y_{i_{n}}),$
\end{center}

\strut 

for all $(i_{1},...,i_{n})\in \{1,...,s\}^{n},$ $n\in \Bbb{N},$ where $%
b_{i_{j}},b_{i_{n}}^{\prime }\in B$ arbitrary. Then

\strut 

$coef_{i_{1},...,i_{n}}\left( R_{x_{1},...,x_{s}}^{C_{1}}\,\,\frame{*}%
\,\,R_{y_{1},...,y_{s}}^{C_{2}}\right) $

\strut 

$=\underset{\pi \in NC(n)}{\sum }\left( \widehat{c_{x}}\oplus \widehat{c_{y}}%
\right) (\pi \cup Kr(\pi ))(x_{i_{1}}\otimes y_{i_{1}}\otimes
b_{i_{2}}x_{i_{2}}\otimes b_{i_{2}}^{\prime }y_{i_{2}}\otimes ...\otimes
b_{i_{n}}x_{i_{n}}\otimes b_{i_{n}}^{\prime }y_{i_{n}})$

$\strut $

$\overset{denote}{=}\underset{\pi \in NC(n)}{\sum }\left( k_{\pi
}^{C_{1}}\oplus k_{Kr(\pi )}^{C_{2}}\right)
(x_{i_{1}},y_{i_{1}},...,x_{i_{n}}y_{i_{n}}),$

\strut \strut 

where $\widehat{c_{x}}\oplus \widehat{c_{y}}=\widehat{c}\mid
_{A_{x}*_{B}A_{y}},$ $A_{x}=A\lg \left( \{x_{i}\}_{i=1}^{s},B\right) $ and $%
A_{y}=A\lg \left( \{y_{i}\}_{i=1}^{s},B\right) $ and where $\pi \cup Kr(\pi )
$ is an alternating union of partitions in $NC(2n)$
\end{definition}

\strut

\begin{proposition}
(See [16])\strut Let $(A,\varphi )$ be a NCPSpace over $B$ and let $%
x_{1},...,x_{s},y_{1},...,y_{s}\in (A,\varphi )$ be $B$-valued random
variables ($s\in \Bbb{N}$). If $\{x_{1},...,x_{s}\}$ and $\{y_{1},...,y_{s}\}
$ are free in $(A,\varphi ),$ then we have

\strut 

$k_{n}(x_{i_{1}}y_{i_{1}},...,x_{i_{n}}y_{i_{n}})$

\strut 

$=\underset{\pi \in NC(n)}{\sum }\left( \widehat{c_{x}}\oplus \widehat{c_{y}}%
\right) (\pi \cup Kr(\pi ))(x_{i_{1}}\otimes y_{i_{1}}\otimes
b_{i_{2}}x_{i_{2}}\otimes y_{i_{2}}\otimes ...\otimes
b_{i_{n}}x_{i_{n}}\otimes y_{i_{n}})$

\strut 

$\overset{denote}{=}\underset{\pi \in NC(n)}{\sum }\left( k_{\pi }\oplus
k_{Kr(\pi )}^{symm(1_{B})}\right)
(x_{i_{1}},y_{i_{1}},...,x_{i_{n}},y_{i_{n}}),$

\strut 

for all $(i_{1},...,i_{n})\in \{1,...,s\}^{n},$ $\forall n\in \Bbb{N},$ $%
b_{i_{2}},...,b_{i_{n}}\in B,$ arbitrary, where $\widehat{c_{x}}\oplus 
\widehat{c_{y}}=\widehat{c}\mid _{A_{x}*_{B}A_{y}},$ $A_{x}=A\lg \left(
\{x_{i}\}_{i=1}^{s},B\right) $ and $A_{y}=A\lg \left(
\{y_{i}\}_{i=1}^{s},B\right) .$ \ $\square $
\end{proposition}

\strut

This shows that ;

\strut

\begin{corollary}
(See [16]) Under the same condition with the previous proposition,

\strut 

\begin{center}
$R_{x_{1},...,x_{s}}\,\,\frame{*}\,%
\,R_{y_{1},...,y_{s}}^{t}=R_{x_{1}y_{1},...,x_{s}y_{s}}.$
\end{center}

$\square $
\end{corollary}

\strut

\strut 

\strut

\subsection{$B$-valued Even Random Variables}

\strut

\strut

\strut In this section, we will consider the $B$-evenness. Let $(A,\varphi )$
be a NCPSpace over $B.$

\strut

\begin{definition}
Let $a\in (A,\varphi )$ be a $B$-valued random variable. We say that this
random variable $a$ is $B$-even if

\strut 

\begin{center}
$\varphi (ab_{2}a...b_{m}a)=0_{B},$ whenever $m$ is odd,
\end{center}

\strut 

where $b_{2},...,b_{m}\in B$ are arbitrary. In particular, if $a$ is $B$%
-even, then $\varphi (a^{m})=0_{B},$ whenever $m$ is odd. But the converse
is not true, in general.
\end{definition}

\strut

Recall that in the $*$-probability space model, the $B$-evenness guarantees
the self-adjointness (See [16]). But the above definition is more general.
By using the M\"{o}bius inversion, we have the following characterization ;

\strut

\begin{proposition}
Let $a\in (A,\varphi )$ be a $B$-valued random variable. Then $a$ is $B$%
-even if and only if

\strut 

\begin{center}
$k_{m}\left( \underset{m-times}{\underbrace{a,.......,a}}\right) =0_{B},$
whenever $m$ is odd.
\end{center}
\end{proposition}

\strut

\strut

The above proposition says that $B$-evenness is easy to veryfy when we are
dealing with either $B$-moments or $B$-cumulants. Now, define a subset $%
NC^{(even)}(2k)$ of $NC(2k),$ for any $k\in \Bbb{N}$ ;

\strut

\begin{center}
$NC^{(even)}(2k)=\{\pi \in NC(2k):\pi $ does not contain odd blocks$\}.$
\end{center}

\strut

We have that ;

\strut

\begin{proposition}
Let $k\in \Bbb{N}$ and let $a\in (A,\varphi )$ be $B$-even. Then

\strut 

\begin{center}
$k_{2k}\left( \underset{2k-times}{\underbrace{a,.......,a}}\right) =%
\underset{\pi \in NC^{(even)}(2k)}{\sum }\widehat{\varphi }(\pi )\left(
a\otimes b_{2}a\otimes ...\otimes b_{2k}a\right) \mu (\pi ,1_{2k})$
\end{center}

\strut 

equivalently,

\strut 

\begin{center}
$\varphi \left( ab_{2}a...b_{2k}a\right) =\underset{\pi \in NC^{(even)}(2k)}{%
\sum }\widehat{c}(\pi )\left( a\otimes b_{2}a\otimes ...\otimes
b_{2k}a\right) .$
\end{center}
\end{proposition}

\strut

\begin{proof}
By the previous proposition, it is enough to show one of the above two
formuli. Fix $k\in \Bbb{N}.$ Then

\strut

$\ \ k_{2k}\left( a,...,a\right) =c^{(2k)}\left( a\otimes b_{2}a\otimes
...\otimes b_{2k}a\right) $

\strut

$\ \ \ \ \ \ \ \ \ \ \ \ \ \ \ \ \ \ \ \ =\underset{\pi \in NC(2k)}{\sum }%
\widehat{\varphi }(\pi )\left( a\otimes b_{2}a\otimes ...\otimes
b_{2k}a\right) \mu (\pi .1_{2k}).$

\strut

Now, suppose that $\theta \in NC(2k)$ and $\theta $ contains its odd block $%
V_{o}\in \pi (o)\cup \pi (i).$ Then

\strut

\strut (2.2.1)

\begin{center}
$\widehat{\varphi }(\theta )\left( a\otimes b_{2}a\otimes ...\otimes
b_{2k}a\right) =0_{B}.$
\end{center}

\strut

Define

\strut

\begin{center}
$NC^{(odd)}(2k)=\{\pi \in NC(2k):\pi $ contains at least one odd block$\}.$
\end{center}

\strut

Then, for any $\theta \in NC^{(odd)}(2k),$ the formular (2.2.1) holds. So,

\strut

\begin{center}
$k_{2k}(a,...,a)=\underset{\pi \in NC(2k)\,\,\setminus \,\,NC^{(odd)}(2k)}{%
\sum }\widehat{\varphi }(\pi )(a\otimes b_{2}a\otimes ...\otimes b_{2k}a)\mu
(\pi ,1_{2k}).$
\end{center}

\strut

It is easy to see that, by definition,

\strut

\begin{center}
$NC^{(even)}(2k)=NC(2k)\,\setminus \,NC^{(odd)}(2k).$
\end{center}
\end{proof}

\strut

\begin{proposition}
Let $a_{1}$ and $a_{2}$ be $B$-even elements in $(A,\varphi ).$ If $a_{1}$
and $a_{2}$ are free over $B,$ then $a_{1}+a_{2}\in (A,\varphi )$ is $B$%
-even, again.
\end{proposition}

\strut \strut 

\strut

\strut

\section{Free Probability Theory on $L(F_{2})*_{L(F_{1})}L(F_{2})$ over $%
L(F_{1})$}

\strut \strut

\strut

In this chapter, we will consider the free group $W^{*}$-algebras $%
B=L(F_{1}) $ and $A=L(F_{2}),$ where $F_{N}$ is a free group with $N$%
-generators ($N\in \Bbb{N}$). i.e

\strut

\begin{center}
$B=\{\underset{h\in F_{1}}{\sum }t_{h}h:\,t_{h}\in \Bbb{C}\}$
\end{center}

and

\begin{center}
$A=\{\underset{g\in F_{2}}{\sum }t_{g}g:t_{g}\in \Bbb{C}\}.$
\end{center}

Recall that there is a map $E:A\rightarrow B$ defined by

\strut

\begin{center}
$E\left( \underset{g\in F_{2}}{\sum }t_{g}g\right) =\underset{h\in F_{1}}{%
\sum }t_{h}h.$
\end{center}

\strut

Notice that $E:A\rightarrow B$ is a conditional expectation ($B$-functional)
and hence $(A,E)$ is a NCPSpace over $B.$

\strut

Now, for any $N\in \Bbb{N},$\ define the canonical trace $\varphi $ on $%
L(F_{N})$ ;

\strut

\begin{center}
$\varphi \left( \underset{g\in F_{N}}{\sum }t_{g}g\right) =t_{e_{F_{N}}},$
\end{center}

\strut

for all $\underset{g\in F_{N}}{\sum }t_{g}g\in L(F_{N}),$ where $e_{F_{N}}$
is the identity of $F_{N}.$ For the convenience of using notation, we will
denote $e_{F_{N}}$ by $e.$

\strut

In this paper, we will concentrate on finding scalar-valued moments of $%
a+b+a^{-1}+b^{-1}+c+d+c^{-1}+d^{-1}\in L(F_{2})*_{L(F_{1})}L(F_{2})$,

\strut

\begin{center}
$\tau \left( (a+b+a^{-1}+b^{-1}+c+d+c^{-1}+d^{-1})^{n}\right) ,$
\end{center}

\strut

where $L\left( <a,b>\right) =L(F_{2})=L(<c,d>)$ and $\tau
:L(F_{2})*_{L(F_{1})}L(F_{2})\rightarrow \Bbb{C}$ is the trace$,$ for all $%
n\in \Bbb{N},$ defined by

\strut

\begin{center}
$\tau \left( y_{1}...y_{n}\right) =t_{e},$ for all $y_{j}\in L(F_{2}),$ $%
j=1,...,n.$
\end{center}

\strut

Remark that, since each $y_{j}$ has the form, $y_{j}=\underset{g\in F_{2}}{%
\sum }t_{g}^{(j)}g,$ in $L(F_{2}),$ we can find the coefficien of $%
e=1_{B}=1_{A}$ in $y_{1}...y_{n}\in L(F_{2})*_{L(F_{1})}L(F_{2}).$ \ So, to
find moments of an element in $L(F_{2})*_{L(F_{1})}L(F_{2})$ is to find $e$%
-terms of the element in $L(F_{2})*_{L(F_{1})}L(F_{2}).$ To directely
compute this moments is very complicated. So, later, we will use the
compatibility and $B$-freeness. Also, later, we will denote this linear
functional $\tau $ by $\varphi ,$ because

\strut

\begin{center}
$\tau (x)=\varphi \left( E*E(x)\right) ,$ for all $x\in
L(F_{2})*_{L(F_{1})}L(F_{2}).$
\end{center}

\strut \strut

\strut

\section{Compatibility of $\left( L(F_{2})*_{L(F_{1})}L(F_{2}),\varphi
\right) $ and $\left( L(F_{2})*_{L(F_{1})}L(F_{2}),F\right) $}

\strut

\strut

\strut

From now on, by $A$ and $B,$ we will denote $L(F_{2})$ and $L(F_{1}),$ $A$
and $B,$ respectively. By the very definitions of $E:A\rightarrow B$ and $%
\tau :A\rightarrow \Bbb{C},$

\strut

\begin{center}
$E\left( \underset{g\in F_{2}}{\sum }t_{g}g\right) =\underset{h\in F_{1}}{%
\sum }t_{h}h$ \ and \ $\tau \left( \underset{g\in F_{2}}{\sum }t_{g}g\right)
=t_{e},$
\end{center}

\strut \strut

a NCPSpace $(A,\varphi )$ and an amalgamated NCPSpace over $B,$ $(A,E)$ are
compatible. In this section, we will show that

\strut

\begin{center}
$\varphi (x)=\varphi \circ E(x),$ for all $x\in A.$
\end{center}

\strut \strut \strut

We can regard $e$ as the identity element in $B$ and $A.$ i.e

\strut

\begin{center}
$1_{A}=e=1_{B}.$
\end{center}

\strut \strut \strut \strut

\begin{lemma}
Let $B\subset A,$ $\varphi $ and $E$ be given as before. Then a NCPSpace $%
(A,\varphi )$ and a NCPSpace over $B,$ $(A,E)$ are compatible.
\end{lemma}

\strut

\strut Note that the trace $\varphi $ on $A_{1}*_{B}A_{2}$ and $\varphi
\circ (E*E)$ coincide. So we have that ;

\strut

\begin{theorem}
$(A_{1}*_{B}A_{2},\varphi )$ and $(A_{1}*_{B}A_{2},\,E*E)$ are compatible.
\end{theorem}

\strut 

\strut 

\subsection{$B$-Evenness and $B$-identically distributedness of $x$ and $y$
in $(A_{1}*_{B}A_{2},E*E)$}

\strut

\strut

By $F,$ we will denote the $B$-functional $E*E:A_{1}*_{B}A_{2}\rightarrow B.$
And in the rest of this paper, we will let

\strut

\begin{center}
$x=a+b+a^{-1}+b^{-1}$ \ and \ $y=c+d+c^{-1}+d^{-1},$
\end{center}

\strut

in $A_{1}*_{B}A_{2}.$

\strut

\begin{lemma}
As a $B$-valued random variable, $x=a+b+a^{-1}+b^{-1}\in (A_{1}*_{B}A_{2},F)$
is $B$-even. $\square $
\end{lemma}

\strut

The above lemma is proved by a straightforward observation. Next section, we
will observe the $B$-evenness of $x,$ in detail.

\strut \strut

\begin{proposition}
Let $x$ and $y$ be given as before, as $B$-valued random variables in $%
\left( A_{1}*_{B}A_{2},F\right) .$ Then $\{x\}$ and $\{y\}$ are free over $B,
$ in $(A_{1}*_{B}A_{2},F)$ and they are identically distributed. i.e

\strut 

\begin{center}
$R_{x}=R_{y}$ \ \ or \ \ $M_{x}=M_{y}.$
\end{center}
\end{proposition}

\strut

\begin{proof}
Since $x\in A_{1}$ and $y\in A_{2}$ in $A_{1}*_{B}A_{2},$ they are free over 
$B,$ in $\left( A_{1}*_{B}A_{2},F\right) .$ By the generating property of $%
\{a,b\}$ and $\{c,d\}$ (i.e they generate same group $F_{2}$), they are
identically distributed. Equivalently,

\strut

\begin{center}
$R_{x}(z)=R_{y}(z)$
\end{center}

$\Longleftrightarrow $

\begin{center}
$
\begin{array}{ll}
K_{n}\left( \underset{n-times}{\underbrace{x\otimes .....\otimes x}}\right)
& =C^{(n)}\left( x\otimes b_{2}x\otimes ...\otimes b_{n}x\right) \\ 
& =C^{(n)}(y\otimes b_{2}y\otimes ...\otimes b_{n}y) \\ 
& =K_{n}\left( \underset{n-times}{\underbrace{y,.....,y}}\right) ,
\end{array}
$
\end{center}

\strut

for all $n\in \Bbb{N}.$ By the M\"{o}bius inversion,

\strut

\begin{center}
$M_{x}(z)=M_{y}(z),$
\end{center}

\strut

as $B$-formal series.
\end{proof}

\strut

\begin{corollary}
Let $x$ and $y$ be given as before, in $(A_{1}*_{B}A_{2},\,F).$ Then

\strut 

\begin{center}
$R_{x+y}(z)=\left( R_{x}+R_{y}\right) (z)=2R_{x}(z).$
\end{center}

$\square $
\end{corollary}

\strut \strut

\begin{corollary}
Let $x$ and $y$ be given as before. Then $x+y$ is $B$-even, too.
\end{corollary}

\strut

\begin{proof}
By the previous lemma, $x$ is $B$-even. Since $y$ is identically distributed
with $x,$ their R-transforms are same and hence $y$ is $B$-even, too. In
[16], we showed that if two $B$-even $B$-valued random variables are $B$%
-free, then the sum of them is also $B$-even. Since our $B$-valued random
variables are $B$-even, $x+y$ is also $B$-even.
\end{proof}

\strut

\strut

\strut

\subsection{\strut Computation of $B$-valued moments of $x,$ $E(x^{n})$}

\strut \strut

\strut

To compute $E(x^{n}),$ we will use some results in [15]. We have that $x\in
\left( A_{1}*_{B}A_{2},F\right) $ is $B$-even. Thus

\strut

\begin{center}
$F(x^{n})=0_{B},$ whenever $n$ is odd.
\end{center}

\strut

So, we concentrate on finding $B$-valued $2n$-th moments of $%
x=a+b+a^{-1}+b^{-1},$ in $(A_{1}*_{B}A_{2},\,F).$ It is known that if we
denote

\strut

\begin{center}
$X_{n}=\underset{\left| w\right| =n}{\sum }w\in \Bbb{C}[F_{N}]$ \ ($N\in 
\Bbb{N}$),
\end{center}

then

\begin{center}
$X_{1}X_{1}=X_{2}+2N\cdot e$
\end{center}

and

\begin{center}
$X_{1}X_{n}=X_{n+1}+(2N-1)X_{n-1},$
\end{center}

\strut

where $e=e_{F_{N}},$ for all $n\in \Bbb{N}$ $\setminus \,\{1\}.$

\strut

In our case, we can regard our $x=a+b+a^{-1}+b^{-1}$ as $X_{1}$ in $\Bbb{C}%
[F_{2}]=A_{1}.$ Thus we have that

\strut

\begin{center}
$X_{1}X_{1}=X_{2}+4e$
\end{center}

and

\begin{center}
$X_{1}X_{n}=X_{n+1}+3X_{n-1}$ \ ($n=2,3,...$).
\end{center}

\strut

By using those two results, we can express $x^{n}$ in terms of $X_{k}$'s ;
For example,

\strut

$x^{2}=x\cdot x=X_{1}X_{1}=X_{2}+4e,$

\strut

$x^{3}=x\cdot x^{2}=X_{1}\left( X_{2}+4e\right)
=X_{1}X_{2}+4X_{1}=X_{3}+3X_{1}+4X_{1}=X_{3}+(3+4)X_{1},$

\strut

continuing

\strut

$x^{4}=X_{4}+(3+3+4)X_{2}+4(3+4)e,$

$x^{5}=X_{5}+(3+3+3+4)X_{3}+\left( 3(3+3+4)+4(3+4)\right) X_{1},$

$x^{6}=X_{6}+\left( 3+3+3+3+4\right) X_{4}$

$\ \ \ \ \ \ \ \ \ \ \ \ \ +\left( 3(3+3+3+4)+3(3+3+4)+4(3+4)\right) X_{2}$

$\ \ \ \ \ \ \ \ \ \ \ \ \ +4\left( 3(3+3+4)+4(3+4)\right) e,$

etc.

\strut

So, we can find a recurrence relation to get $x^{n}$ ($n\in \Bbb{N}$) with
respect to $X_{k}$'s ($k\leq n$). Inductively, we have that $x^{2k-1}$ and $%
x^{2k}$ have their representations in terms of $X_{j}$'s as follows ;

\strut

\begin{center}
$%
x^{2k-1}=X_{1}^{2k-1}=X_{2k-1}+q_{2k-3}^{2k-1}X_{2k-3}+q_{2k-5}^{2k-1}X_{2k-5}+...+q_{3}^{2k-1}X_{3}+q_{1}^{2k-1}X_{1} 
$
\end{center}

\strut \strut

and

\begin{center}
$%
x^{2k}=X_{1}^{2k}=X_{2k}+p_{2k-2}^{2k}X_{2k-2}+p_{2k-4}^{2k}X_{2k-4}+...+p_{2}^{2k}X_{2}+p_{0}^{2k}e, 
$
\end{center}

\strut

where $k\geq 2.$ Also, we have the following recurrence relation ;

\strut

\begin{proposition}
Let's fix $k\in \Bbb{N}\,\setminus \,\{1\}.$ Let $q_{i}^{2k-1}$ and $%
p_{j}^{2k}$ ($i=1,3,5,...,2k-1,....$ and $j=0,2,4,...,2k,...$) be given as
before. If $p_{0}^{2}=4$ and $q_{1}^{3}=3+p_{0}^{2},$ then we have the
following recurrence relations ;

\strut 

(1) Let

\begin{center}
$%
x^{2k-1}=X_{2k-1}+q_{2k-3}^{2k-1}X_{2k-3}+...+q_{3}^{2k-1}X_{3}+q_{1}^{2k-1}X_{1}.
$
\end{center}

Then

\strut 

$x^{2k}=X_{2k}+\left( 3+q_{2k-3}^{2k-1}\right) X_{2k-2}+\left(
3q_{2k-3}^{2k-1}+q_{2k-5}^{2k-1}\right) X_{2k-4}$

$\ \ \ \ \ \ \ \ \ \ \ \ \ \ \ +\left(
3q_{2k-5}^{2k-1}+q_{2k-7}^{2k-1}\right) X_{2k-6}+$

$\ \ \ \ \ \ \ \ \ \ \ \ \ \ \ +...+\left( 3q_{3}^{2k-1}+q_{1}^{2k-1}\right)
X_{2}+4q_{1}^{2k-1}e.$

\strut i.e,

\strut 

$p_{2k-2}^{2k}=3+q_{2k-3}^{2k-1},$ $\
p_{2k-4}^{2k}=3q_{2k-3}^{2k-1}+q_{2k-5}^{2k-1},...,$ $%
p_{2}^{2k}=3q_{3}^{2k-1}+q_{1}^{2k-1}$ and $p_{0}^{2k}=4q_{1}^{2k-1}.$

\strut 

(2) Let

\begin{center}
$x^{2k}=X_{2k}+p_{2k-2}^{2k}X_{2k-2}+...+p_{2}^{2k}X_{2}+p_{0}^{2k}e.$
\end{center}

Then

\strut 

$x^{2k+1}=X_{2k+1}+\left( 3+p_{2k-2}^{2k}\right) X_{2k-1}+\left(
3p_{2k-2}^{2k}+p_{2k-4}^{2k}\right) X_{2k-3}$

$\ \ \ \ \ \ \ \ \ \ \ \ \ \ \ \ \ \ \ \ \ +\left(
3p_{2k-4}^{2k}+p_{2k-6}^{2k}\right) X_{2k-5}+$

$\ \ \ \ \ \ \ \ \ \ \ \ \ \ \ \ \ \ \ \ \ +...+\left(
3p_{4}^{2k}+p_{2}^{2k}\right) X_{3}+\left( 3p_{2}^{2k}+p_{0}^{2k}\right)
X_{1}.$

i.e,

\strut 

$q_{2k-1}^{2k+1}=3+p_{2k-2}^{2k},$ \ $%
q_{2k-3}^{2k+1}=3p_{2k-2}^{2k}+p_{2k-4}^{2k},...,$ \ $%
q_{3}^{2k+1}=3p_{4}^{2k}+p_{2}^{2k}$ and $%
q_{1}^{2k+1}=3p_{2}^{2k}+p_{0}^{2k}.$ \ \ \ $\square $
\end{proposition}

\strut

The above recurrence relations give us an algorithm, how to find the $2n$-th 
$B$-valued moments of $x^{2n}\in (A_{1}*_{B}A_{2},\,F).$

\strut

\begin{example}
Let $p_{0}^{2}=4$ and $q_{1}^{3}=3+p_{0}^{2}=3+4=7.$ Put

\strut 

\begin{center}
$x^{8}=X_{8}+p_{6}^{8}X_{6}+p_{4}^{8}X_{4}+p_{2}^{8}X_{4}+p_{0}^{8}e.$
\end{center}

\strut 

Then, by the previous proposition, we have that

\strut 

\begin{center}
$p_{6}^{8}=3+q_{5}^{7},$ \ \ $p_{4}^{8}=3q_{5}^{7}+q_{3}^{7},$ \ \ $%
p_{2}^{8}=3q_{3}^{7}+q_{1}^{7}$ \ \ and \ $p_{0}^{8}=4q_{1}^{7}.$
\end{center}

\strut 

Similarly, by the previous proposition,

\strut 

\ \ \ \ $\ q_{5}^{7}=3+p_{4}^{6},$ \ \ \ $q_{3}^{7}=3p_{4}^{6}+p_{2}^{6}$ \
\ \ \ and \ \ $\ q_{1}^{7}=3p_{2}^{6}+p_{0}^{6},$

$\ \ \ \ \ \ p_{4}^{6}=3+q_{3}^{5},$ \ \ \ $p_{2}^{6}=3q_{3}^{5}+q_{1}^{5}$
\ \ \ \ \ and \ \ \ \ $p_{0}^{6}=4q_{1}^{5},$

\strut 

$\ \ \ \ \ \ q_{3}^{5}=3+p_{2}^{4}$ \ \ \ \ \ \ \ and \ \ \ \ \ \ $%
q_{1}^{5}=3p_{2}^{4}+p_{2}^{4},$

\strut 

$\ \ \ \ \ \ p_{2}^{4}=3+q_{1}^{3}$ \ \ \ \ \ \ \ \ \ and \ \ \ \ \ \ \ \ \ $%
p_{0}^{4}=4q_{1}^{3},$

\strut 

and

\strut 

$\ \ \ \ \ \ \ q_{1}^{3}=3+p_{0}^{2}=7.$

\strut 

Therefore, combining all information,

\strut 

\begin{center}
$x^{8}=X_{8}+22\,X_{6}+202\,X_{4}+744\,X_{2}+1316\,e.$
\end{center}
\end{example}

\strut

We have the following diagram with arrows which mean that

\begin{center}
$\swarrow \swarrow $ \ : \ $3+[$former term$]$

$\searrow $ \ \ \ \ : \ \ $3\cdot [$former term$]$

$\swarrow $ \ \ \ \ : \ \ $\cdot +[$former term$]$
\end{center}

and

\begin{center}
$\searrow \searrow $ \ : \ \ $4\cdot [$former term$].$
\end{center}

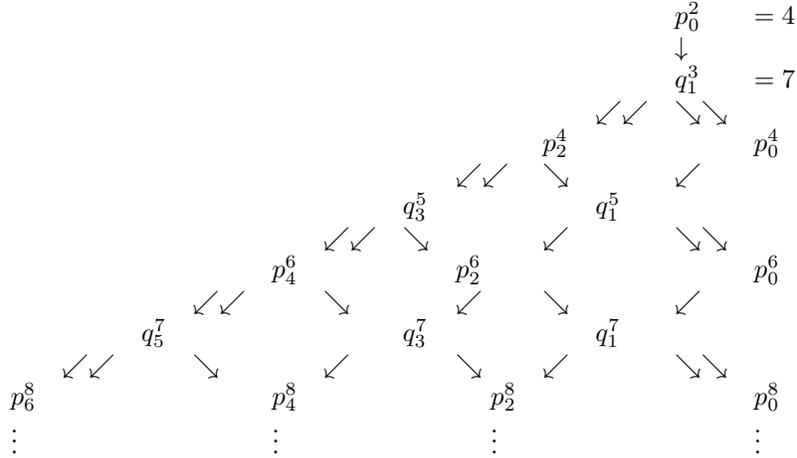
\strut

\begin{center}
$
\begin{array}{llllllllllll}
&  &  &  &  &  &  &  &  &  & p_{0}^{2} & =4 \\ 
&  &  &  &  &  &  &  &  &  & \downarrow &  \\ 
&  &  &  &  &  &  &  &  &  & q_{1}^{3} & =7 \\ 
&  &  &  &  &  &  &  &  & \swarrow \swarrow & \searrow \searrow &  \\ 
&  &  &  &  &  &  &  & p_{2}^{4} &  &  & p_{0}^{4} \\ 
&  &  &  &  &  &  & \swarrow \swarrow & \searrow &  & \swarrow &  \\ 
&  &  &  &  &  & q_{3}^{5} &  &  & q_{1}^{5} &  &  \\ 
&  &  &  &  & \swarrow \swarrow & \searrow &  & \swarrow &  & \searrow
\searrow &  \\ 
&  &  &  & p_{4}^{6} &  &  & p_{2}^{6} &  &  &  & p_{0}^{6} \\ 
&  &  & \swarrow \swarrow &  & \searrow &  & \swarrow & \searrow &  & 
\swarrow &  \\ 
&  & q_{5}^{7} &  &  &  & q_{3}^{7} &  &  & q_{1}^{7} &  &  \\ 
& \swarrow \swarrow &  & \searrow &  & \swarrow &  & \searrow & \swarrow & 
& \searrow \searrow &  \\ 
p_{6}^{8} &  &  &  & p_{4}^{8} &  &  & \text{ \ \ \ }p_{2}^{8} &  &  &  & 
p_{0}^{8} \\ 
\vdots &  &  &  & \vdots &  &  & \text{ \ \ \ }\vdots &  &  &  & \vdots
\end{array}
$
\end{center}

\strut

\strut

Now, recall that $h=aba^{-1}b^{-1}$ (or $h=cdc^{-1}d^{-1}$) in $\left(
A_{1}*_{B}A_{2},\,F\right) $ ,where $F_{1}=\,<h>.$ So, since $F_{1}=\,<h>$
is a cyclic group, WLOG, we denote $\underset{g\in F_{1}}{\sum }t_{g}g\in
B\hookrightarrow A_{1}*_{B}A_{2}$ \ by \ $\sum_{n=-\infty }^{\infty
}t_{n}h^{n}\in B,$ where $t_{g},$ $t_{n}\in \Bbb{C},$ with $t_{0}=t_{e}.$
Therefore, we can let

\strut

\begin{center}
$\varphi \left( \sum_{n=-\infty }^{\infty }t_{n}h^{n}\right) =t_{0}\equiv
t_{e}$
\end{center}

\strut

\strut And hence, to find a $B$-value moment of $x$ is to find $h^{n}$%
-terms, where $n\in \Bbb{Z}$.Note that $h$ and $h^{-1}$ are words with their
length 4. Therefore, $X_{4k}$ contains $h^{k}$-terms and $h^{-k}$-terms, for 
$k\in \Bbb{N}$ !

\strut \strut \strut

\begin{theorem}
Fix $k\in \Bbb{N}.$ Let $h=aba^{-1}b^{-1}\in A_{1}*_{B}A_{2}$ with $h^{0}=e.$

\strut 

(1) $E\left( x^{4k}\right) =\left( h^{k}+h^{-k}\right)
+\sum_{j=1}^{k-1}p_{4k-4j}^{4k}\left( h^{k-j}+h^{-(k-j)}\right)
+p_{0}^{4k}h^{0},$

\strut 

where $p_{0}^{4}=28.$

\strut 

(2) If $4\nmid 2k$ and if there are $X_{4l_{1}},...,X_{4l_{p}}$ terms in $%
x^{2k},$ then

\strut 

\begin{center}
$E(x^{2k})=\sum_{j=1}^{k-1}p_{(2k-2)-4j}^{2k}\left( h^{\frac{k-1}{2}%
-2j}+h^{-(\frac{k-1}{2}-2j)}\right) +p_{0}^{2k}h^{0},$

\strut 
\end{center}

where $p_{0}^{2}=4.$
\end{theorem}

\strut

\begin{proof}
(1) By the straightforward computation using the previous proposition, we
have that

\strut

$E\left( x^{4k}\right) $

\strut

$\ =E\left(
X_{4k}+p_{4k-2}^{4k}X_{4k-2}+p_{4k-4}^{4k}X_{4k-4}+...+p_{4}^{4k}X_{4}+p_{2}^{4k}X_{2}+p_{0}^{4k}h^{0}\right) 
$

\strut

$\ =E(X_{4k})+p_{4k-2}^{4k}E(X_{4k-2})+p_{4k-4}^{4k}E(X_{4k-4})+$

\begin{center}
$...+p_{4}^{4k}E(X_{4})+p_{2}^{4k}E(X_{2})+p_{0}^{4k}h^{0}.$
\end{center}

\strut

Since $h^{p}$ and $h^{-p}$ terms are in $X_{4p},$ for any $p\in \Bbb{N},$ we
can continue the above computation ;

\strut

$\
=E(X_{4k})+p_{4k-4}^{4k}E(X_{4k-4})+...+p_{4}^{4k}E(X_{4})+p_{0}^{4k}h^{0} $

\strut

$\ =\left( h^{k}+h^{-k}\right)
+p_{4k-4}^{4k}(h^{k-1}+h^{-(k-1)})+...+p_{4}^{4k}(h+h^{-1})+p_{0}^{4k}h^{0}.$

\strut

(2) If $4\nmid 2k,$ then $k=1,3,5,....$. If $k=1,$ then the above formula
holds true ;

\strut

\begin{center}
$E(x^{2})=E\left( X_{2}+4h^{0}\right) =4h^{0}.$
\end{center}

\strut

If $k\neq 1$ is odd, then

\strut

$E(x^{2k})$

\strut

$\
=E(X_{2k}+p_{2k-2}^{2k}X_{2k-2}+p_{2k-4}^{2k}X_{2k-4}+p_{2k-6}^{2k}X_{2k-6}+$

$\ \ \ \ \ \ \ \ \ \ \ \ \ \ \ \ \ \ \ \ \ \ \ \ \ \ \ \ \ \ \ \ \ \ \ \ \
...+p_{4}^{2k}X_{4}+p_{2}^{2k}X_{2}+p_{0}^{2k}h^{0})$

\strut

$\
=E(X_{2k})+p_{2k-2}^{2k}E(X_{2k-2})+p_{2k-4}^{2k}E(X_{2k-4})+p_{2k-6}^{2k}E(X_{2k-6})+ 
$

$\ \ \ \ \ \ \ \ \ \ \ \ \ \ \ \ \ \ \ \ \ \ \ \ \ \ \ \ \ \ \ \ \ \ \ \ \
...+p_{4}^{2k}E(X_{4})+p_{2}^{2k}E(X_{2})+p_{0}^{2k}h^{0}$

\strut

$\ =0_{B}+p_{2k-2}^{2k}\left( h^{k-1}+h^{-(k-1)}\right)
+0_{B}+p_{2k-6}^{2k}\left( h^{k-3}+h^{-(k-3)}\right) +$

$\ \ \ \ \ \ \ \ \ \ \ \ \ \ \ \ \ \ \ \ \ \ \ \ \ \ \ \ \ \ \ \ \ \ \ \ \
...+p_{4}^{2k}(h+h^{-1})+0_{B}+p_{0}^{2k}h^{0},$

\strut

since $X_{2k-2},$ $X_{2k-6},...,X_{4}$ contain $h^{p}$-terms and $h^{-p}$%
-terms.
\end{proof}

\strut

\strut

\strut

\strut

\strut

\strut

\section{The Amalgamated R-transform of $x+y$}

\strut

\strut

\strut

Throughout this section, we will use the same notations used in the previous
sections. To compute $F\left( (x+y)^{n}\right) ,$ we will consider the $%
R_{x+y}^{t}.$ Since $x$ and $y$ are free over $B,$ we have that

\strut

\begin{center}
$R_{x+y}^{t}=R_{x}^{t}+R_{y}^{t}.$
\end{center}

\strut

And since $x$ and $y$ are identically distributed,

\strut

\begin{center}
$R_{x+y}^{t}=R_{x}^{t}+R_{y}^{t}=2R_{x}^{t}$ \ \ or \ \ $2R_{y}^{t}.$
\end{center}

\strut

The above paragraph shows that why we need to observe $R_{x+y}^{t},$ to get
a $n$-th coefficients of $M_{x+y}.$ By the $B$-freeness of $x$ and $y,$ we
can compute $n$-th coefficients of $R_{x+y}^{t}=2R_{x}^{t},$ relatively
easier than to compute $n$-th coefficients of $M_{x+y}^{t},$ directly.
Moreover, since we have the recurrence relation for $F(x^{n})=E(x^{n}),$ $%
n\in \Bbb{N},$ we have a tool for computing $n$-th coefficients of $%
2R_{x}^{t}.$ But there is a difficulty to compute them, because of the
insertion property of noncrossing partitions (different from the
scalar-valued case). Hence we need to find other recurrence relation related
to this insertion property. By the $B$-evenness of $x\in \left(
A_{1}*_{B}A_{2},F\right) ,$ we have that

\strut

\begin{center}
$coef_{n}\left( R_{x}^{t}\right) =k_{n}^{t}\left( \underset{n-times}{%
\underbrace{x,......,x}}\right) =0_{B},$ whenever $n$ is odd.
\end{center}

\strut

So, we will only consider the even coefficient of $R_{x}^{t}.$

\strut

\begin{lemma}
For $k\in \Bbb{N},$

\strut 

$coef_{2k}\left( R_{x}^{t}\right) =\underset{l_{1},...,l_{p}\in
\{2,4,...,2k\},\,\,l_{1}+l_{2}+..+l_{p}=2k}{\sum }\,$

\begin{center}
$\underset{\theta \in NC_{l_{1},...,l_{p}}(2k)}{\sum }\widehat{F}(\theta
)\left( x\otimes ...\otimes x\right) \mu (\theta ,1_{2k}).$
\end{center}
\end{lemma}

\strut

\strut \strut

The above lemma shows that we need to construct a recurrence relation for

\strut

\begin{center}
$F\left(
x^{m_{1}}{}^{k_{1}}x^{m_{2}}h^{k_{2}}x^{m_{3}}...h^{k_{n}}x^{m_{n}}\right)
=E\left(
x^{m_{1}}h^{k_{1}}x^{m_{2}}h^{k_{2}}x^{m_{3}}...h^{k_{n}}x^{m_{n}}\right) ,$
\end{center}

\strut

where $m_{j}\in \Bbb{N}$ and $k_{j}\in \Bbb{Z},$ $j=1,...,n,$ for all $n\in 
\Bbb{N}.$ This recurrence relation can explain the computation of
partition-dependent $B$-valued moments with respect to $E.$ By the
observation of Section 3.2.2, it suffices to find the recurrence relation for

\strut

\begin{center}
$F\left(
X_{q_{1}}h^{k_{1}}X_{q_{2}}h^{k_{2}}X_{q_{3}}...h^{k_{n}}X_{q_{n}}\right)
=E\left(
X_{q_{1}}h^{k_{1}}X_{q_{2}}h^{k_{2}}X_{q_{3}}...h^{k_{n}}X_{q_{n}}\right) ,$
\end{center}

\strut

where $q_{1},...,q_{n}\in \Bbb{N}$ and $X_{N}=\underset{\left| w\right| =N}{%
\sum }w,$ where $w$ is a word in $\{a,b,a^{-1},b^{-1}\},$ for $N\in \Bbb{N}.$

\strut

\strut

\strut

\subsection{Recurrence Relation For $E\left( X_{m}h^{k}X_{n}\right) $}

\strut

\strut

In this section, we will consider the recurrence relation for $%
E(X_{m}hX_{n}).$ Then we can generalize this case to $E(X_{m}h^{k}X_{n}),$
where $k\in \Bbb{Z}$ and $m,n\in \Bbb{N}.$ We have that

\strut

\begin{center}
$E(X_{1}hX_{3})=e=E(X_{3}hX_{1})$ \ \ \ and\ \ \ \ \ $E(X_{2}hX_{2})=e$
\end{center}

\strut

and

\strut

\begin{proposition}
Let $m,n\in \Bbb{N}$ and $k\in \Bbb{Z}.$ Then

\strut 

(1)

\begin{center}
$E\left( h^{k}X_{n}\right) =\left\{ 
\begin{array}{lll}
h^{k}h^{\frac{n}{4}}=h^{k+\frac{n}{4}} &  & \text{if \ }4\mid n \\ 
&  &  \\ 
0_{B} &  & \text{otherwise.}
\end{array}
\right. $
\end{center}

\strut 

(2)

\begin{center}
$E(X_{m}h^{k})=\left\{ 
\begin{array}{lll}
h^{\frac{m}{4}}h^{k}=h^{\frac{m}{4}+k} &  & \text{if \ }4\mid m \\ 
&  &  \\ 
0_{B} &  & \text{otherwise.}
\end{array}
\right. $
\end{center}
\end{proposition}

\strut

\begin{proof}
Since $E:A_{j}\rightarrow B$ is a conditional expectation ($j=1,2$), we have
that

\strut \strut

(1)

\begin{center}
$E(h^{k}X_{n})=h^{k}E(X_{n})=\left\{ 
\begin{array}{lll}
h^{k}h^{\frac{n}{4}}=h^{k+\frac{n}{4}} &  & \text{if \ }4\mid n \\ 
&  &  \\ 
0_{B} &  & \text{otherwise.}
\end{array}
\right. $
\end{center}

\strut

(2)

\begin{center}
$E(X_{m}h^{k})=E(X_{m})h^{k}=\left\{ 
\begin{array}{lll}
h^{\frac{m}{4}}h^{k}=h^{\frac{m}{4}+k} &  & \text{if \ }4\mid m \\ 
&  &  \\ 
0_{B} &  & \text{otherwise.}
\end{array}
\right. $
\end{center}
\end{proof}

\strut \strut

Now fix the sufficiently big numbers $m\,$\ and $n$ in $\Bbb{N}.$ Then we
can have that

\strut (3.3.1.1)

\begin{center}
$X_{m}hX_{n}=\left( \underset{\left| w\right| =m}{\sum }w\right) h\left( 
\underset{\left| w^{\prime }\right| =n}{\sum }w^{\prime }\right) ,$
\end{center}

\strut

where $w$ and $w^{\prime }$ are words in $\{a,b,a^{-1},b^{-1}\}.$ Recall
that $h=aba^{-1}b^{-1}$ is a word with length 4. Hence, by the possible
cancellation, we can rewrite (3.3.1.1) as

\strut

$X_{m}hX_{n}$

$\strut \strut $

$\ \ \ =W_{m+n+4}+W_{m+n+2}+W_{m+n}+W_{m+n-2}$

$\strut \strut $

$\ \ \ \ \ \ \ \ \ \ \ \ \ \ \ +\left( \underset{\left| w\right|
=m-4,\,\,End(w)\neq b}{\sum }w\right) \left( X_{n}\right) +X_{m}\left( 
\underset{\left| w\right| =n-4,\,\,Init(w)\neq a}{\sum }w\right) $ \ \ \ \
(3.3.1.2)

\strut

$\ \ \ \ \ \ \ \ \ \ \ \ \ \ \ +\left( \underset{\left| w\right|
=m-3,\,\,End(w)\neq a^{-1}}{\sum }w\right) \left( \underset{\left| w\right|
=n-1,\,\,Init(w)\neq b^{-1}}{\sum }w\right) $ \ \ \ \ \ \ \ \ \ \ \ \ \ \ \
\ \ (3.3.1.3)

$\ \ \ \ \ \ \ \ \ \ \ \ \ \ \ +\left( \underset{\left| w\right|
=m-1,\,\,End(w)\neq a}{\sum }w\right) \left( \underset{\left| w\right|
=n-3,\,\,Init(w)\neq b}{\sum }w\right) $ \ \ \ \ \ \ \ \ \ \ \ \ \ \ \ \ \ \
\ \ \ \ (3.3.1.4)

$\ \ \ \ \ \ \ \ \ \ \ \ \ \ \ +\left( \underset{\left| w\right|
=m-2,\,\,End(w)\neq a}{\sum }w\right) \left( \underset{\left| w\right|
=n-2,\,\,Init(w)\neq a^{-1}}{\sum }w\right) ,$ \ \ \ \ \ \ \ \ \ \ \ \ \ \ \
\ \ \ \ (3.3.1.5)

\strut

where $W_{m+n+k}$ is the subsum of words with length $m+n+k.$ In the above
formula, (3.3.1.2) is gotten from the full cancellation of $X_{m}$ and $h,$
and the full cancellation of $h$ and $X_{n}.$ (3.3.1.3) (resp. (3.3.1.4)) is
gotten from the $3$-letter-cancellation of $X_{m}$ and $h$ (resp. $1$%
-letter-cancellation of $X_{m}$ and $h$) and the 1-letter-cancellation of $%
X_{n}$ and $h$ (resp. 3-letter-cancellation of $X_{n}$ and $h$). Similarly,
(3.3.1.5) is gotten from the 2-letter-cancellation from the left and right
of $h.$ Similar to the full $h$-cancellation, (3.3.1.2) $\sim $ (3.3.1.5),
we can rewrite that

\strut

(3.3.1.6)

\strut

$\ \ W_{m+n+2}=\left( \underset{\left| w\right| =m-1,\,\,End(w)\neq a}{\sum }%
w\right) (ba^{-1}b^{-1})X_{n}+X_{m}(aba^{-1})\left( \underset{\left|
w\right| =n-1,\,\,Init(w)\neq b^{-1}}{\sum }w\right) ,$

\strut \strut

(3.3.1.7)

\strut

$\ \ W_{m+n}=\left( \underset{\left| w\right| =m-2,\,\,End(w)\neq b}{\sum }%
w\right) (a^{-1}b^{-1})X_{n}+X_{m}(ab)\left( \underset{\left| w\right|
=n-2,\,\,Init(w)\neq a^{-1}}{\sum }w\right) $

$\ \ \ \ \ \ \ \ \ \ \ \ \ \ \ \ \ \ \ \ \ \ \ \ \ \ \ \ \ \ \ \ \ \ \ \ \ \
\ +\left( \underset{\left| w\right| =m-1,\,\,End(w)\neq a}{\sum }w\right)
ba^{-1}\left( \underset{\left| w\right| =n-1,\,\,Init(w)\neq b^{-1}}{\sum }%
w\right) $

\strut

and

\strut

(3.3.1.8)

\strut

$\ \ \ W_{m+n-2}=\left( \underset{\left| w\right| =m-3,\,\,End(w)\neq a^{-1}%
}{\sum }w\right) b^{-1}X_{n}+X_{m}a\left( \underset{\left| w\right|
=n-3,\,\,Init(w)\neq b}{\sum }w\right) .$

\strut

Now, we will define a function $F_{pq}$ from $\Bbb{N\times N}$ to $B$.

\strut

\begin{definition}
Define a function $F:\Bbb{N}\times \Bbb{N}\rightarrow B$ by

\strut 

\begin{center}
$F_{pq}(k,l)=E\left( \left( \underset{\left| w\right| =k,\,\,End(w)=p}{\sum }%
w\right) \left( \underset{\left| w^{\prime }\right| =l,\,Init(w^{\prime })=q%
}{\sum }w^{\prime }\right) \right) ,$
\end{center}

\strut 

where $p,q\in \{a,b,a^{-1},b^{-1}\},$ where $End(w)$ and $Init(w)$ mean the
end letter of the word $w$ and initial letter of the word $w,$ respectively.
\end{definition}

\strut

\begin{definition}
Let $p,q\in \{a,a^{-1},b,b^{-1}\}$ and let $w=p_{1}...p_{k}$ be a word with
length $k$ in $\{a,a^{-1},b,b^{-1}\}.$ We define the following relation
denoted by ''$\lhd \,$'' ;

\strut 

\begin{center}
$pq\lhd w=p_{1}...p_{k}\overset{def}{\Longleftrightarrow }\exists \,j\in
\{1,...,k-1\}$ s.t $pq=p_{j}p_{j+1}$ \ and \ $pq\neq e$
\end{center}
\end{definition}

\strut

For example,

\strut

\begin{center}
$pq\lhd h=aba^{-1}b^{-1}$
\end{center}

if and only if

\begin{center}
$pq=ab$ \ or \ $pq=ba^{-1}$ \ or \ $pq=a^{-1}b^{-1}.$
\end{center}

\strut

Again recall that, for all $n\in \Bbb{N},$

\strut

\begin{center}
$E(X_{n})=\left\{ 
\begin{array}{lll}
h^{\frac{n}{4}}+h^{-\,\frac{n}{4}} &  & \text{if \ }4\mid n \\ 
&  &  \\ 
0_{B} &  & \text{otherwise.}
\end{array}
\right. $
\end{center}

\strut

Then we have the following lemmas ;

\strut

\begin{lemma}
Let $p,q\in \{a,b,a^{-1},b^{-1}\}$ and $k,l\in \Bbb{N}$ (sufficiently big).
Then

\strut 

\begin{center}
$F_{pq}(k,l)=0_{B},$ whenever $pq\ntriangleleft h=aba^{-1}b^{-1}.$
\end{center}
\end{lemma}

\strut

\begin{proof}
Suppose that $pq\ntriangleleft aba^{-1}b^{-1}.$ Then

\strut

\begin{center}
$F_{pq}(k,l)=E\left( \left( \underset{\left| w\right| =k-1,\,\,End(w)\neq
p^{-1}}{\sum }w\right) pq\left( \underset{\left| w^{\prime }\right|
=l-1,\,\,Init(w^{\prime })\neq q^{-1}}{\sum }\right) \right) =0_{B},$
\end{center}

\strut

since every word $W=l_{i_{1}}....l_{i_{k-1}}pql_{j_{1}}...l_{j_{l-1}}$
cannot be $h^{\frac{k+l}{4}},$ where $l_{i_{1}}...l_{i_{k-1}}$ is the word
with length $k-1$ such that $l_{i_{k-1}}\neq p^{-1}$ and $%
l_{j_{1}}...l_{j_{l-1}}$ is the word with length $l-1$ such that $%
l_{j_{1}}\neq q^{-1},$ we can get the above equality.
\end{proof}

\strut \strut \strut \strut

\begin{lemma}
We have the following equalities ;

\strut 

(1)

\begin{center}
$F_{ab}(k,l)=\left\{ 
\begin{array}{lll}
h^{\frac{k+l}{4}}+h^{-\,\frac{k+l}{4}} &  & \text{if \ }4\mid (k-1)\text{
and }4\mid k+l \\ 
&  &  \\ 
0_{B} &  & \text{otherwise.}
\end{array}
\right. $
\end{center}

\strut 

(2)

\begin{center}
$F_{ba^{-1}}(k,l)=\left\{ 
\begin{array}{lll}
h^{\frac{k+l}{4}}+h^{-\,\frac{k+l}{4}} &  & \text{if \ }4\mid k\text{ and }%
4\mid k+l \\ 
&  &  \\ 
0_{B} &  & \text{otherwise.}
\end{array}
\right. $
\end{center}

\strut 

(3)

\begin{center}
$F_{a^{-1}b^{-1}}(k,l)=\left\{ 
\begin{array}{lll}
h^{\frac{k+l}{4}}+h^{-\,\frac{k+l}{4}} &  & \text{if \ }4\mid (l-1)\text{
and }4\mid k+l \\ 
&  &  \\ 
0_{B} &  & \text{otherwise.}
\end{array}
\right. $
\end{center}
\end{lemma}

\strut

\strut \strut

\begin{lemma}
We have the following recurrence relation for $F_{p,q}(k,l),$ where $p,q\in
\{a,b,a^{-1},b^{-1}\}$ and $k,l\in \Bbb{N}$ (sufficiently large) ;

\strut 

(4) $F_{aa^{-1}}(k,l)=F_{aa^{-1}}(k-1,l-1)+F_{ab}(k-1,l-1)$

$\ \ \ \ \ \ \ \ \ \ \ \ \ \ \ \
+F_{ba^{-1}}(k-1,l-1)+F_{bb^{-1}}(k-1,l-1)+F_{b^{-1}b}(k-1,l-1).$

\strut 

(5) $F_{bb^{-1}}(k,l)=F_{bb^{-1}}(k-1,l-1)+F_{aa^{-1}}(k-1,l-1)$

$\ \ \ \ \ \ \ \ \ \ \ \ \ \ \
+F_{a^{-1}a}(k-1,l-1)+F_{ba^{-1}}(k-1,l-1)+F_{a^{-1}b^{-1}}(k-1,l-1).$

\strut 

(6) $F_{a^{-1}a}(k,l)=F_{bb^{-1}}(k-1,l-1)+F_{a^{-1}a}(k-1,l-1)$

$\ \ \ \ \ \ \ \ \ \ \ \ \ \ \ \ \ \ \ \ \ \ \ \ \
+F_{b^{-1}b}(k-1,l-1)+F_{a^{-1}b^{-1}}(k-1,l-1).$

\strut 

(7) $F_{b^{-1}b}(k,l)=$\strut $F_{aa^{-1}}(k-1,l-1)+F_{a^{-1}a}(k-1,l-1)$

$\ \ \ \ \ \ \ \ \ \ \ \ \ \ \ \ \ \ \ \ \ \ \ \ \
+F_{b^{-1}b}(k-1,l-1)+F_{ab}(k-1,l-1).$
\end{lemma}

\strut

By the previous lemmas, we can get three equalities (1), (2) and (3) and
four recurrence relations (4) $\sim $ (7). Again, by the previous lemmas, we
can conclude that

\strut

\strut

\begin{theorem}
Let $k,\,l\in \Bbb{N}$ be sufficeintly big and let $p,q\in
\{a,b,a^{-1},b^{-1}\}.$ Then

\strut 

(i) \ If $pq\ntriangleleft h=aba^{-1}b^{-1},$ then $F_{pq}(k,l)=0_{B}.$

\strut 

(ii) \ If $pq\lhd h=aba^{-1}b^{-1},$ then

$\strut $

\begin{center}
$F_{pq}(k,l)=\left\{ 
\begin{array}{ll}
\,\,\,h^{\frac{k+l}{4}}, & \,\,\,\text{if }pq=ab,\text{ }4\mid (k-1)\text{
\& }4\mid (k+l) \\ 
\,\,\,h^{\frac{k+l}{4}}, & \,\,\,\text{if }pq=ba^{-1},\text{ }4\mid k\text{
\& }4\mid k+l \\ 
\begin{array}{l}
h^{\frac{k+l}{4}}, \\ 
0_{B}
\end{array}
& 
\begin{array}{l}
\text{if }pq=a^{-1}b^{-1},\text{ }4\mid (l-1)\text{ \& }4\mid (k+l) \\ 
\text{otherwise.}
\end{array}
\end{array}
\right. $
\end{center}

\strut 

(iii) If $pq=e,$ then we have the following recurrence relation ;

$\strut $

\begin{center}
$F_{pq}(k,l)=\underset{r,s\in \{a,b,a^{-1},b^{-1}\},\,\,(r,s)\neq (q,p)}{%
\sum }F_{rs}(k-1,l-1).$
\end{center}
\end{theorem}

\strut

\begin{proof}
(i) is proved by Lemma 3.12 and (ii) is proved by Lemma 3.13. Now, by using
the results (i) and (ii), we can characterize Lemma 3.14. The
characterization is the statement (iii) ;

\strut

Suppose that $pq=e=aa^{-1}=a^{-1}a=bb^{-1}=b^{-1}b.$ Then

\strut

$F_{pq}(k,l)=E\left( \left( \underset{\left| w\right| =k,\,\,End(w)=p}{\sum }%
w\right) \left( \underset{\left| w^{\prime }\right| =l,\,\,Init(w^{\prime
})=q}{\sum }w\right) \right) $

\strut

$\ =E\left( \left( \underset{\left| w\right| =k-1,\,\,End(w)\neq p^{-1}}{%
\sum }w\right) \left( \underset{\left| w^{\prime }\right|
=l-1,\,\,Init(w^{\prime })=q^{-1}}{\sum }w\right) \right) $

\strut

$\ =\underset{r,s\in \{a,b,a^{-1},b^{-1}\},\,r\neq p^{-1},\,s\neq q^{-1}}{%
\sum }E\left( \left( \underset{\left| w\right| =k-1,\,\,End(w)=r}{\sum }%
w\right) \left( \underset{\left| w^{\prime }\right| =l-1,\,\,Init(w^{\prime
})=s}{\sum }w\right) \right) $

\strut

(3.15.1)

\strut

$\ =\underset{r,s\in \{a,b,a^{-1},b^{-1}\},\,r\neq p^{-1},\,s\neq q^{-1}}{%
\sum }F_{rs}(k-1,l-1).$

\strut

Since $pq=e,$ $q=p^{-1}$ and hence

\strut

\begin{center}
$r\neq p^{-1}$ \ and \ $s\neq q^{-1}\Longleftrightarrow r\neq q$ \ and \ $%
s\neq p.$
\end{center}

\strut

Therefore, the formula (3.15.1) is equivalent to

\strut

(3.15.2)

\begin{center}
$\underset{r,s\in \{a,b,a^{-1},b^{-1}\},\,\,(r,s)\neq (q,p)}{\sum }%
F_{rs}(k-1,l-1).$
\end{center}
\end{proof}

\strut

Now, we will consider the case $X_{m}h^{-1}X_{n}.$ But this case will be
very similar to the previous case.

\strut \strut

\begin{theorem}
Let $k,\,l\in \Bbb{N}$ be sufficeintly big and let $p,q\in
\{a,b,a^{-1},b^{-1}\}.$ Then

\strut 

(i) \ If $pq\ntriangleleft h^{-1}=bab^{-1}a^{-1},$ then $F_{pq}(k,l)=0_{B}.$

\strut 

(ii) \ If $pq\lhd h^{-1}=bab^{-1}a^{-1},$ then

$\strut $

\begin{center}
$F_{pq}(k,l)=\left\{ 
\begin{array}{ll}
\,\,\,h^{-\frac{k+l}{4}}, & \,\,\,\text{if }pq=ba,\text{ }4\mid (k-1)\text{
\& }4\mid (k+l) \\ 
\,\,\,h^{-\frac{k+l}{4}}, & \,\,\,\text{if }pq=ab^{-1},\text{ }4\mid k\text{
\& }4\mid k+l \\ 
\begin{array}{l}
h^{-\frac{k+l}{4}}, \\ 
0_{B}
\end{array}
& 
\begin{array}{l}
\text{if }pq=b^{-1}a^{-1},\text{ }4\mid (l-1)\text{ \& }4\mid (k+l) \\ 
\text{otherwise.}
\end{array}
\end{array}
\right. $
\end{center}

\strut 

(iii) If $pq=e,$ then we have the following recurrence relation ;

$\strut $

\begin{center}
$F_{pq}(k,l)=\underset{r,s\in \{a,b,a^{-1},b^{-1}\},\,\,(r,s)\neq (q,p)}{%
\sum }F_{rs}(k-1,l-1).$
\end{center}

$\square $
\end{theorem}

\strut

Now, by using the above equalities and recurrence relation, we can compute
the $E\left( X_{m}hX_{n}\right) $. We need the following definition ; Recall
that

\strut

$X_{m}hX_{n}$

$\strut \strut $

$\ \ \ =W_{m+n+4}+W_{m+n+2}+W_{m+n}+W_{m+n-2}$

$\strut \strut $

$\ \ \ \ \ \ \ \ \ \ \ \ \ \ \ +\left( \underset{\left| w\right|
=m-4,\,\,End(w)\neq b}{\sum }w\right) \left( X_{n}\right) +X_{m}\left( 
\underset{\left| w\right| =n-4,\,\,Init(w)\neq a}{\sum }w\right) $ \ \ \ \
(3.3.1.2)

\strut

$\ \ \ \ \ \ \ \ \ \ \ \ \ \ \ +\left( \underset{\left| w\right|
=m-3,\,\,End(w)\neq a^{-1}}{\sum }w\right) \left( \underset{\left| w\right|
=n-1,\,\,Init(w)\neq b^{-1}}{\sum }w\right) $ \ \ \ \ \ \ \ \ \ \ \ \ \ \ \
\ \ (3.3.1.3)

$\ \ \ \ \ \ \ \ \ \ \ \ \ \ \ +\left( \underset{\left| w\right|
=m-1,\,\,End(w)\neq a}{\sum }w\right) \left( \underset{\left| w\right|
=n-3,\,\,Init(w)\neq b}{\sum }w\right) $ \ \ \ \ \ \ \ \ \ \ \ \ \ \ \ \ \ \
\ \ \ \ (3.3.1.4)

$\ \ \ \ \ \ \ \ \ \ \ \ \ \ \ +\left( \underset{\left| w\right|
=m-2,\,\,End(w)\neq a}{\sum }w\right) \left( \underset{\left| w\right|
=n-2,\,\,Init(w)\neq a^{-1}}{\sum }w\right) ,$ \ \ \ \ \ \ \ \ \ \ \ \ \ \ \
\ \ \ \ (3.3.1.5)

\strut

where

\strut

(3.3.1.6)

\strut

$\ \ W_{m+n+2}=\left( \underset{\left| w\right| =m-1,\,\,End(w)\neq a}{\sum }%
w\right) (ba^{-1}b^{-1})X_{n}+X_{m}(aba^{-1})\left( \underset{\left|
w\right| =n-1,\,\,Init(w)\neq b^{-1}}{\sum }w\right) ,$

\strut \strut

(3.3.1.7)

\strut

$\ \ W_{m+n}=\left( \underset{\left| w\right| =m-2,\,\,End(w)\neq b}{\sum }%
w\right) (a^{-1}b^{-1})X_{n}+X_{m}(ab)\left( \underset{\left| w\right|
=n-2,\,\,Init(w)\neq a^{-1}}{\sum }w\right) $

$\ \ \ \ \ \ \ \ \ \ \ \ \ \ \ \ \ \ \ \ \ \ \ \ \ \ \ \ \ \ \ \ \ \ \ \ \ \
\ +\left( \underset{\left| w\right| =m-1,\,\,End(w)\neq a}{\sum }w\right)
ba^{-1}\left( \underset{\left| w\right| =n-1,\,\,Init(w)\neq b^{-1}}{\sum }%
w\right) $

\strut

and

\strut

(3.3.1.8)

\strut

$\ \ \ W_{m+n-2}=\left( \underset{\left| w\right| =m-3,\,\,End(w)\neq a^{-1}%
}{\sum }w\right) b^{-1}X_{n}+X_{m}a\left( \underset{\left| w\right|
=n-3,\,\,Init(w)\neq b}{\sum }w\right) .$

\strut

Now, define a certain generalization of $F_{pq}(k,l)$ ;

\strut

\begin{definition}
Let $p_{1},...,p_{N}\in \{a,b,a^{-1},b^{-1}\}$ and let $d,\,N<M\in \Bbb{N}.$
Define the relation $\lhd $ by

\strut 

\begin{center}
$p_{1}...p_{N}\lhd l_{1}.....l_{M}\overset{def}{\Longleftrightarrow }\exists
\,j\in \{1,...,M-N-1\}$ s.t $p_{1}...p_{N}=l_{j+1}...l_{j+N}.$
\end{center}

\strut 

Also, define a map $F_{p_{1}...p_{j}]\_[p_{j+1}...p_{N}}:\Bbb{N\times N}%
\rightarrow B,$ for all $j=1,...,N-1$ by

\strut 

$F_{p_{1}...p_{j}]\_[p_{j+1}...p_{N}}(k,l)=E\left( \left( \underset{\left|
w\right| =k,\,\,End(w)=p_{1}...p_{j}}{\sum }w\right) \left( \underset{\left|
w^{\prime }\right| =l,\,\,Init(w^{\prime })=p_{j+1}...p_{N}}{\sum }w^{\prime
}\right) \right) .$
\end{definition}

\strut

For example,

\begin{center}
\strut $aba^{-1}\lhd h$ \ or \ $bab^{-1}\lhd h$
\end{center}

and

\begin{center}
$a^{2}b\ntriangleleft h$
\end{center}

etc.

\strut

By using the above new definition, we can rewrite that

\strut

\begin{center}
$F_{pq}(k,l)=F_{p]\_[q}(k,l),$
\end{center}

\strut

for all $p,q\in \{a,b,a^{-1},b^{-1}\},$ where $k$ and $l$ are sufficiently
large natural numbers.

\strut \strut

\begin{proposition}
Suppose that $k,l$ and $d>N\geq 3$ in $\Bbb{N}.$ Then

\strut 

$F_{p_{1}...p_{j}]\_[p_{j+1}...p_{N}}(k,l)$

\strut 

\begin{center}
$=\left\{ 
\begin{array}{ll}
\left\{ 
\begin{array}{lll}
h^{\frac{k+l}{4}} &  & \text{if \ }p_{1}=a,\,\text{\ }4\mid (k-j)\text{ \& }%
4\mid (k+l) \\ 
h^{\frac{k+l}{4}} &  & \text{if \ }p_{1}=b,\text{ \ }4\mid (k-j+1)\text{ \& }%
4\mid (k+l) \\ 
h^{\frac{k+l}{4}} &  & \text{if \ }p_{1}=a^{-1},\text{ }4\mid (k-j+2)\text{
\& }4\mid (k+l) \\ 
h^{\frac{k+l}{4}} &  & \text{if }p_{1}=b^{-1},\text{ }4\mid (k-j+3)\text{ \& 
}4\mid (k+l). \\ 
0_{B} &  & \text{otherwise.}
\end{array}
\right\}  & \text{if }p_{1}...p_{N}\lhd h^{d} \\ 
&  \\ 
0_{B} & \text{if }p_{1}...p_{N}\ntriangleleft h^{d}
\end{array}
\right. $
\end{center}

\strut 

where $j\in \{1,...,N-1\}$. Moreover, if $\ 4\mid (k-j+i)$ and $4\mid (k+l),$
for the fixed $i\in \{0,1,...,3\},$ then, $4\nmid (k-j+i^{\prime }),$
whenever $i^{\prime }\neq i.$ (i.e, if \ $4\mid (k-j),$ then \ $4\nmid
(k-j+i),$ for all $i=1,2,3.$)
\end{proposition}

\strut

\strut

\begin{corollary}
Suppose that $d>N\geq 3$ in $\Bbb{N}.$ Then

\strut 

$F_{p_{1}...p_{j}]\_[p_{j+1}...p_{N}}(k,l)=$

\strut 

\begin{center}
$\left\{ 
\begin{array}{ll}
\left\{ 
\begin{array}{lll}
h^{-\,\frac{k+l}{4}} &  & \text{if \ }p_{1}=b,\,\text{\ }4\mid (k-j)\text{
\& }4\mid (k+l) \\ 
h^{-\,\frac{k+l}{4}} &  & \text{if \ }p_{1}=a,\text{ \ }4\mid (k-j+1)\text{
\& }4\mid (k+l) \\ 
h^{-\,\frac{k+l}{4}} &  & \text{if \ }p_{1}=b^{-1},\text{ }4\mid (k-j+2)%
\text{ \& }4\mid (k+l) \\ 
h^{-\,\frac{k+l}{4}} &  & \text{if }p_{1}=a^{-1},\text{ }4\mid (k-j+3)\text{
\& }4\mid (k+l). \\ 
0_{B} &  & \text{otherwise.}
\end{array}
\right\}  & \text{if }p_{1}...p_{N}\lhd h^{-d} \\ 
&  \\ 
0_{B} & \text{if }p_{1}...p_{N}\ntriangleleft h^{-d}
\end{array}
\right. $

\strut 
\end{center}

where $j\in \{1,...,N-1\}$. \ $\square $
\end{corollary}

\strut \strut

\begin{corollary}
Let $p,q\in \{a,b,a^{-1},b^{-1}\}$ and $k,l\in \Bbb{N},$ sufficiently large.
If $pq\neq e,$ then

\strut 

$F_{pq}(k,l)=F_{p]\_[q}(k,l)$

\strut 

\begin{center}
$=\left\{ 
\begin{array}{lll}
h^{\frac{k+l}{4}} &  & \text{if \ }p=a,\,\text{\ }4\mid (k-1)\text{ \& }%
4\mid (k+l) \\ 
h^{\frac{k+l}{4}} &  & \text{if \ }p=b,\text{ \ }4\mid k\text{ \& }4\mid
(k+l) \\ 
h^{\frac{k+l}{4}} &  & \text{if \ }p=a^{-1},\text{ }4\mid (k+1)\text{ \& }%
4\mid (k+l) \\ 
h^{\frac{k+l}{4}} &  & \text{if }p=b^{-1},\text{ }4\mid (k+2)\text{ \& }%
4\mid (k+l). \\ 
0_{B} &  & \text{otherwise.}
\end{array}
\right. $

\strut 
\end{center}

and

\begin{center}
\strut 
\end{center}

$F_{pq}(k,l)=F_{p]\_[q}(k,l)$

\begin{center}
\strut 

$=\left\{ 
\begin{array}{lll}
h^{-\,\frac{k+l}{4}} &  & \text{if \ }p=b,\,\text{\ }4\mid (k-1)\text{ \& }%
4\mid (k+l) \\ 
h^{-\,\frac{k+l}{4}} &  & \text{if \ }p=a,\text{ \ }4\mid k\text{ \& }4\mid
(k+l) \\ 
h^{-\,\frac{k+l}{4}} &  & \text{if \ }p=b^{-1},\text{ }4\mid (k+1)\text{ \& }%
4\mid (k+l) \\ 
h^{-\,\frac{k+l}{4}} &  & \text{if }p=a^{-1},\text{ }4\mid (k+2)\text{ \& }%
4\mid (k+l). \\ 
0_{B} &  & \text{otherwise.}
\end{array}
\right. $

\strut 

In particular, this is the generalization of (ii) of the Theorem 3.6. $%
\square $
\end{center}
\end{corollary}

\strut

\strut

Now, consider the case of $F_{p_{1}^{-1}\_p_{1}...p_{N}}(k,l)$ and $%
F_{p_{1}...p_{N}\_p_{N}^{-1}}(k,l).$ Of course, we will restrict our
interests to the case when $p_{1}...p_{N}\lhd h^{d}.$

\strut

\begin{proposition}
Let $k,l$ and $d>N>2$ be in $\Bbb{N}.$ Assume that $p_{1}...p_{N}\lhd h^{d}.$
Then

\strut 

\begin{center}
$F_{p_{1}^{-1}]\_[p_{1}...p_{N}}(k,l)=F_{p_{2}^{-1}]%
\_[p_{2}...p_{N}}(k-1,l-1)$
\end{center}

and

\begin{center}
$F_{p_{1}...p_{N}]\_[p_{N}^{-1}}(k,l)=F_{p_{1}...p_{N-1}]%
\_[p_{N-1}^{-1}}(k-1,l-1).$
\end{center}
\end{proposition}

\strut

\strut \strut

Now, we need the following new function from $\Bbb{N\times N}$ to $B$ ;

\strut

\begin{definition}
Let $p_{1},...,p_{M}\in \{a,b,a^{-1},b^{-1}\}$ and assume that $%
p_{1}...p_{M}\lhd h^{d}.$ Let $p^{1},...,p^{i},p^{j},...,p^{N}\in
\{a,b,a^{-1},b^{-1}\}$ and let $p^{1}...p^{i}$ and $p^{j}...p^{N}$ be words
in $\{a,b,a^{-1},b^{-1}\}.$ Always assume that $\left|
p^{1}...p^{i}p^{j}...p^{N}\right| <\left| p_{1}...p_{M}\right| .$ Define a
map

\strut 

\begin{center}
$F_{p^{1}...p^{i}]\_<p_{1}...p_{M}>\_[p^{j}...p^{N}}:\Bbb{N\times N}%
\rightarrow B$
\end{center}

by

\strut 

$F_{p^{1}...p^{i}]\_<p_{1}...p_{M}>\_[p^{j}...p^{N}}(k,M,l)$

\strut 

\begin{center}
$=E\left( \left( \underset{\left| w\right| =k,\,\,End(w)=p^{1}...p^{i}}{\sum 
}w\right) (p_{1}...p_{M})\left( \underset{\left| w^{\prime }\right|
=l,\,Init(w^{\prime })=p^{j}...p^{M}}{\sum }w^{\prime }\right) \right) ,$
\end{center}

\strut 

for all $k,l\,\,\in \Bbb{N}.$
\end{definition}

\strut

For example,

\begin{center}
$F_{p^{1}...p^{i}]\_<e>\_[p^{j}...p^{N}}(k,0,l)=F_{p^{1}...p^{i}]%
\_[p^{j}...p^{N}}(k,l).$
\end{center}

\strut

\begin{proposition}
By using the same notations in the previous definition, we have that

\strut 

(1) If $p^{1}...p^{i}p_{1}...p_{M}p^{j}...p^{N}\ntriangleleft \,h^{d},$ then

\strut 

\begin{center}
$F_{p^{1}...p^{i}]\_<p_{1}...p_{M}>\_[p^{j}...p^{N}}(k,M,l)=0_{B}.$

\strut 
\end{center}

\ \ \ (Remaind that $\left| p^{1},,,p^{i}p^{j}...p^{N}\right| <\left|
p_{1}...p_{M}\right| \ $and hence $p^{1}...p^{i}p_{1}...p_{M}p^{j}...p^{N}%
\neq e.$)

\strut 

\ \ \ \ \strut If $p^{1}...p^{i}p_{1}...p_{M}p^{j}...p^{N}\lhd \,h^{d},$ then

\strut 

\ \ \ \ \ $F_{p^{1}...p^{i}]\_<p_{1}...p_{M}>\_[p^{j}...p^{N}}(k,M,l)$

\strut 

\begin{center}
$=\left\{ 
\begin{array}{lll}
h^{\frac{k+M+l}{4}} &  & \text{if \ }p_{1}=a,\text{ }4\mid (k-i)\text{ \& }%
4\mid (k+M+l) \\ 
h^{\frac{k+M+l}{4}} &  & \text{if \ }p_{1}=b,\text{ }4\mid (k-j+1)\text{ \& }%
4\mid (k+M+l) \\ 
h^{\frac{k+M+l}{4}} &  & \text{if \ }p_{1}=a^{-1},\text{ }4\mid (k-j+2)\text{
\& }4\mid (k+M+l) \\ 
h^{\frac{k+M+l}{4}} &  & \text{if }p_{1}=b^{-1},\text{ }4\mid (k-j+3)\text{
\& }4\mid (k+M+l) \\ 
0_{B} &  & \text{otherwise.}
\end{array}
\right. $
\end{center}

\strut 

(2) Assume that $p^{i}=p_{1}^{-1}.$ Then

\strut 

\begin{center}
$F_{p^{1}...p^{i}]\_<p_{1}...p_{M}>%
\_[p^{j}...p^{N}}(k,M,l)=F_{p_{1}...p^{i-1}]\_<p_{2}...p_{M}>%
\_[p^{j}...p^{N}}(k-1,M-1,l).$

\strut 
\end{center}

In particular, if $p^{1}...p^{i}p_{1}...p_{M}p^{j}...p^{N}\,\lhd \,h^{d},$
then $F_{p^{1}...p^{i}]\_<p_{1}...p_{M}>\_[p^{j}...p^{N}}(k,M,l)=0_{B}.$

\strut 

(3) Assume that $p^{j}=p_{M}^{-1}.$ Then

\strut 

\begin{center}
$F_{p^{1}...p^{i}]\_<p_{1}...p_{M}>\_[p^{j}...p^{N}}(k,M,l)=F_{p^{1}...p^{i}%
\_p_{1}...p_{M-1}\_p^{j+1}...p^{N}}(k,M-1,l-1).$
\end{center}

\strut 

In particular, if $p^{1}...p^{i}p_{1}...p_{M}p^{j}...p^{N}\,\lhd \,h^{d},$
then $F_{p^{1}...p^{i}]\_<p_{1}...p_{M}>\_[p^{j}...p^{N}}(k,M,l)=0_{B}.$

\strut 

(4) Let $p^{i}=p_{1}^{-1}$ and $p^{j}=p_{M}^{-1}.$ Then

\strut 

\begin{center}
$F_{p^{1}...p^{i}]\_<p_{1}...p_{M}>%
\_[p^{j}...p^{N}}(k,M,l)=F_{p^{1}...p^{i-1}]\_<p_{2}...p_{M-1}>%
\_[p^{j+1}...p^{N}}(k-1,M-1,l-1).$
\end{center}

\strut 

In particular, if $p^{1}...p^{i}p_{1}...p_{M}p^{j}...p^{N}\,\lhd \,h^{d},$
then $F_{p^{1}...p^{i}]\_<p_{1}...p_{M}>\_[p^{j}...p^{N}}(k,M,l)=0_{B}.$
\end{proposition}

\strut

\strut \strut

\begin{proposition}
Let $p_{1}...p_{i}p_{i+1}...p_{M}\lhd \,h^{d}.$ Then

\strut 

(1) $F_{p_{M}^{-1}...p_{1}^{-1}]\_<p_{1}...p_{M}>\_[p^{j}...p^{N}}(k,M,l)$

\strut 

\begin{center}
$=\underset{p_{0}\in \{a,b,a^{-1},b^{-1}\},\,\,p_{0}\neq
p_{M},\,\,p_{0}p^{j}...p^{N}\lhd \,h^{d}}{\sum }F_{p_{0}]%
\_[p^{j}...p^{N}}(k-M,l).$
\end{center}

\strut 

where $p_{0}\neq p_{M}$ and $p_{0}p^{j}...p^{N}\lhd \,h^{d}.$

\strut 

(2) $F_{p^{1}...p^{i}]\_<p_{1}...p_{M}>\_[p_{M}^{-1}...p_{1}^{-1}}(k,M,l)$

\strut 

\begin{center}
$=\underset{p_{0}\in \{a,b,a^{-1},b^{-1}\},\,\,p_{0}\neq
p_{1},\,p^{1}...p^{i}p_{0}\lhd \,h^{d}}{\sum }F_{p_{0}]%
\_[p^{j}...p^{N}}(k-M,l).$
\end{center}

\strut 

(3) $F_{p_{i}^{-1}...p_{1}^{-1}]\_<p_{1}...p_{i}p_{i+1}...p_{M}>%
\_[p_{M}^{-1}...p_{i+1}^{-1}}(k,M,l)$

\strut 

\begin{center}
$=\underset{p\neq p_{i},\,\,q\neq p_{i+1}}{\sum }F_{p]\_[q}(k-i,\,l-(M-i)).$%
\strut 
\end{center}
\end{proposition}

\strut

\strut \strut \strut 

Now, we will compute $E\left( X_{m}hX_{n}\right) $ ;

\strut

\begin{lemma}
We have that

\strut 

(1) $E\left( W_{m+n+4}\right) =\left\{ 
\begin{array}{lll}
h^{m+n+1} &  & \text{if \ }4\mid m\text{ and }4\mid n \\ 
&  &  \\ 
0_{B} &  & \text{otherwise.}
\end{array}
\right. $

\strut 

(2) $E\left( W_{m+n+2}\right)
=F_{b^{-1}]\_[ba^{-1}b^{-1}}(m-1,n+3)+F_{aba^{-1}]\_[a}(m+3,n-1).$

\strut 

(3) $E\left( W_{m+n}\right) =\underset{p\neq b,\,\,q\neq a^{-1}}{\sum }%
F_{p]\_[q}(m-2,\,n-2).$

\strut 

(4) $E(W_{m+n-2})=\underset{p\neq a^{-1}}{\sum }F_{p]\_[b^{-1}}(m-3,n+1)+%
\underset{p\neq b}{\sum }F_{a]\_[p}(m+1,n-3)$

$\ \ \ \ \ \ \ +\underset{p\neq b,\,q\neq b^{-1}}{\sum }F_{p]\_<a>%
\_[q}(m-2,1,n-1)+\underset{p\neq a,\,q\neq a^{-1}}{\sum }F_{p]\_<b>%
\_[q}(m-1,1,n-2).$
\end{lemma}

\strut \strut 

Now, we will consider the general case $E(X_{m}h^{d}X_{n}),$ where $d\in 
\Bbb{N}.$

\strut

\begin{lemma}
Let $m,n,d\in \Bbb{N}.$ Then

\strut 

\begin{center}
$E\left( W_{m+n+4d}\right) =\left\{ 
\begin{array}{lll}
h^{\frac{m+n+4d}{4}} &  & \text{if }4\mid m\text{ and }4\mid n \\ 
&  &  \\ 
0_{B} &  & \text{otherwise.}
\end{array}
\right. $
\end{center}
\end{lemma}

\strut \strut \strut \strut 

\strut The above case is the NO-cancellation case of $X_{m}h^{d}X_{n}.$

\strut

\begin{lemma}
Let $m,n\in N$ be sufficiently big and $d\in \Bbb{N}.$ Then

\strut 

$E\left( W_{m+n+(4d-2j)}\right) =\underset{p\neq p_{2j}}{\sum }%
F_{p]\_[p_{2j+1}...p_{4d}}(m-2j,n+4d)$

\strut 

$\ \ \ \ \ \ \ \ \ \ +\underset{i,j\in
\{1,...,2d-1\},\,(m-i)+(n+j)=m+n+(4d-2j)}{\sum }$

\begin{center}
$\underset{p\neq p_{i},\,p^{\prime }\neq p_{j}}{\sum }F_{p]%
\_<p_{i+1}...p_{j-1}>\_[p^{\prime }}(m-i,\,j-1-i,\,n+j)$
\end{center}

\strut 

$\ \ \ \ \ \ \ \ \ \ +\underset{p\neq p_{4d-(2j+1)}}{\sum }%
F_{p_{1}...p_{4d-2j}]\_[p}(m+4d,n-2j).$
\end{lemma}

\strut

\strut We will call the above cancellation case a $(i+j)$-cancellation from
the left and from the right.

\strut \strut

\begin{lemma}
Let $m,n,d\in \Bbb{N}.$ Then

\strut 

\begin{center}
$
\begin{array}{ll}
E(W_{m+n-4d}) & =\underset{p\in \{a,b,a^{-1}\},p^{\prime }\in
\{a,b,a^{-1},b^{-1}\}}{\sum }F_{p]\_[p^{\prime }}(m-4d,n) \\ 
& \,\,\,\,\,\,\,\,+\underset{p\in \{a,b,a^{-1},b^{-1}\},p^{\prime }\in
\{b,a^{-1},b^{-1}\}}{\sum }F_{p]\_[p^{\prime }}(m,n-4d) \\ 
& \,\,\,\,\,\,\,\,+\sum_{i=1}^{4d}\,\underset{p\neq p_{i},\,p^{\prime }\neq
p_{i+1}}{\sum }F_{p]\_[p^{\prime }}(m-i,n-(4d-i)).
\end{array}
$
\end{center}
\end{lemma}

\strut \strut

\strut

\begin{lemma}
Let $m,n,N\in \Bbb{N}.$ If $p_{i_{1}}...p_{i_{N}}\lhd h^{d},$ then

\strut 

$\ \ \ \sum_{j=0}^{2d-1}E(W_{m+n-2j})=\sum_{j=0}^{2d-1}\,\underset{%
k+N+l=m+n-2j}{\sum }\,\underset{i_{1}\neq i_{2}\neq ...\neq i_{N}\in
\{1,...,4d\}}{\sum }$

\strut 

$\ \ \ \ \ \ \ \ \ \ \ \ \ \ \ \ \ \ \ \ \ \ \ \ \ \ \underset{p,q\in
\{a,b,a^{-1},b^{-1}\},\,pp_{i_{1}}...p_{i_{N}}\lhd h^{d}}{\sum }%
\,\,F_{p]\_[p_{i_{1}}...p_{i_{N}}}(k,N,l).$

$\square $
\end{lemma}

\strut

Finally, we can get the $B$-functional value $E(X_{m}h^{d}X_{n})$ ;

\strut \strut

\begin{theorem}
Let $m,n,d\in \Bbb{N}.$ We have that

\strut 

\begin{center}
$
\begin{array}{ll}
E\left( X_{m}h^{d}X_{n}\right)  & =E(W_{m+n+4d})+%
\sum_{j=1}^{2d-1}E(W_{m+n+(4d-2j)}) \\ 
& \,\,\,\,\,\,\,\,\,\,\,\,+\sum_{j=0}^{2d-1}E(W_{m+n-2j})+E(W_{m+n-4d}),
\end{array}
$
\end{center}

\strut \strut 

where, by putting $h^{d}=p_{1}...p_{4d},$

\strut 

$\ \ \ \ E(W_{m+n+4d})=\left\{ 
\begin{array}{lll}
h^{\frac{m+n+4d}{4}} &  & \text{if }4\mid m\text{ and }4\mid n \\ 
&  &  \\ 
0_{B} &  & \text{otherwise,}
\end{array}
\right. $

\strut 

$\ \ \ \ \sum_{j=1}^{2d-1}E(W_{m+n+(4d-2j)})=\sum_{j=1}^{2d-1}\,\,\,\,(%
\underset{p\neq p_{2j}}{\sum }F_{p]\_[p_{2j+1}...p_{4d}}(m-2j,n+4d)$

\strut 

$\ \ \ \ \ \ \ \ \ \ +\underset{i,j\in
\{1,...,2d-1\},\,(m-i)+(n+j)=m+n+(4d-2j)}{\sum }$

\begin{center}
$\underset{p\neq p_{i},\,p^{\prime }\neq p_{j}}{\sum }F_{p]%
\_<p_{i+1}...p_{j-1}>\_[p^{\prime }}(m-i,\,j-1-i,\,n+j)$
\end{center}

\strut 

$\ \ \ \ \ \ \ \ \ \ +\underset{p\neq p_{4d-(2j+1)}}{\sum }%
F_{p_{1}...p_{4d-2j}]\_[p}(m+4d,n-2j)\,\,\,\,\,\,)$

\strut 

$\sum_{j=0}^{2d-1}E(W_{m+n-2j})=\sum_{j=0}^{2d-1}\,\underset{k+N+l=m+n-2j}{%
\sum }\,\underset{i_{1}\neq i_{2}\neq ...\neq i_{N}\in \{1,...,4d\}}{\sum }$

\strut \strut 

\begin{center}
$\underset{p,q\in \{a,b,a^{-1},b^{-1}\},\,pp_{i_{1}}...p_{i_{N}}\lhd h^{d}}{%
\sum }\,\,F_{p]\_[p_{i_{1}}...p_{i_{N}}}(k,N,l).$
\end{center}

\strut \strut 

and

\strut 

$\ \ \ \ 
\begin{array}{ll}
E(W_{m+n-4d}) & =\underset{p\in \{a,b,a^{-1}\},p^{\prime }\in
\{a,b,a^{-1},b^{-1}\}}{\sum }F_{p]\_[p^{\prime }}(m-4d,n) \\ 
& \,\,\,\,\,\,\,\,+\underset{p\in \{a,b,a^{-1},b^{-1}\},p^{\prime }\in
\{b,a^{-1},b^{-1}\}}{\sum }F_{p]\_[p^{\prime }}(m,n-4d) \\ 
& \,\,\,\,\,\,\,\,+\sum_{i=1}^{4d}\,\underset{p\neq p_{i},\,p^{\prime }\neq
p_{i+1},\,p_{i}p_{i+1}\lhd \,h^{d}}{\sum }F_{p]\_[p^{\prime }}(m-i,n-(4d-i)).
\end{array}
$

$\square $
\end{theorem}

\strut \strut \strut

For considering $h^{-d}$ ($d\in \Bbb{N}$), we have the following result,
like the former theorem ;

\strut \strut

\begin{theorem}
Let $m,n,d\in \Bbb{N}.$ We have that

\strut 

\begin{center}
$
\begin{array}{ll}
E\left( X_{m}h^{-d}X_{n}\right)  & =E(W_{m+n+4d})+%
\sum_{j=1}^{2d-1}E(W_{m+n+(4d-2j)}) \\ 
& \,\,\,\,\,\,\,\,\,\,\,\,+\sum_{j=0}^{2d-1}E(W_{m+n-2j})+E(W_{m+n-4d}),
\end{array}
$
\end{center}

\strut \strut 

where, by putting $h^{-d}=p_{1}...p_{4d},$

\strut 

$\ \ \ \ E(W_{m+n+4d})=\left\{ 
\begin{array}{lll}
h^{-\,\frac{m+n+4d}{4}} &  & \text{if }4\mid m\text{ and }4\mid n \\ 
&  &  \\ 
0_{B} &  & \text{otherwise,}
\end{array}
\right. $

\strut 

$\ \ \ \ \sum_{j=1}^{2d-1}E(W_{m+n+(4d-2j)})=\sum_{j=1}^{2d-1}\,\,\,\,(%
\underset{p\neq p_{2j}}{\sum }F_{p]\_[p_{2j+1}...p_{4d}}(m-2j,n+4d)$

\strut 

$\ \ \ \ \ \ \ \ \ \ +\underset{i,j\in
\{1,...,2d-1\},\,(m-i)+(n+j)=m+n+(4d-2j)}{\sum }$

\begin{center}
$\underset{p\neq p_{i},\,p^{\prime }\neq p_{j}}{\sum }F_{p]%
\_<p_{i+1}...p_{j-1}>\_[p^{\prime }}(m-i,\,j-1-i,\,n+j)$
\end{center}

\strut 

$\ \ \ \ \ \ \ \ \ \ +\underset{p\neq p_{4d-(2j+1)}}{\sum }%
F_{p_{1}...p_{4d-2j}]\_[p}(m+4d,n-2j)\,\,\,\,\,\,)$

\strut 

$\sum_{j=0}^{2d-1}E(W_{m+n-2j})=\sum_{j=0}^{2d-1}\,\underset{k+N+l=m+n-2j}{%
\sum }\,\underset{i_{1}\neq i_{2}\neq ...\neq i_{N}\in \{1,...,4d\}}{\sum }$

\strut \strut 

\begin{center}
$\underset{p,q\in \{a,b,a^{-1},b^{-1}\},\,pp_{i_{1}}...p_{i_{N}}\lhd h^{d}}{%
\sum }\,\,F_{p]\_[p_{i_{1}}...p_{i_{N}}}(k,N,l).$
\end{center}

\strut \strut 

and

\strut 

$\ \ \ \ 
\begin{array}{ll}
E(W_{m+n-4d}) & =\underset{p\in \{a,a^{-1},b^{-1}\},p^{\prime }\in
\{a,b,a^{-1},b^{-1}\}}{\sum }F_{p]\_[p^{\prime }}(m-4d,n) \\ 
& \,\,\,\,\,\,\,\,+\underset{p\in \{a,b,a^{-1},b^{-1}\},p^{\prime }\in
\{a,a^{-1},b^{-1}\}}{\sum }F_{p]\_[p^{\prime }}(m,n-4d) \\ 
& \,\,\,\,\,\,\,\,+\sum_{i=1}^{4d}\,\underset{p\neq p_{i},\,p^{\prime }\neq
p_{i+1},\,p_{i}p_{i+1}\lhd \,h^{d}}{\sum }F_{p]\_[p^{\prime }}(m-i,n-(4d-i)).
\end{array}
$

$\square $
\end{theorem}

\strut \strut \strut

\strut

\strut

\subsection{Recerrence Relation For $E\left(
X_{m_{1}}h^{d_{2}}X_{m_{2}}h^{d_{3}}...h^{d_{n}}X_{m_{n}}\right) $}

\strut

\strut

In this section, we will compute more general form

\strut

\begin{center}
$E\left( X_{m_{1}}h^{d_{2}}X_{m_{2}}h^{d_{3}}...h^{d_{n}}X_{m_{n}}\right) ,$
\end{center}

\strut

where $m_{1},...,m_{n}\in \Bbb{N}$ and $d_{2},...,d_{n}\in \Bbb{Z}.$ First,
we will consider the most simple form, among them,

\strut

\begin{center}
$E\left( X_{m_{1}}h^{d_{2}}X_{m_{2}}h^{d_{3}}X_{m_{3}}\right) ,$
\end{center}

\strut

for $m_{1},m_{2},m_{3}\in \Bbb{N}$ and $d_{2},d_{3}\in \Bbb{Z}.$ Notice
that, by the evenness of $x=a+b+a^{-1}+b^{-1},$ we have the following
trivial condition ;

\strut

\begin{center}
$m_{1}+...+m_{n}\in N$ should be even.
\end{center}

\strut

Throughout this paper, we will assume that $m_{1}+...+m_{n}$ are even!
(Notice that each $m_{j}$ need not be even. For instance, we can have that $%
m_{1}=1,m_{2}=3$ and $m_{3}=2.$) Recall that, by the previous section, we
have that, for $m_{1},m_{2}\in \Bbb{N}$ and $d_{2}\in \Bbb{N},$

\strut

\begin{center}
$
\begin{array}{ll}
E\left( X_{m_{1}}h^{d_{2}}X_{m_{2}}\right) & =E(W_{m_{1}+m_{2}+4d_{2}})+%
\sum_{j=1}^{2d_{2}-1}E(W_{m_{1}+m_{2}+(4d_{2}-2j)}) \\ 
& \,\,\,\,\,\,\,\,\,\,\,\,+%
\sum_{j=0}^{2d_{2}-1}E(W_{m+n-2j})+E(W_{m+n-4d_{2}}),
\end{array}
$
\end{center}

\strut \strut

where each summands are determined recurssively.

\strut \strut

\strut

But, in $X_{m_{1}}h^{d_{2}}X_{m_{2}}h^{d_{3}}X_{m_{3}},$ there will be much
more terms which we have to consider. Also, different from the $%
X_{m_{1}}h^{d_{2}}X_{m_{2}},$ we have to consider the case when both $d_{2}$
and $d_{3}$ are positive integers or both $d_{2}$ and $d_{3}$ are negative
integers or either $d_{2}$ or $d_{3}$ is a positive integer and the other is
negative. So, we need the following new definition.

\strut

\begin{definition}
Let $m_{1},m_{2},m_{3},N_{2},N_{3}\in \Bbb{N}$. Define a map $\Phi _{0}$
from $\Bbb{N\times N\times N\times N\times N}$ to $B$

\strut 

\begin{center}
$F_{p_{m_{1}}^{E}]\_<p_{11}...p_{1N_{2}}>\_[p_{m_{2}}^{E}\_p_{m_{2}}^{I}]%
\_<p_{21}...p_{2N_{3}}>\_[p_{m_{3}}^{I}}(m_{1},N_{2},m_{2},N_{3},m_{3}),$
\end{center}

\strut 

denoted by $\Phi _{0}(m_{1},N_{2},m_{2},N_{3},m_{3}),$ by

\strut 

$\ \ \ \ \ E\,\left( \underset{\left| w\right| =m_{1},\,End(w)=p_{m_{1}}^{E}%
}{\sum }w\right) (p_{11}...p_{1N_{2}})\left( \underset{\left| w\right|
=m_{2},\,\,Init(w)=p_{m_{2}}^{I},\,End(w)=p_{m_{2}}^{E}}{\sum }w\right) $

$\ \ \ \ \ \ \ \ \ \ \ \ \ \ \ \ \ \ \ \ \ \ \ \ \ \ \ \ \ \ \ \ \ \ \ \ \ \
\ \ \ \ \ \ \ \ \ \ \ \ \ \ \ \ \ \ \cdot (p_{21}...p_{2N_{3}})\left( 
\underset{\left| w\right| =m_{3},\,\,Init(w)=p_{m_{3}}^{I}}{\sum }w\right) .$

\strut 

If $p_{11}...p_{1N_{2}}=e$ or $p_{21}...p_{2N_{3}}=e,$ then we have that

\strut 

$\ \Phi _{0}(m_{1},0,m_{2},N_{3},m_{3})$

$\ \ \ \ \ \ \ \ \ \
=F_{p_{m_{1}}^{E}]\_[p_{m_{2}}^{E}\_p_{m_{2}}^{I}]\_<p_{21}...p_{2N_{3}}>%
\_[p_{m_{3}}^{I}}(m_{1},0,m_{2},N_{3},m_{3})$

or

\strut 

$\ \Phi _{0}(m_{1},N_{2},m_{2},0,m_{3})$

$\ \ \ \ \ \ \ \ \ \
=F_{p_{m_{1}}^{E}]\_<p_{11}...p_{1N_{2}}>\_[p_{m_{2}}^{E}\_p_{m_{2}}^{I}]%
\_[p_{m_{3}}^{I}}(m_{1},N_{2},m_{2},0,m_{3}).$

\strut 

And if both of them are $e,$ then we can define

\strut 

$\ \Phi _{0}(m_{1},0,m_{2},0,m_{3})$

$\ \ \ \ \ \ \ \ \ \
=F_{p_{m_{1}}^{E}]\_[p_{m_{2}}^{E}\_p_{m_{2}}^{I}]%
\_[p_{m_{3}}^{I}}(m_{1},0,m_{2},0,m_{3}).$
\end{definition}

\strut

In the above definition of $\Phi _{0}(m_{1},N_{2},m_{2},N_{3},m_{3}),$ the
half-open bracket '' $[p$ '' and '' $p]$ '' mean that words with initial
letter $p$ and the words with ending letter $p.$ Also the bracket $<p...q>$
means the word $p...q$ and $[p^{I}\_p^{E}]$ means $\underset{%
Init(w)=p^{I},\,\,End(w)=p^{E}}{\sum }w.$

\strut \strut

Now, we will observe the above function $\Phi
_{0}(m_{1},N_{2},m_{2},N_{3},m_{3}).$

\strut

\begin{lemma}
Let $m_{1},m_{2},m_{3},N_{2},N_{3}\in \Bbb{N}.$ Assume that [$%
p_{m_{1}}^{E}\neq p_{11}^{-1}$ and $p_{m_{2}}^{I}\neq p_{1N_{2}}^{-1}$] or [$%
p_{m_{2}}^{E}\neq p_{21}^{-1}$ and $p_{m_{3}}^{I}\neq p_{2N_{3}}^{-2}$]$.$
Then

\strut 

\begin{center}
$\Phi _{0}(m_{1},N_{2},m_{2},N_{3},m_{3})=0_{B},$
\end{center}

\strut 

whenever $p_{11}...p_{1N_{2}}\ntriangleleft \,h^{d_{2}}$ or $%
p_{21}...p_{2N_{3}}\ntriangleleft \,h^{d_{3}},$ where $d_{2},d_{3}\in \Bbb{Z}%
\,\setminus \,\{0\}.$
\end{lemma}

\strut

\begin{proof}
By defintion, we have that

\strut

$\Phi _{0}(m_{1},N_{2},m_{2},N_{3},m_{3})$

\strut

$=E\,(\left( \underset{\left| w\right| =m_{1},\,End(w)=p_{m_{1}}^{E}}{\sum }%
w\right) (p_{11}...p_{1N_{2}})\left( \underset{\left| w\right|
=m_{2},\,\,Init(w)=p_{m_{2}}^{I},\,End(w)=p_{m_{2}}^{E}}{\sum }w\right) $

$\ \ \ \ \ \ \ \ \ \ \ \ \ \ \ \ \ \ \ \ \ \ \ \ \ \ \ \ \ \ \ \ \ \ \ \ \ \
\ \ \ \ \ \ \ \ \ \ \ \ \ \ \ \ \ \ \cdot (p_{21}...p_{2N_{3}})\left( 
\underset{\left| w\right| =m_{3},\,\,Init(w)=p_{m_{3}}^{I}}{\sum }w\right)
). $

\strut

Under the hypothesis, since $p_{11}...p_{1N_{2}}\ntriangleleft \,h^{d_{2}}$
or $p_{21}...p_{2N_{3}}\ntriangleleft \,h^{d_{3}},$ it vanishs.
\end{proof}

\strut

From now assume that

\strut

(A) \ \ \ \ \ \ \ \ \ \ \ \ \ \ \ \ \ \ \ \ \ \ \ \ \ \ \ $%
p_{11}...p_{1N_{2}}\lhd \,h^{d_{2}}$ \ \ and \ \ $p_{21}...p_{2N_{3}}\lhd
\,h^{d_{3}},$

\strut

where $N_{2},N_{3}\in \Bbb{N}$ and $d_{2},d_{3}\in \Bbb{Z}.$ Observe that if
(A) is satisfied, then

$\strut $

$\left( \underset{\left| w\right| =m_{1},\,End(w)=p_{m_{1}}^{E}}{\sum }%
w\right) (p_{11}...p_{1N_{2}})\left( \underset{\left| w\right|
=m_{2},\,\,Init(w)=p_{m_{2}}^{I},\,End(w)=p_{m_{2}}^{E}}{\sum }w\right) $

$\ \ \ \ \ \ \ \ \ \ \ \ \ \ \ \ \ \ \ \ \ \ \ \ \ \ \ \ \ \ \ \ \ \ \ \ \ \
\ \ \ \ \ \ \ \ \ \ \ \ \ \ \cdot (p_{21}...p_{2N_{3}})\left( \underset{%
\left| w\right| =m_{3},\,\,Init(w)=p_{m_{3}}^{I}}{\sum }w\right) $

\strut

$\ =(\sum_{j_{2}=0}^{N_{2}}\left( \underset{\left| w\right|
=m_{1},\,End(w)=p_{m_{1}}^{E}}{\sum }w\right) (p_{11}...p_{1j_{2}})$

$\ \ \ \ \ \ \ \ \ \ \ \ \ \ \ \ \ \ \ \ \ \cdot \
(p_{1(j_{2}+1)}...p_{1N_{2}})\left( \underset{\left| w\right|
=m_{2},\,\,Init(w)=p_{m_{2}}^{I},\,End(w)=p_{m_{2}}^{E}}{\sum }w\right) )$

$\ \ \ \ \ \ \ \ \ \ \ \ \ \ \ \ \ \ \ \ \ \cdot (p_{21}...p_{2N_{3}})\left( 
\underset{\left| w\right| =m_{3},\,\,Init(w)=p_{m_{3}}^{I}}{\sum }w\right) $

\strut

(3.3.2.1)

\strut

$\ =\sum_{j_{2}=0}^{N_{2}-1}\sum_{j_{3}=0}^{N_{3}}\left( \underset{\left|
w\right| =m_{1},\,End(w)=p_{m_{1}}^{E}}{\sum }w\right) (p_{11}...p_{1j_{2}})$

$\ \ \ \ \ \ \ \ \ \ \ \ \ \ \ \ \ \ \ \ \ \cdot \
(p_{1(j_{2}+1)}...p_{1N_{2}})\left( \underset{\left| w\right|
=m_{2},\,\,Init(w)=p_{m_{2}}^{I},\,End(w)=p_{m_{2}}^{E}}{\sum }w\right) $

$\ \ \ \ \ \ \ \ \ \ \ \ \ \ \ \ \ \ \ \ \ \cdot
(p_{21}...p_{2j_{3}})(p_{j_{3}+1}...p_{N_{3}})\left( \underset{\left|
w\right| =m_{3},\,\,Init(w)=p_{m_{3}}^{I}}{\sum }w\right) ,$

\strut

where $p_{10}=e$ and $p_{20}=e.$ Above, the case [$j_{2}=0$] (resp. [$%
j_{2}=N_{2}$]) means the case when there is no cancellation (resp. full
cancellation) for $\left( \underset{\left| w\right|
=m_{1},\,End(w)=p_{m_{1}}^{E}}{\sum }w\right) $ and $p_{11}...p_{1N_{2}}.$
Similarly, the case [$j_{3}=0$] (resp. [$j_{3}=N_{3}$]) means the case when
there is no cancellation (resp. full cancellation) for $p_{21}...p_{2N_{3}}$
and $\left( \underset{\left| w\right| =m_{3},\,\,Init(w)=p_{m_{3}}^{I}}{\sum 
}w\right) .$

\strut

Consider the summand,

\strut

$\ \ \ \ \ \ \ \ S_{m_{1},(j_{2},N_{2}),m_{2},(j_{3},N_{3}),m_{3}}=\left( 
\underset{\left| w\right| =m_{1},\,End(w)=p_{m_{1}}^{E}}{\sum }w\right)
(p_{11}...p_{1j_{2}})$

$\ \ \ \ \ \ \ \ \ \ \ \ \ \ \ \ \ \ \ \ \ \ \ \ \cdot \
(p_{1(j_{2}+1)}...p_{1N_{2}})\left( \underset{\left| w\right|
=m_{2},\,\,Init(w)=p_{m_{2}}^{I},\,End(w)=p_{m_{2}}^{E}}{\sum }w\right) $

$\ \ \ \ \ \ \ \ \ \ \ \ \ \ \ \ \ \ \ \ \ \ \ \ \cdot
(p_{21}...p_{2j_{3}})(p_{2(j_{3}+1)}...p_{2N_{3}})\left( \underset{\left|
w\right| =m_{3},\,\,Init(w)=p_{m_{3}}^{I}}{\sum }w\right) .$

\strut

\begin{lemma}
If we have the assumption (A), then the formula (3.3.2.1) contains the
following $h$-terms ;

\strut 

$\left( F_{p_{m_{1}}^{E}]\_<p_{11}...p_{1j_{2}}>}\left(
m_{1},(j_{2}-1)\right) \right) $

$\ \ \cdot \left(
F_{<p_{1}(j_{2}+1)...p_{1N_{2}}>\_[p_{m_{2}}^{I}...p_{m_{2}}^{E}]%
\_<p_{21},...,p_{2j_{3}}>}\left( (N_{2}-(j_{2}+1)),m_{2},(j_{3}-1)\right)
\right) $

$\ \ \cdot \left( F_{<p_{2(j_{3}+1)}...p_{2N_{3}}>\_[p_{m_{3}}^{I}}\left(
(N_{3}-(j_{3}+1)),m_{3}\right) \right) $,

\strut 

where

\strut 

$F_{<p_{1}...p_{n}>\_[p^{I}\_p^{E}]\_<q_{1}...q_{m}>}(n,k,m)$

\strut 

\begin{center}
$\overset{def}{=}E\left( (p_{1}...p_{n})\left( \underset{\left| w\right|
=k,\,Init(w)=p^{I},End(w)=p^{E}}{\sum }w\right) (q_{1}...q_{m})\right) .$
\end{center}
\end{lemma}

\strut

\begin{proof}
To find the $h$-terms in the given summand (3.3.2.1) is equivalent to compute

\strut

\begin{center}
$E(S_{m_{1},(j_{2},N_{2}),m_{2},(j_{3},N_{3}),m_{3}}).$
\end{center}

\strut Then

\strut

$E\left( S_{m_{1},(j_{2},N_{2}),m_{2},(j_{3},N_{3}),m_{3}}\right) $

\strut \strut

$\ \ \ =\left( F_{p_{m_{1}}^{E}]\_<p_{11}...p_{1j_{2}}>}\left(
m_{1},(j_{2}-1)\right) \right) $

\ \ $\ \ \cdot \left(
F_{<p_{1}(j_{2}+1)...p_{1N_{2}}>\_[p_{m_{2}}^{I}...p_{m_{2}}^{E}]%
\_<p_{21},...,p_{2j_{3}}>}\left( (N_{2}-(j_{2}+1)),m_{2},(j_{3}-1)\right)
\right) $

\ \ $\ \ \cdot \left(
F_{<p_{2(j_{3}+1)}...p_{2N_{3}}>\_[p_{m_{3}}^{I}}\left(
(N_{3}-(j_{3}+1)),m_{3}\right) \right) .$\strut
\end{proof}

\strut

So, we can conclude that ;

\strut

\begin{theorem}
Suppose that (A) is satisfied. Then

\strut 

\begin{center}
$\Phi
_{0}(m_{1},N_{2},m_{2},N_{3},m_{3})=\sum_{j_{2}=0}^{N_{2}-1}%
\sum_{j_{3}=0}^{N_{3}-1}E\left(
S_{m_{1},(j_{2},N_{2}),m_{2},(j_{3},N_{3}),m_{3}}\right) .$
\end{center}

$\square $
\end{theorem}

\strut

Now, we will apply the above theorem to our case ;

\strut

\begin{theorem}
Let $m_{1},m_{2},m_{3}\in \Bbb{N}$ and $d_{2},d_{3}\in \Bbb{Z}$. Then

\strut 

$E\left( X_{m_{1}}h^{d_{2}}X_{m_{2}}h^{d_{3}}X_{m_{3}}\right) =\underset{%
p_{1}^{E},\,p_{2}^{I}\in \{a,b,a^{-1},b^{-1}\}}{\sum }%
\sum_{j_{2}=0}^{4d_{2}}\sum_{j_{3}=0}^{4d_{3}}$

\strut 

$\ \ \ \ \ (F_{p_{m_{1}}^{E}]\_<p_{21}...p_{2j_{2}}>}\left(
m_{1},(j_{2}-1)\right) $

\strut 

$\ \ \ \ \ \cdot
F_{<p_{2(j_{2}+1)}...p_{2(4d_{2})}>\_[p_{m_{2}}^{I}...p_{m_{2}}^{E}]%
\_<p_{31}...p_{3j_{3}}>}\left( (4d_{2}-(j_{2}+1)),m_{2},(j_{3}-1)\right) $

\strut 

$\ \ \ \ \ \ \cdot
F_{<p_{3(j_{3}+1)}...p_{3(4d_{3})}>\_[p_{m_{3}}^{I}}\left(
(4d_{3}-(j_{3}+1)),m_{3}\right) ),$

\strut 

where

\strut 

\begin{center}
$h^{d_{2}}=p_{21}...p_{2(4d_{2})}$ \ \ and \ \ $%
h^{d_{3}}=p_{31}...p_{3(4d_{3})}.$
\end{center}
\end{theorem}

\strut

Based on the above results, we can extend our interests to the general $%
E(X_{m_{1}}h^{d_{2}}X_{m_{2}}...h^{d_{n}}X_{m_{n}})$-case.

\strut

\begin{definition}
Let $n\in \Bbb{N}$, $m_{1},...,m_{n}\in \Bbb{N}$ and $N_{2},...,N_{n}\in 
\Bbb{N}\cup \{0\}.$ Define a map $\Phi $ from $\underset{(2n-1)-times}{%
\underbrace{\Bbb{N\times }.....\Bbb{\times N}}}$ \ to $B$ by

\strut 

$\Phi \left( m_{1},N_{2},m_{2},...,N_{n},m_{n}\right) $

\strut 

$\ :=\underset{%
p_{m_{1}}^{E},p_{m_{2}}^{I},p_{m_{2}}^{E},p_{m_{3}}^{I},p_{m_{3}}^{E},...,p_{m_{n-1}}^{I},p_{m_{n-1}}^{E},p_{m_{n}}^{I}\in \{a,b,a^{-1},b^{-1}\}%
}{\sum }$

\strut 

$\ \ \ \ \
F_{p_{m_{1}}^{E}]\_<p_{21}...p_{2N_{2}}>\_[p_{m_{2}}^{I}...p_{m_{2}}^{E}]%
\_<p_{31}...p_{3N_{3}}>\_[p_{m_{3}}^{I}...p_{m_{3}}^{E}]\_...%
\_<p_{n1}...p_{nN_{n}}>\_[p_{m_{n}}^{I}}$

$\ \ \ \ \ \ \ \ \ \ \ \ \ \ \ \ \ \ \ \ \ \ \ \ \ \ \ \ \ \ \ \ \ \ \ \ \ \
\ \ \ \ \ \ \ \ \ \ \ \ \ \ \ \ \ \ \ \ \ \ \ \ \ \ \ \
(m_{1},N_{2},m_{2},...N_{n},m_{n})$

\strut 

$\ =E(\left( \underset{\left| w\right| =m_{1}}{\sum }w\right)
(p_{21}...p_{2N_{2}})\left( \underset{\left| w\right| =m_{2}}{\sum }w\right) 
$

$\ \ \ \ \ \ \ \ \ \ \ \ \ \ \ \ \cdot (p_{31}...p_{3N_{3}})\left( \underset{%
\left| w\right| =m_{3}}{\sum }w\right) (p_{41}...p_{4N_{4}})$

$\ \ \ \ \ \ \ \ \ \ \ \ \ \ \ \ \ \cdot .....(p_{n1}...p_{nN_{n}})\left( 
\underset{\left| w\right| =m_{n}}{\sum }w\right) ),$

\strut 

where $p_{ij},\,p_{m_{k}}^{I},\,p_{m_{k}}^{E}\in \{a,b,a^{-1},b^{-1}\}.$
\end{definition}

\strut \strut

For the fixed $p_{m_{1}}^{E},$ $p_{m_{n}}^{I},\,\ p_{m_{j}}^{I},$ $%
p_{m_{j}}^{E}\in \{a,b,a^{-1},b^{-1}\},$ for $\ j=2,...,n-1,$ define

\strut

$S_{m_{1},\,(j_{2},N_{2}),m_{2},(j_{3},N_{3}),m_{3},...,(j_{n},N_{n}),m_{n}}$

\strut

$\ \ =\left( \underset{\left| w\right| =m_{1},\,End(w)=p_{m_{1}}^{E}}{\sum }%
w\right) (p_{21}...p_{2j_{2}})$

$\ \ \ \ \ \ \ \ \ \ \ \ \cdot \ (p_{2(j_{2}+1)}...p_{2N_{2}})\left( 
\underset{\left| w\right|
=m_{2},\,\,Init(w)=p_{m_{2}}^{I},\,End(w)=p_{m_{2}}^{E}}{\sum }w\right) $

$\ \ \ \ \ \ \ \ \ \ \ \ \cdot
(p_{31}...p_{3j_{3}})(p_{3(j_{3}+1)}...p_{3N_{3}})\left( \underset{\left|
w\right| =m_{3},\,\,Init(w)=p_{m_{3}}^{I},End(w)=p_{m_{3}}^{E}}{\sum }%
w\right) $

$\ \ \ \ \ \ \ \ \ \ \ \ \cdot
(p_{41}...p_{4j_{4}})(p_{4(j_{4}+1)}...p_{4N_{4}})\left( \underset{\left|
w\right| =m_{4},\,\,Init(w)=p_{m_{4}}^{I},End(w)=p_{m_{3}}^{E}}{\sum }%
w\right) $

$\ \ \ \ \ \ \ \ \ \ \ \ \cdot \cdot \cdot \cdot $

\strut $\ \ \ \ \ \ \ \ \ \ \ \ \cdot \cdot \cdot \cdot $

$\ \ \ \ \ \ \ \ \ \ \ \ \cdot
(p_{n1}...p_{nj_{n}})(p_{n(j_{n}+1)}...p_{nN_{n}})\left( \underset{\left|
w\right| =m_{n},\,\,Init(w)=p_{m_{n}}^{I}}{\sum }w\right) $

\strut \strut \strut

Also, similar to the former discussion, we will assume that

\strut

(AA) \ \ \ \ \ \ \ \ \ \ \ \ \ \ \ \ \ \ \ \ \ \ \ \ \ \ \ \ \ \ \ \ $%
p_{k_{1}}...p_{kN_{k}}\lhd \,h^{d_{k}},$ for \textbf{all} $k=1,...,n,$

\strut

where $d_{k}\in \Bbb{Z}.$

\strut

\begin{theorem}
Suppose that the condition (AA) is satisfied. Then

\strut 

$E\left(
S_{m_{1},\,(j_{2},N_{2}),m_{2},(j_{3},N_{3}),m_{3},...,(j_{n},N_{n}),m_{n}}%
\right) $

\strut 

$=\left( F_{p_{m_{1}}^{E}]\_<p_{21}...p_{2j_{2}}>}(m_{1},j_{2}-1)\right) $

\strut 

$\ \ \cdot \left(
\prod_{k=2}^{n-1}F_{<p_{kj_{k+1})}...p_{kN_{k}}>%
\_[p_{m_{k}}^{I}...p_{m_{k}}^{E}]%
\_<p_{(k+1)1}...p_{(k+1)j_{k+1}}>}(N_{k}-(j_{k}+1),m_{k},j_{k+1}-1)\right) $

\strut 

$\ \ \cdot \left(
F_{<p_{n(j_{n}+1)}...p_{nN_{n}}>\_p_{m_{n}}^{I}}(N_{n}-(j_{n}+1),m_{n})%
\right) .$
\end{theorem}

\strut

\begin{proof}
If (AA) is satisfied, then

\strut

$E\left(
S_{m_{1},\,(j_{2},N_{2}),m_{2},(j_{3},N_{3}),m_{3},...,(j_{n},N_{n}),m_{n}}%
\right) $

\strut

$=E(\left( \underset{\left| w\right| =m_{1},\,End(w)=p_{m_{1}}^{E}}{\sum }%
w\right) (p_{21}...p_{2j_{2}})$

$\ \ \ \ \ \ \ \ \ \ \ \ \ \cdot \ (p_{2(j_{2}+1)}...p_{2N_{2}})\left( 
\underset{\left| w\right|
=m_{2},\,\,Init(w)=p_{m_{2}}^{I},\,End(w)=p_{m_{2}}^{E}}{\sum }w\right) $

$\ \ \ \ \ \ \ \ \ \ \ \ \ \ \cdot
(p_{31}...p_{3j_{3}})(p_{3(j_{3}+1)}...p_{3N_{3}})\left( \underset{\left|
w\right| =m_{3},\,\,Init(w)=p_{m_{3}}^{I},End(w)=p_{m_{3}}^{E}}{\sum }%
w\right) $

$\ \ \ \ \ \ \ \ \ \ \ \ \ \ \cdot
(p_{41}...p_{4j_{4}})(p_{4(j_{4}+1)}...p_{4N_{4}})\left( \underset{\left|
w\right| =m_{4},\,\,Init(w)=p_{m_{4}}^{I},\,End(w)=p_{m_{4}}^{E}}{\sum }%
w\right) $

$\ \ \ \ \ \ \ \ \ \ \ \ \ \ \cdot \cdot \cdot \cdot $

$\ \ \ \ \ \ \ \ \ \ \ \ \ \ \cdot \cdot \cdot \cdot $

$\ \ \ \ \ \ \ \ \ \ \ \ \ \ \cdot
(p_{n1}...p_{nj_{n}})(p_{n(j_{n}+1)...p_{nN_{n}}})\left( \underset{\left|
w\right| =m_{n},\,\,Init(w)=p_{m_{n}}^{I}}{\sum }w\right) )$

\strut

$=\left( F_{p_{m_{1}}^{E}]\_<p_{21}...p_{2j_{2}}>}(m_{1},j_{2}-1)\right) $

$\cdot \left(
F_{<p_{2(j_{2}+1)}...p_{2N_{2}}>\_[p_{m_{2}}^{I}...p_{m_{2}}^{E}]%
\_<p_{31}...p_{3j_{3}}>}(N_{2}-(j_{2}+1),m_{2},j_{3}-1)\right) $

$\cdot \left(
F_{<p_{3j_{3+1})}...p_{3N_{3}}>\_[p_{m_{3}}^{I}...p_{m_{3}}^{E}]%
\_<p_{41}...p_{4j_{4}}>}(N_{2}-(j_{2}+1),m_{2},j_{3}-1)\right) $

$\cdot \cdot \cdot $

$\cdot \cdot \cdot $

$\cdot \left(
F_{<p_{(n-1)j_{n-1+1})}...p_{(n-1)N_{(n-1)}}>%
\_[p_{m_{n-1}}^{I}...p_{m_{n-1}}^{E}]%
\_<p_{(n)1}...p_{nj_{n1}}>}(N_{n-1}-(j_{n-1}+1),m_{n-1},j_{n}-1)\right) $

$\cdot \left(
F_{<p_{n(j_{n}+1)}...p_{nN_{n}}>\_p_{m_{n}}^{I}}(N_{n}-(j_{n}+1),m_{n})%
\right) $

\strut

$=\left( F_{p_{m_{1}}^{E}]\_<p_{21}...p_{2j_{2}}>}(m_{1},j_{2}-1)\right) $

\strut

$\ \ \cdot \left(
\prod_{k=2}^{n-1}F_{<p_{kj_{k+1})}...p_{kN_{k}}>%
\_[p_{m_{k}}^{I}...p_{m_{k}}^{E}]%
\_<p_{(k+1)1}...p_{(k+1)j_{k+1}}>}(N_{k}-(j_{k}+1),m_{k},j_{k+1}-1)\right) $

\strut

$\ \ \cdot \left(
F_{<p_{n(j_{n}+1)}...p_{nN_{n}}>\_p_{m_{n}}^{I}}(N_{n}-(j_{n}+1),m_{n})%
\right) .$
\end{proof}

\strut

By using the above theorem, we have that ;

\strut \strut \strut

\begin{theorem}
Suppose that the condition (AA) is satisfied. Then

\strut 

$\Phi \left( m_{1},N_{2},m_{2},...,N_{n},m_{n}\right) $

\strut 

$\ =\underset{p_{m_{1}}^{E},p_{m_{n}}^{I},p_{m_{j}}^{I},p_{m_{j}}^{E}\in
\{a,b,a^{-1},b^{-1}\},\,\,j=2,...,n-1}{\sum }\sum_{j_{2}=0}^{N_{2}}%
\sum_{j_{3}=0}^{N_{3}}\cdot \cdot \cdot \sum_{j_{n}=0}^{N_{n}}$

\begin{center}
$E\left(
S_{m_{1},\,(j_{2},N_{2}),m_{2},(j_{3},N_{3}),m_{3},...,(j_{n},N_{n}),m_{n}}%
\right) .$
\end{center}

$\square $
\end{theorem}

\strut

Applying the above two theorems, we have that ;

\strut

\begin{theorem}
Let $n\in \Bbb{N},$ $m_{1},...,m_{n}\in \Bbb{N}$ and $d_{2},...,d_{n}\in 
\Bbb{Z}.$ Then

\strut 

$\ E\left(
X_{m_{1}}h^{d_{2}}X_{m_{2}}h^{d_{3}}X_{m_{3}}...h^{d_{n}}X_{m_{n}}\right) $

\begin{center}
\strut 

$=\Phi \left( m_{1},4d_{2},m_{3},4d_{3},...,4d_{n},m_{n}\right) $,
\end{center}

\strut 

by the triangular relation ''$\,\lhd \,$''. $\square $
\end{theorem}

\strut

\begin{theorem}
Let $n\in \Bbb{N},$ $m_{1},...,m_{n}\in \Bbb{N}$ and $d_{2},...,d_{n}\in 
\Bbb{N}\cup \{0\}.$ Then

\strut 

$E\left(
X_{m_{1}}(h^{d_{2}}+h^{-d_{2}})X_{m_{2}}(h^{d_{3}}+h^{-d_{3}})X_{m_{3}}...(h^{d_{n}}+h^{d_{n}})X_{m_{n}}\right) 
$

\begin{center}
\strut 
\end{center}

$\ \ =\underset{r_{2}\in \{\pm d_{2}\},\,r_{3}\in \{\pm d_{3}\},...,r_{n}\in
\{\pm d_{n}\}}{\sum }\Phi \left(
m_{1},4r_{2},m_{3},4r_{3},...,4r_{n},m_{n}\right) .$
\end{theorem}

\strut

\strut

\strut

\subsection{Computing Trivial Cumulants of $x$}

\strut

\strut

\strut

Let $k\in \Bbb{N}.$ Consider the $2k$-th trivial cumulants of $%
x=a+b+a^{-1}+b^{-1}$,

\strut

$K_{2k}^{t}\left( \underset{2k-times}{\underbrace{x,......,x}}\right) =%
\underset{\pi \in NC^{(even)}(2k)}{\sum }\widehat{E}(\pi )\left( x\otimes
...\otimes x\right) \mu (\pi ,1_{2k})$

\strut

$=\underset{l_{1},...,l_{p}\in 2\Bbb{N},\,\,l_{1}+...+l_{p}=2k}{\sum }\,%
\underset{\pi \in NC_{l_{1},...,l_{p}}(2k)}{\sum }\widehat{E}(\pi )\left(
x\otimes ...\otimes x\right) \mu (\pi ,1_{2k}),$

\strut

where

\strut

\begin{center}
$NC_{l_{1},...,l_{p}}(2k)=\{\pi \in NC^{(even)}(2k):V\in \pi \Leftrightarrow
\left| V\right| =l_{j},\,j=1,...,p\}.$
\end{center}

\strut

For example,

\begin{center}
$NC^{(even)}(8)=NC_{2,2,2,2}(8)\cup NC_{2,2,4}(8)\cup NC_{2,6}(8)\cup
NC_{4,4}(8)\cup NC_{8}(8).$
\end{center}

\strut

\strut

\begin{lemma}
(See Section 3.2.2) Fix $k\in \Bbb{N}.$ Let $h=aba^{-1}b^{-1}\in
A_{1}*_{B}A_{2}$ with $h^{0}=e.$

\strut 

(1) If $4\mid 2k,$ then

\strut 

\begin{center}
$E\left( x^{2k}\right) =\sum_{j=0}^{\frac{k-2}{2}}p_{2k-4j}^{2k}\left( h^{%
\frac{k}{2}-j}+h^{-(\frac{k}{2}-j)}\right) +p_{0}^{2k}h^{0},$

\strut 
\end{center}

where $p_{0}^{4}=28$ and $p_{2k}^{2k}=1$.

\strut 

(2) If $4\nmid 2k$ and if there are $X_{4l_{1}},...,X_{4l_{p}}$ terms in $%
x^{2k},$ then

\strut 

\begin{center}
$E(x^{2k})=\sum_{j=0}^{\frac{k-3}{2}}p_{(2k-2)-4j}^{2k}\left( h^{\frac{k-1}{2%
}-2j}+h^{-(\frac{k-1}{2}-2j)}\right) +p_{0}^{2k}h^{0},$

\strut 
\end{center}

where $p_{0}^{2}=4.$ \ $\square $
\end{lemma}

\strut

\begin{definition}
Let $n\in \Bbb{N}$ and let $\pi \in NC(n).$ Let $%
V=(v_{1},...,v_{k}),W=(w_{1},...,w_{l})\in \pi $ and assume that there
exists $j\in \{1,...,k\}$ such that

\strut 

\begin{center}
$1\leq v_{1}<...<v_{j}<w_{1}<...<w_{l}<v_{j+1}<...<v_{k}\leq n,$
\end{center}

\strut 

in $\{1,...,n\}.$ Then we say that the block $W$ is a subblock of the block $%
V.$
\end{definition}

\strut

Suppose that $V\in \pi (i)$ is an inner block of a partition $\pi $ with its
outer block $W\in \pi (o).$ Then $V$ is a subblock of $W.$\strut

\strut

\strut

\begin{definition}
Let $\pi \in NC_{l_{1},...,l_{p}}(2k).$ Let $V^{o}\in \pi (o)$ and let $%
V,W\in \pi $ be subblocks in $V^{o}.$ We say that $V$ is inner in $W,$ if $V$
is a subblok of $W.$ It has the following pictorial expression ;

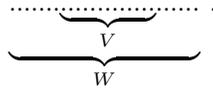
\strut 

\begin{center}
$\underset{W}{\underbrace{.......\underset{V}{\underbrace{.............}}%
......}}$ .
\end{center}

\strut 

Notice that it does not mean $V\in \pi (i)$ with its outer block $W\in \pi
(o).$ i.e we can have the following pictorial expression of subblocks $%
V^{\prime },V,W,W^{\prime }$ in $V^{o}\in \pi (o)$ ;

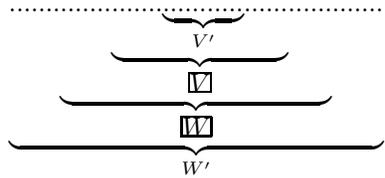
\strut 

\begin{center}
$\underset{W^{\prime }}{\underbrace{.......\underset{\frame{$W$}}{%
\underbrace{.......\underset{\frame{$V$}}{\underbrace{.......\underset{%
V^{\prime }}{\underbrace{...........}}......}}......}}.......}},$
\end{center}

\strut 

where $V^{\prime }$ and $W^{\prime }$ are other blocks in $\pi $ which are
subblocks of $V^{o}\in \pi (o).$ We also say that $V$\textbf{\ is a deepest
subblock} if there is no subblcok $V^{\prime }$ which is inner in $V.$ We
define that $V\in \pi (o)$ containing no inner blocks is also a deepest
block. We will express pictorially,

\strut 

\begin{center}
$\underset{W}{\underbrace{.......\underset{V}{\underbrace{.............}}%
......}}$ \ \ \ \ \ by \ \ \ \ \ $\underset{W}{\underbrace{\frame{$m_{1}$}%
\,\,\underset{V}{\underbrace{.............}}\,\frame{$m_{2}$}}}$ ,
\end{center}

\strut 

where $m_{1}=\left| W\right| ^{(1)}$ and $m_{2}=\left| W\right| ^{(2)}$ such
that $\left| W\right| =m_{1}+m_{2}.$
\end{definition}

\strut

\begin{definition}
Suppose that $\pi \in NC_{l_{1},...,l_{p}}(2k).$ Then $\pi
=\{V_{l_{1}},...,V_{l_{p}}\}$ with $\left| V_{l_{j}}\right| =l_{j},$ $%
\forall \,j=1,...,p.$ Define a map $\Lambda $ from

\strut 

\begin{center}
$\cup _{k=1}^{\infty }\left( \underset{l_{1},...,l_{p}\in 2\Bbb{N}%
,\,\,l_{1}+...+l_{p}=2k}{\cup }NC_{l_{1},...,l_{p}}(2k)\right) $
\end{center}

to

\begin{center}
$\cup _{n=1}^{\infty }\left( \underset{n-times}{\underbrace{\Bbb{N}\times
...\times \Bbb{N}}}\right) ,$

\strut 
\end{center}

by the following rules with respect to subblocks \textbf{in each outer block}
of $\pi $ ;

\strut \strut 

(1) Let $\pi \in NC_{l_{1},...,l_{p}}(2k)$ and $V\in \pi (o)$ and let $%
V_{1},...,V_{k},W$ be subblocks in $V.$ Let $V_{1}$, $V_{2},$ ..., $V_{k}$
be deepest subblocks in $W,$ which is outer of $V_{1},...,V_{k}.$ If

\strut 

\begin{center}
$\underset{W}{\underbrace{m_{1}\,\underset{V_{1}}{\underbrace{............}}%
m_{2}\underset{V_{2}}{\underbrace{............}}m_{3}...m_{k}\underset{V_{k}%
}{\underbrace{............}}\,\,m_{k+1}}}$ \ \ \ 
\end{center}

with \ \ \ 

\begin{center}
$\left| V_{1}\right| =n_{1},...,\left| V_{k}\right| =n_{k}$
\end{center}

\strut 

then the restriction of $\Lambda $ to $\left( V_{1},...,V_{k},W\right) $ is

\strut 

\begin{center}
$\Lambda \mid _{(V_{1},...,V_{k},W)}(\pi )=\left(
m_{1},[n_{1}],m_{2},[n_{2}],m_{3},...,m_{k},[n_{k}],m_{k+1}\right) .$
\end{center}

\strut 

(2) Let $V\in \pi (o),\,W\in \pi (i)$ be given as in the step (1). Suppose
that the subblock of $V,$ $W_{1}\in \pi $ is outer of $W$ such that

\strut 

\begin{center}
$\underset{W_{1}}{\underbrace{p_{1}\,\,\underset{W}{\underbrace{..........}}%
\,\,p_{2}}}$ \ \ \ \ \ \ \ with \ \ \ \ $\left| W_{1}\right| =p_{1}+p_{2}.$
\end{center}

\strut 

Then the restriction of $\Lambda $ to $(W,W_{1})$ is

\strut 

\begin{center}
$\Lambda \mid _{(W,\,W_{1})}(\pi )=$ \ $\left( p_{1},\,[\Lambda \mid
_{(V_{1},...,V_{k},W)}(\pi )],\,p_{2}\right) .$
\end{center}

\strut 

If $W_{2}$ is outer of $W_{1},$ then, similarly, we can define $\Lambda \mid
_{(W_{1},W_{2})}(\pi ),$ inductively. This is the insertion property of the
map $\Lambda .$ Suppose that there are restrictions of the map $\Lambda ,$ $%
\Lambda \mid _{(V_{1}^{1},...,V_{k_{1}}^{1},W^{1})}(\pi ),...,\Lambda \mid
_{(V_{1}^{l},...,V_{k}^{l},W^{l})}(\pi )$ \ given as in (1) and assume that $%
W^{1},...,W^{l}$ are inner in $W_{1}$ and $W_{1}$ is outer of them. Then
similar to (1), we have that

\strut 

$\Lambda _{(W^{1},...,W^{l},W_{1})}(\pi )$

\strut 

\begin{center}
$=\left( q_{1},[\Lambda \mid _{(V_{1}^{1},...,V_{k_{1}}^{1},W^{1})}(\pi
)],q_{2},...,q_{l},[\Lambda \mid _{(V_{1}^{l},...,V_{k}^{l},W^{l})}(\pi
)],q_{l+1}\right) $
\end{center}

\strut 

(3) Let $V^{o}\in \pi (o)$ be an outer block having the above insertion
property (1) and (2). Then the restriction of the map $\Lambda $ to $V^{o}$, 
$\Lambda \mid _{V^{o}}(\pi ),$ is defined like in (2), inductively.

\strut 

(4) Let $\pi (o)=\{V_{1}^{o},...,V_{t}^{o}\}.$ Define the map $\Lambda $ for 
$\pi \in NC_{l_{1},...,l_{p}}(2k)$ by

\strut 

\begin{center}
$\Lambda (\pi )=\Lambda \mid _{V_{1}^{o}}(\pi )\times ...\times \Lambda \mid
_{V_{t}^{o}}(\pi ),$
\end{center}

\strut 

where the Cartesian product ``$\times $'' is set-theoratically determined.
Recall that each $\Lambda \mid _{V_{j}^{o}}(\pi )$ has the insertion
property as in (2) and (3).

The map $\Lambda ,$ on even noncrossing partitions is called the
``(partition-dependent) Numbering'' map. Remark that the image of this
numbering map contains the number with rectangular bracket $[.]$ which
represents the length of the deepest blocks in the outer blocks of
partitions.
\end{definition}

\strut \strut \strut

\begin{example}
Let $\pi \in NC_{2,4,4,4,4,6}(24)$ be given as follows ;

\strut 

$\pi =\,\underset{W_{1}}{\underbrace{\circ \,\,\underset{V_{1}}{\underbrace{%
\circ \circ \circ \circ }}\,\underset{V_{2}}{\underbrace{\circ \circ \circ
\circ }}\circ \circ \circ }}\,\underset{W_{2}}{\underbrace{\,\circ \circ \,%
\underset{V_{4}}{\underbrace{\circ \underset{V_{3}}{\underbrace{\circ \circ
\circ \circ \circ \circ }}\circ }}\,\circ \circ }}.$

\strut 

Then we can reexpress it by

\strut 

$\pi =\,\underset{W_{1}}{\underbrace{\frame{$1$}\,\,\underset{V_{1}}{%
\underbrace{\circ \circ \circ \circ }}\,\underset{V_{2}}{\underbrace{\circ
\circ \circ \circ }}\,\frame{$3$}}}\,\underset{W_{2}}{\underbrace{\,\frame{$2
$}\,\underset{V_{4}}{\underbrace{\frame{$1$}\underset{V_{3}}{\underbrace{%
\circ \circ \circ \circ \circ \circ }}\frame{$1$}}}\,\frame{$2$}}}$

\strut 

and hence, by the map $\Lambda $, we have that

\strut 

\begin{center}
$\Lambda (\pi )=(1,[4],[4],3)\times (2,(1,[6],1),2).$
\end{center}
\end{example}

\strut \strut

\begin{definition}
Fix $k\in \Bbb{N}$ and $\pi \in NC_{l_{1},...,l_{p}}(2k).$ Let $(A,E)$ be a
NCPSpace over $B$ and let $x_{0}\in (A,E)$ be a $B$-even random variable.
Define a map

\strut 

\begin{center}
$\Psi _{x_{0}}:\cup _{k=1}^{\infty }\left( \underset{l_{1},...,l_{p}\in
2N,\,\,l_{1}+...+l_{p}=2k}{\cup }NC_{l_{1},...,l_{p}}(2k)\right) \rightarrow
B$
\end{center}

by

\strut \strut 

\begin{center}
$\Psi _{x_{0}}(\pi )=\Phi _{x_{0}}\circ \Lambda (\pi )\overset{def}{=}%
\underset{V\in \pi (o)}{\prod }\Phi _{x_{0}}\left( \Lambda (V)\right) .$
\end{center}

\strut 

Remark that $\Lambda (V)$ satisfies the insertion property and hence $\Phi
_{x_{0}}\left( \Lambda (V)\right) $ is also defined by the insertion
property. And hence $\Psi _{x_{0}}$ is defined by the insertion property.
\end{definition}

\strut

\begin{example}
Let $\pi \in NC_{2,4,4,4,4,6}(24)$ be given as follows ;

\strut \strut 

\begin{center}
$\pi =\,\underset{W_{1}}{\underbrace{\frame{$1$}\,\,\underset{V_{1}}{%
\underbrace{\circ \circ \circ \circ }}\,\underset{V_{2}}{\underbrace{\circ
\circ \circ \circ }}\,\frame{$3$}}}\,\underset{W_{2}}{\underbrace{\,\frame{$2
$}\,\underset{V_{4}}{\underbrace{\frame{$1$}\underset{V_{3}}{\underbrace{%
\circ \circ \circ \circ \circ \circ }}\frame{$1$}}}\,\frame{$2$}}}$
\end{center}

\strut 

and hence,

\strut 

$\Psi _{x_{0}}\left( \pi \right) =\left( \Psi _{x_{0}}(W_{1})\right) \left(
\Psi _{x_{0}}(W_{2})\right) $

\strut 

$\ \ \ \ \ \ \ \ \ \ \ =\left( \Phi _{x_{0}}(\Lambda (W_{1}))\right) \left(
\Phi _{x_{0}}(\Lambda (W_{2}))\right) $

\strut 

$\ \ \ \ \ \ \ \ \ \ \ =\left( \Phi _{x_{0}}(1,[4],[4],3)\right) \left( \Phi
_{x_{0}}(2,(1,[6],1),2)\right) $

\strut 

$\ \ \ \ \ \ \ \ \ \ \ =E\left(
x_{0}E(x_{0}^{4})E(x_{0}^{4})x_{0}^{3}\right) \cdot E\left(
x_{0}^{2}E(x_{0}E(x_{0}^{6})x_{0})x_{0}^{2}\right) .$

\strut 

for the (arbitrary) fixed $B$-valued random variable $x_{0}$ in (some)
NCPSpace over $B,$ $(A,E).$
\end{example}

\strut

Now, let's go back to our problem and observe the following ;

\strut \strut

\begin{theorem}
Let $B$-valued random variables, $x=a+b+a^{-1}+b^{-1}$ and $%
y=c+d+c^{-1}+d^{-1}$ in $\left( A_{1}*_{B}A_{2},F:=E*E\right) $ be given as
before.Then

\strut 

(1) If $R_{x+y}^{t}$ is the trivial $B$-valued R-transform of $x+y,$ then

\strut 

\begin{center}
$
\begin{array}{ll}
coef_{2k}\left( R_{x+y}^{t}\right)  & =K_{2k}^{t}\left(
(x+y),...,(x+y)\right)  \\ 
&  \\ 
& =2\underset{l_{1},...,l_{p}\in 2\Bbb{N},\,l_{1}+...+l_{p}=2k}{\sum }\,%
\underset{\pi \in NC_{l_{1},...,l_{p}}(2k)}{\sum }\mu _{\pi }\cdot \Psi
_{x}(\pi ),
\end{array}
$
\end{center}

\strut 

where $\mu _{\pi }=\mu (\pi ,1_{2k})\in \Bbb{C},$ for all $\pi \in
NC^{(even)}(2k),$ $k\in \Bbb{N}.$

\strut 

(2) If $M_{x+y}^{t}$ is the trivial $B$-valued moment series of $x+y,$ then

\strut 

\begin{center}
$
\begin{array}{ll}
coef_{2k}\left( M_{x+y}^{t}\right)  & =E\left( (x+y)^{2k}\right)  \\ 
&  \\ 
& =\,\underset{\theta \in NC^{(even)}(2k)}{\sum }\,\left( 2^{\left| \theta
\right| }\underset{\pi \in NC^{(even)}(2k),\,\pi \leq \theta }{\sum }\mu
_{\pi }^{\theta }\cdot \Psi _{x}(\pi )\right) ,
\end{array}
$
\end{center}

\strut 

where $\mu _{\pi }^{\theta }=\mu (\pi ,\theta )$, for all $k\in \Bbb{N}.$
\end{theorem}

\strut

\strut

So, to compute $B$-valued $2k$-th cumulants of $x$ (and hence to compute $B$%
-valued moments of $x$), it is sufficeint to compute $\Psi _{x}(\pi ),$ for
each $\pi \in NC^{(even)}(2k),$ for $k\in \Bbb{N}.$ \strut

\strut

\begin{quote}
\frame{\textbf{Notation}} (1) From now, for the convenience of using
notations, we will denote

\strut
\end{quote}

\begin{center}
$E(x^{2k})=\sum_{n=0}^{a(2k)}\alpha _{4n}^{2k}\left( h^{n}+h^{-n}\right) \in
B,$
\end{center}

\begin{quote}
\strut

where $\alpha _{0}^{2k}=p_{0}^{2k}$, $\alpha _{4}^{2k}=p_{4}^{2k},$ $\alpha
_{8}^{2k}=p_{8}^{2k},...,\alpha _{4n}^{2k}=p_{4n}^{2k},....$ Moreover, we
will denote

\strut
\end{quote}

\begin{center}
$(h^{0}+h^{-0})\overset{denote}{\equiv }e.$
\end{center}

\begin{quote}
\strut

Remark that the above equality ''$\equiv $'' is just mean the notation, not
equality ! Also, we will denote that

\strut
\end{quote}

\begin{center}
$a(2k)=\left\{ 
\begin{array}{lll}
\frac{k}{2} &  & \text{if }4\mid 2k \\ 
&  &  \\ 
\frac{k-1}{2} &  & \text{if }4\nmid 2k.
\end{array}
\right. $
\end{center}

\begin{quote}
$\strut $

(2) Let $m\in \Bbb{N}.$ By Section 3.2.2, we know that

\strut
\end{quote}

\begin{center}
$x^{m}=\left\{ 
\begin{array}{lll}
X_{m}+p_{m-2}^{m}X_{m-2}+...+p_{2}^{m}X_{2}+p_{0}^{m}e &  & \text{if }m\in 2%
\Bbb{N} \\ 
&  &  \\ 
X_{m}+q_{m-2}^{m}X_{m-2}+...+q_{3}^{m}X_{3}+q_{1}^{m}X_{1} &  & \text{if }%
m\in 2\Bbb{N}-1
\end{array}
\right. $
\end{center}

\begin{quote}
\strut

We will denote $x^{m},$ at once, by

\strut
\end{quote}

\begin{center}
$x^{m}=\sum_{n=0}^{m}\beta _{n}^{m}X_{n},$
\end{center}

\begin{quote}
\strut where
\end{quote}

\begin{center}
$\beta _{n}^{m}=\left\{ 
\begin{array}{lll}
p_{n}^{m} &  & \text{if }m\in 2\Bbb{N} \\ 
q_{n}^{m} &  & \text{if }m\in 2\Bbb{N}-1 \\ 
0 &  & \text{if }n\neq m-2j.
\end{array}
\right. \strut $
\end{center}

\begin{quote}
$\strut $Notice that if $n\neq m-2j,$ then $\beta _{n}^{m}=0$ in $\Bbb{C}.$
\ $\square $
\end{quote}

\strut

Recall that, there are two recurrence relations for computing $E(x^{2k})$,
the first is the case when $4\mid 2k$ \ and the second is the case when \ $%
4\nmid 2k.$ The above notation is used because we want to consider $%
E(x^{2k}),$ at once.\strut \strut i.e we want to avoid the situation where
we have to observe case-by-case.

\strut

\begin{proposition}
Let $m_{1},...,m_{n}\in \Bbb{N}$ and $d_{2},...,d_{n}\in \Bbb{N}\cup \{0\}.$
Then

\strut 

$E\left(
x^{m_{1}}(h^{d_{2}}+h^{-d_{2}})x^{m_{2}}(h^{d_{3}}+h^{-d_{3}})x^{m_{3}}...(h^{d_{n}}+h^{-d_{n}})x^{m_{n}}\right) 
$

\strut 

$=\sum_{i_{1}=0}^{m_{1}}\sum_{i_{2}=0}^{m_{2}}\sum_{i_{3}=0}^{m_{3}}\cdot
\cdot \cdot \sum_{i_{n}=0}^{m_{n}}\left( \beta _{i_{1}}^{m_{1}}\beta
_{i_{2}}^{m_{2}}\beta _{i_{3}}^{m_{3}}\cdot \cdot \cdot \beta
_{i_{n}}^{m_{n}}\right) $

\strut 

$\ \ \ \ \ \ \ \ \ \left( \underset{r_{2}\in \{\pm d_{2}\},\,r_{3}\in \{\pm
d_{3}\},...,r_{n}\in \{\pm d_{n}\}}{\sum }\Phi \left(
[m_{1},[4r_{2}],m_{3},[4r_{3}],...,[4r_{n}],m_{n}]\right) \right) .$
\end{proposition}

\strut \strut 

\strut

\begin{theorem}
Let $k\in \Bbb{N}.$ Let $l_{1},...,l_{p}\in 2\Bbb{N}$ and $%
l_{1}+...+l_{p}=2k.$ Let $\pi \in NC_{l_{1},...,l_{p}}(2k)$ and $V\in \pi .$
If $V_{1},...,V_{n}$ are deepest subblocks in $V$ which is outer of $%
V_{1},...,V_{n}$ and if we have

\strut 

\begin{center}
$\Lambda \mid _{(V_{1},...,V_{n},V)}(\pi )=\left(
m_{1},[l_{j_{2}}],m_{2},[l_{j_{3}}],m_{3},...,[l_{j_{n}}],m_{n}\right) ,$
\end{center}

\strut 

where $l_{j_{2}},...,l_{j_{n}}\in \{l_{1},...,l_{p}\}$ and $%
m_{1},...,m_{n}\in \Bbb{N}\cup \{0\},$ then the restriction of $\Psi _{x}$
to $(V_{1},...,V_{n},V)$ is

\strut 

$\Psi _{x}\mid _{(V_{1},...,V_{n},V)}(\pi
)=\sum_{n_{1}=0}^{m_{1}}\sum_{n_{2}=0}^{m_{2}}\cdot \cdot \cdot
\sum_{n_{n}=0}^{m_{n}}\sum_{i_{2}=0}^{a(l_{j_{2}})}%
\sum_{i_{3}=0}^{a(l_{j_{3}})}\cdot \cdot \cdot \sum_{i_{n}=0}^{a(l_{j_{n}})}$

\strut 

$\ \ \ \ \ \ \ \ \ \ \ \ \ \ \ \ \ \ \ \ \left( \beta _{n_{1}}^{m_{1}}\beta
_{n_{2}}^{m_{2}}\cdot \cdot \cdot \beta _{n_{n}}^{m_{n}}\right) \left(
\alpha _{4i_{2}}^{l_{j_{2}}}\alpha _{4i_{3}}^{l_{j_{3}}}\cdot \cdot \cdot
\alpha _{4i_{n}}^{l_{j_{n}}}\right) $

\strut 

$\ \ \ \ \ \ \ \ \ \ \ \ \ \ \ \ \ \ \left( \underset{r_{2}\in \{\pm
i_{2}\},\,r_{3}\in \{\pm i_{3}\},...,r_{n}\in \{\pm i_{n}\}}{\sum }\,\Phi
\left( n_{1},4r_{2},n_{3},4r_{3},...,4r_{n},n_{n}\right) \right) .$
\end{theorem}

\strut

\begin{proof}
Fix $\pi \in NC_{l_{1},...,l_{p}}(2k)$ and let $V^{o}\in \pi (o)$ and assume
that $\Lambda (V^{o})=\left(
m_{1},[l_{j_{2}}],m_{2},...,[l_{j_{n}}],m_{n}\right) ,$ where $1\leq
j_{1}<j_{2}<...<j_{n}\leq p.$ Then

\strut

$\Psi _{x}\mid _{(V_{1},...,V_{n},V)}(\pi )=\Phi _{x}\circ \Lambda \mid
_{(V_{1},...,V_{n},V)}(\pi )$

\strut

$\ =E\left(
x^{m_{1}}E(x^{l_{j_{2}}})x^{m_{2}}E(x^{l_{j_{3}}})x^{m_{3}}...E(x^{l_{j_{n}}})x^{m_{n}}\right) 
$

\strut

$\ =E(\left( \sum_{n_{1}=0}^{m_{1}}\beta _{n_{1}}^{m_{1}}X_{n_{1}}\right)
\left( \sum_{i_{2}=0}^{a(l_{j_{2}})}\alpha
_{4i_{2}}^{l_{j_{2}}}(h^{i_{2}}+h^{-i_{2}})\right) $

$\ \ \ \ \ \ \ \ \ \ \ \ \ \left( \sum_{n_{2}=0}^{m_{2}}\beta
_{n_{2}}^{m_{2}}X_{n_{2}}\right) \left( \sum_{i_{3}=0}^{a(l_{j_{3}})}\alpha
_{4i_{3}}^{l_{j_{3}}}(h^{i_{3}}+h^{-i_{3}})\right) \left(
\sum_{n_{3}=0}^{m_{3}}\beta _{n_{3}}^{m_{3}}X_{n_{3}}\right) $

$\ \ \ \ \ \ \ \ \ \ \ \ \ \ .........................\left(
\sum_{i_{n}=0}^{a(l_{j_{n}})}\alpha
_{4i_{n}}^{l_{j_{n}}}(h^{i_{n}}+h^{-i_{n}})\right) \left(
\sum_{n_{n}=0}^{m_{n}}\beta _{n_{n}}^{m_{n}}X_{n_{n}}\right) )$

\strut

$\ =\sum_{n_{1}=0}^{m_{1}}\sum_{n_{2}=0}^{m_{2}}\cdot \cdot \cdot
\sum_{n_{n}=0}^{m_{n}}\sum_{i_{2}=0}^{a(l_{j_{2}})}%
\sum_{i_{3}=0}^{a(l_{j_{3}})}\cdot \cdot \cdot \sum_{i_{n}=0}^{a(l_{j_{n}})}$

$\ \ \ \ \ \ \ \ \ \ \ \ \ \ \ \ \ \ \ \ \left( \beta _{n_{1}}^{m_{1}}\beta
_{n_{2}}^{m_{2}}\cdot \cdot \cdot \beta _{n_{n}}^{m_{n}}\right) \left(
\alpha _{4i_{2}}^{l_{j_{2}}}\alpha _{4i_{3}}^{l_{j_{3}}}\cdot \cdot \cdot
\alpha _{4i_{n}}^{l_{j_{n}}}\right) $

$\ \ \ \ \ \ \ \ \ \ \ \ \ \ \ \ \ \ \ \ E\left(
X_{n_{1}}(h^{i_{2}}+h^{-i_{2}})X_{n_{2}}(h^{i_{3}}+h^{-i_{3}})X_{n_{3}}...(h^{i_{n}}+h^{-i_{n}})X_{n_{n}}\right) 
$

\strut

$\ =\sum_{n_{1}=0}^{m_{1}}\sum_{n_{2}=0}^{m_{2}}\cdot \cdot \cdot
\sum_{n_{n}=0}^{m_{n}}\sum_{i_{2}=0}^{a(l_{j_{2}})}%
\sum_{i_{3}=0}^{a(l_{j_{3}})}\cdot \cdot \cdot \sum_{i_{n}=0}^{a(l_{j_{n}})}$

$\ \ \ \ \ \ \ \ \ \ \ \ \ \ \ \ \ \ \ \ \left( \beta _{n_{1}}^{m_{1}}\beta
_{n_{2}}^{m_{2}}\cdot \cdot \cdot \beta _{n_{n}}^{m_{n}}\right) \left(
\alpha _{4i_{2}}^{l_{j_{2}}}\alpha _{4i_{3}}^{l_{j_{3}}}\cdot \cdot \cdot
\alpha _{4i_{n}}^{l_{j_{n}}}\right) $

$\ \ \ \ \ \ \ \ \ \ \ \ \ \ \ \ \ \ \left( \underset{r_{2}\in \{\pm
i_{2}\},\,r_{3}\in \{\pm i_{3}\},...,r_{n}\in \{\pm i_{n}\}}{\sum }\,\Phi
\left( n_{1},4r_{2},n_{3},4r_{3},...,4r_{n},n_{n}\right) \right) .$
\end{proof}

\strut \strut \strut

\begin{proposition}
Let $k\in \Bbb{N}$ and let $\pi \in NC^{(even)}(2k).$ If $%
x=a+b+a^{-1}+b^{-1}\in A_{1}*_{B}A_{2}$ is our $B$-valued random variable,
then

\strut

\begin{center}
$\widehat{E}(\pi )\left( x\otimes ...\otimes x\right) =\Psi _{x}(\pi ).$
\end{center}

$\square $
\end{proposition}

\strut

The above proposition is proved by only using the definitions of
partition-dependent $B$-moments of $x$ and of $\Psi _{x_{0}}.$ \strut Now,
we will put all our information ;

\strut

\begin{theorem}
Let $k\in \Bbb{N}$ and let $x=a+b+a^{-1}+b^{-1}$ and $y=c+d+c^{-1}+d^{-1}$
be given as before. Then

\strut

$K_{2k}^{t}\left( (x+y),...,(x+y)\right) $

\strut

$\ =2\cdot \underset{l_{1},...,l_{p}\in 2\Bbb{N},\,\,l_{1}+...+l_{p}=2k}{%
\sum }\,\underset{\pi \in NC_{l_{1},...,l_{p}}(2k)}{\sum }\mu _{\pi }\,\cdot
\,\Psi _{x}(\pi )$

\strut

$\ =2\cdot \underset{l_{1},...,l_{p}\in 2\Bbb{N},\,\,l_{1}+...+l_{p}=2k}{%
\sum }\,\underset{\pi \in NC_{l_{1},...,l_{p}}(2k)}{\sum }\mu _{\pi }\,\cdot
\left( \underset{V^{o}\in \pi (o)}{\prod }\left( \Psi _{x_{0}}\mid
_{V^{o}}(\pi )\right) \right) ,$

\strut

by the previous proposition. $\square $
\end{theorem}

\strut

\strut

\begin{theorem}
Let $k\in \Bbb{N}$ and let $x=a+b+a^{-1}+b^{-1}$ and $y=c+d+c^{-1}+d^{-1}$
be given as before. Then

\strut

\begin{center}
$
\begin{array}{ll}
coef_{2k}\left( M_{x+y}^{t}\right) & =\underset{\theta \in NC^{(even)}(2k)}{%
\sum }2^{\left| \theta \right| }\cdot \,\underset{\pi \in
NC^{(even)}(2k),\,\pi \leq \theta }{\sum }\mu _{\pi }^{\theta }\cdot \Psi
_{x}(\pi ).
\end{array}
\ \ \ \ \ \ \ \ \ \ \ \ \ $
\end{center}
\end{theorem}

\strut \strut

\begin{proof}
Fix $k\in \Bbb{N}.$ Then we already observed that

\strut

$E\left( (x+y)^{2k}\right) =\underset{\theta \in NC^{(even)}(2k)}{\sum }%
\widehat{C_{E}}(\theta )\left( (x+y)\otimes ...\otimes (x+y)\right) $

\strut

where $\widehat{C}=(C^{(n)})_{n=1}^{\infty }\in I\left(
L(F_{2}),L(F_{1})\right) ^{c}$ is the cumulant multiplicative bimodule map
induced by $E:L(F_{2})\rightarrow L(F_{1})$

\strut

$\ =\underset{\theta \in NC^{(even)}(2k)}{\sum }2^{\left| \theta \right|
}\cdot \,\widehat{C}(\theta )\left( x\otimes x\otimes ...\otimes x\right) $

\strut

by the $B$-freeness of $x$ and $y$ and by the identically $B$%
-distributedness of $x$ and $y$

\strut

$\ =\underset{\theta \in NC^{(even)}(2k)}{\sum }2^{\left| \theta \right|
}\cdot \left( \underset{\pi \in NC^{(even)}(2k),\,\,\pi \leq \theta }{\sum }%
\widehat{E}(\pi )\left( x\otimes x\otimes ...\otimes x\right) \mu _{\pi
}^{\theta }\right) $

\strut

where $\mu _{\pi }^{\theta }=\mu (\pi ,\theta ).$\strut
\end{proof}

\strut

\begin{remark}
We have that

\strut

\begin{center}
$K_{2k}^{t}\left( (x+y),...,(x+y)\right) =2\cdot \underset{\pi \in
NC^{(even)}(2k)}{\sum }\mu _{\pi }\,\cdot \,\Psi _{x}(\pi )$
\end{center}

\strut

and

\strut

\begin{center}
$E\left( (x+y)^{2k}\right) =\underset{\theta \in NC^{(even)}(2k)}{\sum }%
\,2^{\left| \theta \right| }\cdot \underset{\pi \in NC^{(even)}(2k),\,\pi
\leq \theta }{\sum }\mu _{\pi }^{\theta }\cdot \Psi _{x}(\pi ),$
\end{center}

\strut

for all $k\in \Bbb{N},$ where $\Psi _{x}=\Phi _{x}^{0}\circ \Lambda $ such
that

\strut

\begin{center}
$\Psi _{x}(\pi )=\widehat{E}(\pi )\left( x\otimes ...\otimes x\right) ,$ \
for all \ $\pi \in NC^{(even)}(2k).$
\end{center}
\end{remark}

\strut

\begin{remark}
The above theorems says that the $B$-valued moments and cmulants of $x+y$
are determined by certain recurrence relations, introduced in Section 4.2,
which are depending on partitions and images of numbering map of partitions.
\end{remark}

\strut

\strut

\subsection{The Trivial R-transform of $x+y$ and The Trivial Moment Series
of $x+y$}

\strut

\strut

By the previous section, we can conclude that

\strut

\begin{corollary}
Let $x=a+b+a^{-1}+b^{-1}$ and $y=c+d+c^{-1}+d^{-1}$ be given as before. Then

\strut

(1) $R_{x+y}^{t}(z)=\sum_{k=1}^{\infty }\left( 2\cdot \underset{\pi \in
NC^{(even)}(2k)}{\sum }\mu _{\pi }\,\cdot \,\Psi _{x}(\pi )\right) z^{2k}$

\strut

(2) $M_{x+y}^{t}(z)=\sum_{k=1}^{\infty }\left( \underset{\theta \in
NC^{(even)}(2k)}{\sum }\,2^{\left| \theta \right| }\cdot \underset{\pi \in
NC^{(even)}(2k),\,\pi \leq \theta }{\sum }\mu _{\pi }^{\theta }\cdot \Psi
_{x}(\pi )\right) z^{2k},$

\strut

where, for the fixed $\pi \in NC_{l_{1},...,l_{p}}(2k)\subset
NC^{(even)}(2k),$

\strut

\begin{center}
$\Psi _{x}(\pi )=\Phi _{x}\circ \Lambda (\pi ).$
\end{center}

$\square $
\end{corollary}

\strut

\strut

\strut

\section{Scalar-Valued Moments of $x+y$ in $\left(
L(F_{2})*_{L(F_{1})}L(F_{2}),\right) $}

\strut

\strut

\strut

\subsection{Scalar-Valued Moments of $x+y$}

\strut

\strut

\strut

In this section, finally, we will compute the scalar-valued moment,

\strut

\begin{center}
$\varphi \left( (x+y)^{n}\right) ,$ \ \ for \ $n\in \Bbb{N}.$
\end{center}

\strut

In Section 2.6, we showed that $x+y$ is even random variable in $%
(A_{1}*_{B}A_{2},\varphi ).$ More generally, if $a\in (A,E)$ is a $B$-even
and if $(A,E)$ and $(A,\varphi )$ are compatible, then $a$ is
(scalar-valued) even in $(A,\varphi ),$ where $B$ is a unital algebra and $A$
is an algebra over $B$ (See [15]). So, it suffices to consider

\strut

\begin{center}
$\varphi \left( (x+y)^{2k}\right) ,$
\end{center}

\strut

for all $k\in \Bbb{N},$ futhermore, by the compatibility, we have that

\strut

\begin{center}
$\varphi \left( (x+y)^{2k}\right) =\varphi \left( E((x+y)^{2k})\right) ,$
for all $k\in \Bbb{N}.$
\end{center}

\strut

In the previous section, we showed that

\strut \strut

\begin{center}
$R_{x+y}^{t}(z)=2R_{x}^{t}(z)=\sum_{k=1}^{\infty }\left( 2\cdot \underset{%
\pi \in NC^{(even)}(2k)}{\sum }\mu _{\pi }\,\cdot \,\Psi _{x}(\pi )\right)
z^{2k}$
\end{center}

\strut

and hence

\strut

\begin{center}
$M_{x+y}^{t}(z)=\sum_{k=1}^{\infty }\left( \underset{\theta \in
NC^{(even)}(2k)}{\sum }\,2^{\left| \theta \right| }\cdot \underset{\pi \in
NC^{(even)}(2k),\,\pi \leq \theta }{\sum }\mu (\pi ,\theta )\cdot \Psi
_{x}(\pi )\right) z^{2k}.$
\end{center}

\strut

By the compatibility of $\left( A_{1}*_{B}A_{2},E\right) $ and $\left(
A_{1}*_{B}A_{2},\varphi \right) ,$ we have that

\strut

\begin{center}
$\varphi \left( (x+y)^{n}\right) =\varphi \left( E((x+y)^{n})\right) ,$ for
all $n\in \Bbb{N}.$
\end{center}

\strut

By the $B$-evenness of $x$ and $y,$ since $x$ and $y$ are $B$-even, $x+y$ is
also $B$-even. Hence $x+y$ is (scalar-valued) even, too (See Section 2.6).
Thus we have that

\strut

\begin{center}
$\varphi \left( (x+y)^{n}\right) =0,$ whenever $n\in 2\Bbb{N}-1.$
\end{center}

\strut

So, we need to observe $\varphi \left( (x+y)^{2k}\right) ,$ for all $k\in 
\Bbb{N}.$ Fix $k\in \Bbb{N}.$ Then

\strut

$\ \ \varphi \left( (x+y)^{2k}\right) =\varphi \left( E\left(
(x+y)^{2k}\right) \right) $

\strut

$\ \ \ \ \ \ \ \ \ \ \ \ \ \ \ \ \ \ \ \ \ \ \ =\varphi \left( \underset{%
\theta \in NC^{(even)}(2k)}{\sum }\,2^{\left| \theta \right| }\cdot 
\underset{\pi \in NC^{(even)}(2k),\,\pi \leq \theta }{\sum }\mu _{\pi
}^{\theta }\cdot \Psi _{x}(\pi )\right) $

\strut

(5.1)

\strut

$\ \ \ \ \ \ \ \ \ \ \ \ \ \ \ \ \ \ \ \ \ \ \ =\underset{\theta \in
NC^{(even)}(2k)}{\sum }\,2^{\left| \theta \right| }\cdot \underset{\pi \in
NC^{(even)}(2k),\,\pi \leq \theta }{\sum }\mu _{\pi }^{\theta }\cdot \varphi
\left( \Psi _{x}(\pi )\right) .$

\strut

\strut \strut

\begin{theorem}
Let $x=a+b+a^{-1}+b^{-1}$ and $y=c+d+c^{-1}+d^{-1}$ in $\left(
L(F_{2})*_{L(F_{1})}L(F_{2}),F\right) ,$ where $F_{2}=\,<a,b>\,=\,<c,d>$ \
and $F_{1}=\,<h=aba^{-1}b^{-1}=cdc^{-1}d^{-1}>.$ Then

\strut

\begin{center}
$\varphi \left( (x+y)^{n}\right) =0,$ whenever $n\in 2\Bbb{N}-1$
\end{center}

and

\strut

\begin{center}
$\varphi \left( (x+y)^{2k}\right) =\underset{\theta \in NC^{(even)}(2k)}{%
\sum }\,2^{\left| \theta \right| }\cdot \underset{\pi \in
NC^{(even)}(2k),\,\pi \leq \theta }{\sum }\mu _{\pi }^{\theta }\cdot \varphi
\left( \Psi _{x}(\pi )\right) .$
\end{center}

$\square $
\end{theorem}

\strut

\strut

\strut

\subsection{Examples}

\strut

\strut \strut

In this section, we will compute the scalar-valued moments, $\tau \left(
(x+y)^{3}\right) $ and $\tau \left( (x+y)^{4}\right) .$

\strut \strut

\textbf{1}. $\tau \left( (x+y)^{6}\right) $ ;

\strut

By the previous result, we have that

\strut

\begin{center}
$\varphi \left( (x+y)^{6}\right) =\underset{\theta \in NC^{(even)}(6)}{\sum }%
\,2^{\left| \theta \right| }\cdot \underset{\pi \in NC^{(even)}(6),\,\pi
\leq \theta }{\sum }\mu _{\pi }^{\theta }\cdot \varphi \left( \Psi _{x}(\pi
)\right) .$
\end{center}

\strut \strut

(Step 1) Observe $NC^{(even)}(6)$ ;

\strut

\begin{center}
$NC^{(even)}(6)=NC_{2,2,2}(6)\cup NC_{2,4}(6)\cup NC_{6}(6).$
\end{center}

\strut

We have that $NC_{6}(6)=\{1_{6}\}$ and

\strut

$\ \ \ \ NC_{2,4}(6)=\{\{(1,2),(3,4,5,6)\},\,\,\{(1,4,5,6),(2,3)\}$

$\ \ \ \ \ \ \ \ \ \ \ \ \ \ \ \ \ \ \ \ \ \ \ \ \ \
\{(1,2,5,6),(3,4)\},\,\,\,\{(1,2,3,6),(4,5)\}$

$\ \ \ \ \ \ \ \ \ \ \ \ \ \ \ \ \ \ \ \ \ \ \ \ \ \
\{(1,2,3,4),(5,6)\},\,\,\,\{(1,6),(2,3,4,5)\}\}.$

\strut

The entries of the above set, $NC_{2,4}(6),$ is gotten from the $\frac{%
360^{o}\,\cdot \,k}{6}$-anticlockwise-rotations of the circular expression
of the first entry $\{(1,2),(3,4,5,6)\},$ where $k=0,1,...,5.$ We will call
this fixed entry, $\{(1,2),(3,4,5,6)\},$ a candidate. (Remark that, for
instance, $\{(1,4,5,6),(2,3)\},$ can be the entry of $NC_{2,4}(6),$ etc.)
This candidate can be determined by the relation, $1_{2}+1_{4}.$ We have that

\strut

$\tau \left( (x+y)^{6}\right) =\underset{\theta \in NC^{(even)}(6)}{\sum }%
\,2^{\left| \theta \right| }\cdot \underset{\pi \in NC^{(even)}(6),\,\pi
\leq \theta }{\sum }\mu _{\pi }^{\theta }\cdot \tau \left( \Psi _{x}(\pi
)\right) $

\strut

$\ \ =\underset{\theta \in NC_{2,2,2}(6)}{\sum }\,2^{\left| \theta \right|
}\cdot \underset{\pi \in NC^{(even)}(6),\,\pi \leq \theta }{\sum }\mu _{\pi
}^{\theta }\cdot \tau \left( \Psi _{x}(\pi )\right) $

$\ \ \ \ \ \ \ \ \ \ \ \ \ \ \ \ \ \ \ +\underset{\theta \in NC_{2,4}(6)}{%
\sum }\,2^{\left| \theta \right| }\cdot \underset{\pi \in
NC^{(even)}(6),\,\pi \leq \theta }{\sum }\mu _{\pi }^{\theta }\cdot \tau
\left( \Psi _{x}(\pi )\right) $

$\ \ \ \ \ \ \ \ \ \ \ \ \ \ \ \ \ \ \ +\underset{\theta \in NC_{6}(6)}{\sum 
}\,2^{\left| \theta \right| }\cdot \underset{\pi \in NC^{(even)}(6),\,\pi
\leq \theta }{\sum }\mu _{\pi }^{\theta }\cdot \tau \left( \Psi _{x}(\pi
)\right) $

\strut

$\ \ =\underset{\theta \in NC_{2,2,2}(6)}{\sum }\,2^{\left| \theta \right|
}\cdot \mu _{\theta }^{\theta }\cdot \tau \left( \Psi _{x}(\pi )\right) $

$\ \ \ \ \ \ \ \ \ \ \ \ \ \ \ \ \ \ \ +\underset{\theta \in NC_{2,4}(6)}{%
\sum }\,2^{\left| \theta \right| }\cdot \underset{\pi \in
NC^{(even)}(6),\,\pi \leq \theta }{\sum }\mu _{\pi }^{\theta }\cdot \tau
\left( \Psi _{x}(\pi )\right) $

$\ \ \ \ \ \ \ \ \ \ \ \ \ \ \ \ \ \ \ +\underset{\theta \in NC^{(even)}(6)}{%
\sum }\,2^{\left| 1_{6}\right| }\cdot \mu _{\theta }\cdot \tau \left( \Psi
_{x}(\pi )\right) ,$

\strut

where $\mu _{\theta }^{\theta }=\mu (\theta ,\theta )=1$ and $\mu _{\theta
}=\mu (\theta ,1_{6}).$ Notice that the last line of the above formular is
nothing but the $6$-th cumulant of $x+y,$ because that $2^{\left|
1_{6}\right| }=2$ and $x$ and $y$ are $B$-valued identically distributed.

\strut

(Step 2) Observe $\underset{\theta \in NC_{2,2,2}(6)}{\sum }2^{\left| \theta
\right| }\cdot \,\mu _{\theta }^{\theta }\cdot \tau \left( \Psi _{x}(\pi
)\right) $ ;

\strut

Fix $\theta \in NC_{2,2,2}(6).$ If $\pi \in NC^{(even)}(6)$ satisfies $\pi
\leq \theta ,$ then $\pi =\theta ,$ by the ordering on $NC(6).$ Therefore,

\strut

$\underset{\theta \in NC_{2,2,2}(6)}{\sum }2^{\left| \theta \right| }\cdot
\,\mu _{\theta }^{\theta }\cdot \tau \left( \Psi _{x}(\pi )\right) $

$\ \ \ \ \ \ \ \ \ \ \ \ \ \ \ \ \ =\underset{\theta \in NC_{2,2,2}(6)}{\sum 
}2^{3}\cdot \tau \left( \Psi _{x}(\theta )\right) =\underset{\theta \in
NC_{2,2,2}(6)}{\sum }8\cdot \tau \left( E(x^{2})^{3}\right) ,$

\strut

since $\left| \theta \right| =3,$ for all $\theta \in NC_{2,2,2}(6).$ So,

\begin{center}
$\underset{\theta \in NC_{2,2,2}(6)}{\sum }8\cdot \tau \left( \Psi
_{x}(\theta )\right) =\left| NC_{2,2,2}(6)\right| \cdot \left( 8\cdot \left(
p_{0}^{2}\right) ^{3}\right) .$
\end{center}

\strut \strut

Therefore, since $\left| NC_{2,2,2}(6)\right| =c_{2}^{3}=5,$

\strut

$\ \ \underset{\theta \in NC_{2,2,2}(6)}{\sum }2^{\left| \theta \right|
}\cdot \,\underset{\pi \in NC^{(even)}(6),\,\pi \leq \theta }{\sum }\mu
_{\pi }^{\theta }\cdot \tau \left( \Psi _{x}(\pi )\right) $

\begin{center}
$=5\cdot 8\cdot \left( p_{0}^{2}\right) ^{3}=2560.$
\end{center}

\strut

(Step 3) Observe $\underset{\theta \in NC_{2,4}(6)}{\sum }2^{\left| \theta
\right| }\cdot \,\underset{\pi \in NC^{(even)}(6),\,\pi \leq \theta }{\sum }%
\mu _{\pi }^{\theta }\cdot \tau \left( \Psi _{x}(\pi )\right) $ ;

\strut

We know all entries, $\theta ,$ of $NC_{2,4}(6),$ from the (Step 1). It is
easy to check that if $\theta \in NC_{2,4}(6)$ and if $\pi \in
NC^{(even)}(6) $ satisfies $\pi \leq \theta ,$ then $\pi =\theta $ or $\pi
\in NC_{2,2,2}(6).$ Moreover, each partition $\theta \in NC_{2,4}(6)$
contains exactly two partitions $\pi $ in $NC_{2,2,2}(6)$ such that $\pi
<\theta .$ For example, if $\theta =\{(1,2),(3,4,5,6)\},$ then we have that $%
\pi _{1}=\{(1,2),(3,4),(5,6)\}$ and $\pi _{2}=\{(1,2),(3,6),(4,5)\}$ in $%
NC_{2,2,2}(6).$

\strut

Similar to (Step 2), we have that

\strut \strut

$\underset{\theta \in NC_{2,4}(6)}{\sum }2^{\left| \theta \right| }\cdot \,%
\underset{\pi \in NC_{2,2,2}(6),\,\pi \leq \theta }{\sum }\mu _{\pi
}^{\theta }\cdot \tau \left( \Psi _{x}(\pi )\right) $

$\ =\left( \left| NC_{2,4}(6)\right| \right) \cdot 2^{2}\cdot \left( \tau
\left( \Psi _{x}(\theta )\right) -\left( \tau \left( \Psi _{x}(\pi
_{1})\right) +\tau \left( \Psi _{x}(\pi _{2})\right) \right) \right) $

$\strut $

$\ =\left( \left| NC_{2,4}(6)\right| \right) \cdot 4\cdot \left( \tau \left(
\Psi _{x}(\theta )\right) -2\left( \tau \left( \Psi _{x}(\pi _{1})\right)
\right) \right) ,$

\strut

where $\pi _{1},\pi _{2}\in NC_{2,2,2}(6)$ such that $\pi _{i}\leq \theta ,$ 
$\theta \in NC_{2,4}(6)$ is arbitrarily fixed. We can get the last line of
the above formular, since all $\Psi _{x}(\pi )$'s are same, for all $\pi \in
NC_{2,2,2}(6).$ Also, we have that

$\strut $

\begin{center}
$
\begin{array}{ll}
\mu _{\pi }^{\theta } & =\mu \left( [0_{1},1_{1}]^{2}\times
[0_{2},1_{2}]\right) \\ 
& =\mu (0_{1},1_{1})^{2}\cdot \mu (0_{2},1_{2})=(-1)^{2-1}c_{2}^{2-1} \\ 
& =-1,
\end{array}
$
\end{center}

\strut

for all pair $\left( \pi _{i},\theta \right) \in NC_{2,2,2}(6)\times
NC_{2,4}(6),$ where $\pi _{i}\leq \theta .$ Futhermore, by (Step 1), we know
all entries of $NC_{2,4}(6).$ Hence we can compute each $\Psi _{x}(\theta ).$

\strut

$\Psi _{x}\left( \{(1,2),(3,4,5,6)\}\right) =\Phi \left( (2)\right) \Phi
\left( (4)\right) =p_{0}^{2}\left( (h+h^{-1})+p_{0}^{4}\right) .$

\strut

$\Psi _{x}\left( \{(1,4,5,6),(2,3)\}\right) =\Phi \left( 1,[1],3\right)
=\left( (h+h^{-1})+p_{0}^{4}\right) p_{0}^{2}.$

\strut

$\Psi _{x}\left( \{(1,2,5,6),(3,4)\}\right) =\Phi \left( (2,[1],2)\right)
=\left( (h+h^{-1})+p_{0}^{4}\right) p_{0}^{2}.$

\strut

$\Psi _{x}\left( \{(1,2,3,6),(4,5)\}\right) =\Phi \left( (3,[1],1)\right)
=\left( (h+h^{-1})+p_{0}^{4}\right) .$

\strut

$\Psi _{x}\left( \{(1,2,3,4),(5,6)\}\right) =\Phi \left( (4)\right) \Phi
\left( (2)\right) =\left( (h+h^{-1})+p_{0}^{4}\right) p_{0}^{2}.$

\strut

$\Psi _{x}\left( \{(1,6),(2,3,4,5)\}\right) =\Phi \left( (1,[4],1)\right) $

\ \ \ \ \ \ \ \ \ \ \ \ \ \ \ \ \ \ \ \ \ \ \ \ \ \ \ \ \ \ \ \ \ \ \ \ \ \ $%
=E\left( x(h+h^{-1})x\right) +p_{0}^{4}p_{0}^{2}=0_{B}+p_{0}^{4}p_{0}^{2}$

\ \ \ \ \ \ \ \ \ \ \ \ \ \ \ \ \ \ \ \ \ \ \ \ \ \ \ \ \ \ \ \ \ \ \ \ \ \ $%
=p_{0}^{2}p_{0}^{4}.$

\strut \strut

And, for any $\pi \in NC^{(even)}(6)$ such that $\pi <\theta ,$ (i.e, $\pi
\in NC_{2,2,2}(6)$ !)

\strut

\begin{center}
$\mu _{\pi }=\mu (\pi ,\theta )=\mu \left( [0_{1},1_{1}]^{2}\times
[0_{2},1_{2}]\right) =-1.$
\end{center}

Therefore,

\strut

$\left( \left| NC_{2,4}(6)\right| \right) \cdot 4\cdot \left( \tau \left(
\Psi _{x}(\theta )\right) -2\left( \tau \left( \Psi _{x}(\pi _{1})\right)
\right) \right) $

\strut

\begin{center}
$=6\cdot 4\cdot \left( 112-128\right) =24\cdot (-16)=-384.$
\end{center}

\strut

i.e,

\strut

\begin{center}
$\underset{\theta \in NC_{2,4}(6)}{\sum }2^{\left| \theta \right| }\cdot \,%
\underset{\pi \in NC_{2,2,2}(6),\,\pi \leq \theta }{\sum }\mu _{\pi
}^{\theta }\cdot \tau \left( \Psi _{x}(\pi )\right) =-384.$
\end{center}

\strut

(Step 4) Observe $\,\underset{\pi \in NC^{(even)}(6)}{\sum }\mu _{\pi }\cdot
\tau \left( \Psi _{x}(\pi )\right) $ ;

\strut

It is easy to check that

\strut

$\underset{\pi \in NC^{(even)}(6)}{\sum }2\cdot \mu _{\pi }\cdot \tau \left(
\Psi _{x}(\pi )\right) =\underset{\pi \in NC_{2,2,2}(6)}{\sum }2\cdot \mu
_{\pi }\cdot \tau \left( \Psi _{x}(\pi )\right) $

\begin{center}
$+\underset{\pi \in NC_{2,4}(6)}{\sum }2\cdot \mu _{\pi }\cdot \tau \left(
\Psi _{x}(\pi )\right) +2\cdot \tau \left( \Psi _{x}(1_{6})\right) .$
\end{center}

\strut

Also, it is easy to see that

\strut

$\underset{\pi \in NC_{2,2,2}(6)}{\sum }2\cdot \mu _{\pi }\cdot \varphi
\left( \Psi _{x}(\pi )\right) =2\cdot 2\cdot \mu _{\pi _{1}}\tau \left( \Psi
_{x}(\pi _{1})\right) +3\cdot 2\cdot \mu _{\pi _{2}}\tau \left( \Psi
_{x}(\pi _{2})\right) $

$\strut $

since there are two kinds of block structures in $NC_{2,2,2}(6)$ ; one kind
is

\strut

\begin{center}
$\{(1,2),(3,4),(5,6)\}$ and its rotations
\end{center}

\strut

(there are two such partitions) and another kind is

\strut

\begin{center}
$\{(1,2),(3,6),(4,5)\}$ and its rotations
\end{center}

\strut

(there are three such partitions)$\strut $, so, we have that,

$\strut $

$\ \ \ \ \ \ \ \ \ \ \ \ \ \ \ =4\cdot 2\cdot (p_{0}^{2})^{3}+6\cdot 1\cdot
(p_{0}^{2})^{3}=8\cdot 64+6\cdot 64=896.$

\strut

Also,

\strut

\begin{center}
$\underset{\pi \in NC_{2,4}(6)}{\sum }2\cdot \mu _{\pi }\cdot \tau \left(
\Psi _{x}(\pi )\right) =6\cdot 2\cdot (-1)\cdot 112=-1344$
\end{center}

\strut

and

\strut \strut

\begin{center}
\strut $2\cdot \tau \left( \Psi _{x}(1_{6})\right) =2\cdot p_{0}^{6}=2\cdot
232=464.$
\end{center}

\strut

Therefore,

\strut

\begin{center}
$\underset{\pi \in NC^{(even)}(6)}{\sum }2\cdot \mu _{\pi }\cdot \tau \left(
\Psi _{x}(\pi )\right) =16.$
\end{center}

\strut

(Step 5) Add all information ;

\strut

\begin{center}
$\tau \left( (x+y)^{6}\right) =2560-384+16=2192.$
\end{center}

\strut

By (Step 1) $\sim $\ (Step 5), we can get that

\strut

\begin{center}
$\tau \left( (x+y)^{6}\right) =2192.$
\end{center}

\strut \strut

\begin{example}
Let $x,y\in (L(F_{2})*_{L(F_{1})}L(F_{2}),F)$ be $L(F_{1})$-valued random
variables such that $x=a+b+a^{-1}+b^{-1}$ and $y=c+d+c^{-1}+d^{-1}$ and let $%
E:L(F_{2})*_{L(F_{1})}L(F_{2})\rightarrow L(F_{1})$ and $\tau
:L(F_{2})*_{L(F_{1})}L(F_{2})\rightarrow \Bbb{C}$ be the conditional
expectation (finding $h$-terms) and the canonical trace, respectively. Then

\strut

\begin{center}
$\tau \left( (x+y)^{6}\right) =\tau \left( E\left( (x+y)^{6}\right) \right)
=2192.$
\end{center}

$\square $
\end{example}

\strut

\begin{remark}
The above result also gotten from the following way ; First, recall that $%
L(F_{2})*_{L(F_{1})}L(F_{2})\simeq L(F_{2}*_{F_{1}}F_{2}).$ Also, we can
regard the group $F_{2}*_{F_{1}}F_{2}$ as a (topological) fundamental group
of torus with genus 2,

\strut

\begin{center}
$G=\,<a,b,c,d:aba^{-1}b^{-1}d^{-1}c^{-1}dc=e>.$
\end{center}

\strut \strut

(Remember that $aba^{-1}b^{-1}=cdc^{-1}d^{-1}$ is our $h$ !) We need to
recognize that $aba^{-1}b^{-1}d^{-1}c^{-1}dc$ is a word with length 8,
without considering the relation in the group $G.$ Again, denote $%
a+b+a^{-1}+b^{-1}$ and $c+d+c^{-1}+d^{-1}$ by $x$ and $y,$ respectively.
Now, define the following trace

\strut

\begin{center}
$\tau _{4}$ $:L(F_{4})\rightarrow \Bbb{C}$
\end{center}

by

\begin{center}
$\tau _{4}\left( \underset{g\in F_{4}}{\sum }\alpha _{g}g\right) =\alpha
_{e_{F_{4}}},$
\end{center}

\strut

where $F_{4}=\,<a,b,c,d>.$ Notice that

\strut

\begin{center}
$\tau _{4}\left( (x+y)^{6}\right) =\tau \left( (x+y)^{6}\right) ,$
\end{center}

\strut

because, in both cases, we cannot make the words with length 8 in $\left(
x+y\right) ^{6}.$ (Of course, in our case, the word $%
aba^{-1}b^{-1}d^{-1}c^{-1}dc$ is $e,$ but this can be come from making the
words with length 8 !) Now, let's compute $\tau _{4}\left( (x+y)^{6}\right)
. $ This can be computed by using the method introduced in [35], as follows
; this method is also used in Section 3.3.

\strut \strut

$\ \tau _{4}\left( (x+y)^{6}\right) =\tau \left(
(a+b+a^{-1}+b^{-1}+c+d+c^{-1}+d^{-1})^{6}\right) $

\strut

$\ \ \ \ \ \ \ \ \ \ \ \ \ \ \ \ \ \ \ \ =\tau _{4}\left( (X_{1})^{6}\right)
,$

\strut

where $X_{1}=$ the sum of length 1 words in $L(F_{4}).$ We have the
following recurrence relations, by [35] ;

\strut

\begin{center}
$X_{1}X_{1}=X_{2}+8e$
\end{center}

and

\begin{center}
$X_{1}X_{N}=X_{N+1}+7X_{N-1},$ for all $N\geq 2.$
\end{center}

\strut

So, to get $\tau _{4}\left( (x+y)^{6}\right) =\tau _{4}\left(
(X_{1})^{6}\right) ,$ we need to compute that

\strut

$\ (x+y)^{2}=(X_{1})^{2}=X_{1}X_{1}=X_{2}+8e,$

\strut

$\ 
\begin{array}{ll}
(x+y)^{3} & =(X_{1})^{3}=X_{1}\left( X_{2}+8e\right) =X_{1}X_{2}+8X_{1} \\ 
& =(X_{3}+7X_{1})+8X_{1}=X_{3}+15X_{1},
\end{array}
\ \ \ \ \ \ $

\strut \strut

$\ 
\begin{array}{ll}
(x+y)^{4} & =(X_{1})^{4}=X_{1}\left( X_{3}+15X_{1}\right)
=X_{1}X_{3}+15X_{1}X_{1} \\ 
& =X_{4}+7X_{2}+15(X_{2}+8e) \\ 
& =X_{4}+22X_{2}+120e,
\end{array}
$

\strut

$\ 
\begin{array}{ll}
(x+y)^{5} & =(X_{1})^{5}=X_{1}\left( X_{4}+22X_{2}+120e\right) \\ 
& =X_{5}+29X_{3}+274X_{1},
\end{array}
$

\strut

and

\strut

$\ 
\begin{array}{ll}
(x+y)^{6} & =(X_{1})^{6}=X_{1}\left( X_{5}+29X_{3}+274X_{1}\right) \\ 
& =X_{6}+36X_{4}+203X_{2}+274X_{2}+2192e.
\end{array}
$

\strut

Thus, we have that

\strut

\begin{center}
$
\begin{array}{ll}
\tau _{4}\left( (x+y)^{6}\right) & =\tau _{4}\left(
X_{6}+36X_{4}+203X_{2}+274X_{2}+2192e\right) \\ 
& =2192.
\end{array}
$
\end{center}

\strut

Therefore, we can conclude that

\strut

\begin{center}
$\tau \left( (x+y)^{6}\right) =2192=\tau _{4}\left( (x+y)^{6}\right) .$
\end{center}
\end{remark}

\strut

\strut

\textbf{2}. $\tau \left( (x+y)^{8}\right) $ ; \ by Section 4.4, we have that

\strut

\begin{center}
$\tau \left( (x+y)^{8}\right) =\underset{\theta \in NC^{(even)}(8)}{\sum }%
\,2^{\left| \theta \right| }\cdot \underset{\pi \in NC^{(even)}(8),\,\pi
\leq \theta }{\sum }\mu _{\pi }^{\theta }\cdot \tau \left( \Psi _{x}(\pi
)\right) .$
\end{center}

\strut

By the separation of $NC^{(even)}(8),$ we have that

\strut

\begin{center}
$NC^{(even)}(8)=NC_{2,2,2,2}(8)\cup NC_{2,2,4}(8)\cup NC_{2,6}(8)\cup
NC_{4,4}(8)\cup \{1_{8}\}.$
\end{center}

\strut

Therefore,

\strut

$\tau \left( (x+y)^{8}\right) $

\strut

\ $\ =\underset{\theta \in NC_{2,2,2,2}(8)}{\sum }\,2^{4}\underset{\pi \in
NC^{(even)}(8),\,\pi \leq \theta }{\sum }\mu _{\pi }^{\theta }\cdot \tau
\left( \Psi _{x}(\pi )\right) $

\strut

$\ \ \ \ \ \ \ \ \ \ \ \ \ \ \ +\underset{\theta \in NC_{2,2,4}(8)}{\sum }%
\,2^{3}\cdot \underset{\pi \in NC^{(even)}(8),\,\pi \leq \theta }{\sum }\mu
_{\pi }^{\theta }\cdot \tau \left( \Psi _{x}(\pi )\right) $

\strut

$\ \ \ \ \ \ \ \ \ \ \ \ \ \ \ +\underset{\theta \in NC_{2,6}(8)}{\sum }%
\,2^{2}\cdot \underset{\pi \in NC^{(even)}(8),\,\pi \leq \theta }{\sum }\mu
_{\pi }^{\theta }\cdot \tau \left( \Psi _{x}(\pi )\right) $

\strut

$\ \ \ \ \ \ \ \ \ \ \ \ \ \ \ +\underset{\theta \in NC_{4,4}(8)}{\sum }%
\,2^{2}\cdot \underset{\pi \in NC^{(even)}(8),\,\pi \leq \theta }{\sum }\mu
_{\pi }^{\theta }\cdot \tau \left( \Psi _{x}(\pi )\right) $

\strut

$\ \ \ \ \ \ \ \ \ \ \ \ \ \ \ +\,2\cdot K_{8}^{t}\left( x,...,x\right) $

\strut

$\ \ \ =\underset{\theta \in NC_{2,2,2,2}(8)}{\sum }\,(16)\cdot \mu _{\theta
}^{\theta }\cdot \tau \left( \Psi _{x}(\theta )\right) $

\strut

(since there is no $\pi \in NC^{(even)}(8)$ such that $\pi \lvertneqq \theta
,$ for $\theta \in NC_{2,2,2,2}(8)$)

\strut

$\ \ \ \ \ \ \ \ \ \ \ +\underset{\theta \in NC_{2,2,4}(8)}{\sum }\,(8)\cdot
\cdot \left( \mu _{\theta }^{\theta }\tau \left( \Psi _{x}(\theta )\right)
+\mu _{\pi _{1}}^{\theta }\tau \left( \Psi _{x}(\pi _{1})\right) +\mu _{\pi
_{2}}^{\theta }\tau \left( \Psi _{x}(\pi _{2})\right) \right) $

\strut

(for each given $\theta \in NC_{2,2,4}(8),$ we have only two proper
partitions $\pi _{1},\pi _{2}$ such that $\pi _{i}\lvertneqq \theta ,$ $%
i=1,2 $)

\strut

$\ \ \ \ \ \ \ \ \ \ \ +\underset{\theta \in NC_{2,6}(8)}{\sum }\,(4)\cdot
(\mu _{\theta }^{\theta }\tau \left( \Psi _{x}(\theta )\right) +6\cdot (\mu
_{\pi }^{\theta }\tau \left( \Psi _{x}(\pi )\right) )$

$\ \ \ \ \ \ \ \ \ \ \ \ \ \ \ \ \ \ \ \ \ \ \ \ \ \ \ \ \ \ \ \ \ +\left(
2\cdot (\mu _{\pi _{1}}^{\theta }\tau \left( \Psi _{x}(\pi _{1})\right)
)+3\cdot (\mu _{\pi _{2}}^{\theta }\tau \left( \Psi _{x}(\pi _{2})\right)
)\right) )$

\strut

(for each given $\theta \in NC_{2,6}(8),$ like in the \textbf{step 5} in 
\textbf{1}., we have the proper partitions $\pi \in NC_{2,2,4}(8)$ and $\pi
_{1},\,\pi _{2}\in NC_{2,2,2,2}(8).$ Notice that $\pi _{1}$ and $\pi _{2}$
have different type of pairings in $NC_{2,2,2}(6)\hookrightarrow
NC_{2,2,2,2}(8)$)

\strut

$\ \ \ \ \ \ \ \ \ \ \ +\underset{\theta \in NC_{4,4}(8)}{\sum }\,(4)\cdot 
\underset{\pi \in NC^{(even)}(8),\,\pi \leq \theta }{\sum }\mu _{\pi
}^{\theta }\cdot \tau \left( \Psi _{x}(\pi )\right) $

\strut

$\ \ \ \ \ \ \ \ \ \ \ +2\cdot \underset{\theta \in NC^{(even)}(8)}{\sum }%
\mu _{\theta }^{1_{8}}\cdot \tau \left( \Psi _{x}(\theta )\right) $

\strut

(\textbf{Step 1}) Compute

\strut

$\ \ \ \ \ \underset{\theta \in NC_{2,2,2,2}(8)}{\sum }\,(16)\cdot \mu
_{\theta }^{\theta }\cdot \tau \left( \Psi _{x}(\theta )\right) =\underset{%
\theta \in NC_{2,2,2,2}(8)}{\sum }\,(16)\cdot \tau \left( \Psi _{x}(\theta
)\right) $

\strut

since $\mu _{\theta }^{\theta }=1\in \Bbb{C}$

\strut

\ \ \ \ \ \ \ \ \ \ \ \ \ \ \ \ \ \ \ \ \ $\ =(14)(16)\varphi \left(
E(x^{2})^{4}\right) =(14)(16)(256)=57344.$

\strut

(\textbf{Step 2}) Compute

\strut

$\ \underset{\theta \in NC_{2,2,4}(8)}{\sum }\,(8)\cdot \left( \mu _{\theta
}^{\theta }\tau \left( \Psi _{x}(\theta )\right) +\mu _{\pi _{1}}^{\theta
}\tau \left( \Psi _{x}(\pi _{1})\right) +\mu _{\pi _{2}}^{\theta }\tau
\left( \Psi _{x}(\pi _{2})\right) \right) $

\strut

\ $\ \ \ \ \ \ \ \ \ \ \ =\underset{\theta \in NC_{2,2,4}(8)}{\sum }%
\,(8)\cdot \left( \tau \left( \Psi _{x}(\theta )\right) -\tau \left( \Psi
_{x}(\pi _{1})\right) -\tau \left( \Psi _{x}(\pi _{2})\right) \right) $

\strut

\ $\ \ \ \ \ \ \ \ \ \ \ =\underset{\theta \in NC_{2,2,4}(8)}{\sum }%
\,(8)\cdot \left(
(p_{0}^{2}p_{0}^{2}p_{0}^{4})-(p_{0}^{2})^{4}-(p_{0}^{2})^{4}\right) $

\strut

\ $\ \ \ \ \ \ \ \ \ \ \ =\underset{\theta \in NC_{2,2,4}(8)}{\sum }%
\,(8)\cdot \left( 448-256-256\right) =\underset{\theta \in NC_{2,2,4}(8)}{%
\sum }\,(8)\cdot (-64)$

\strut

\ $\ \ \ \ \ \ \ \ \ \ \ =\underset{\theta \in NC_{2,2,4}(8)}{\sum }\left(
-512\right) =(28)\cdot (-512)$

\strut

\ $\ \ \ \ \ \ \ \ \ \ \ =-14336.$

\strut

We can get that $\left| NC_{2,2,4}(8)\right| =28.$

\strut

(\textbf{Step 3}) Compute

\strut

$\ \underset{\theta \in NC_{2,6}(8)}{\sum }\,(4)\cdot (\mu _{\theta
}^{\theta }\tau \left( \Psi _{x}(\theta )\right) +6\cdot (\mu _{\pi
}^{\theta }\tau \left( \Psi _{x}(\pi )\right) )$

$\ \ \ \ \ \ \ \ \ \ \ \ \ \ \ \ \ \ \ \ \ \ \ \ \ \ \ \ \ \ \ \ \ +\left(
2\cdot (\mu _{\pi _{1}}^{\theta }\tau \left( \Psi _{x}(\pi _{1})\right)
)+3\cdot (\mu _{\pi _{2}}^{\theta }\tau \left( \Psi _{x}(\pi _{2})\right)
)\right) )$

\strut

$\ \ \ \ \ \ \ \ \ \ =\underset{\theta \in NC_{2,6}(8)}{\sum }\,(4)\cdot
(\tau \left( \Psi _{x}(\theta )\right) -6\cdot \tau \left( \Psi _{x}(\pi
)\right) $

$\ \ \ \ \ \ \ \ \ \ \ \ \ \ \ \ \ \ \ \ \ \ \ \ \ \ \ \ \ \ \ \ +2\cdot
2\cdot \tau \left( \Psi _{x}(\pi _{1})\right) +3\cdot 1\cdot \tau \left(
\Psi _{x}(\pi _{2})\right) )$

\strut

$\ \ \ \ \ \ \ \ \ \ =\underset{\theta \in NC_{2,6}(8)}{\sum }\,(4)\cdot
\left( p_{0}^{2}p_{0}^{6}-6\cdot p_{0}^{2}p_{0}^{2}p_{0}^{4}+4\cdot
(p_{0}^{2})^{4}+3\cdot (p_{0}^{2})^{4}\right) $

\strut

$\ \ \ \ \ \ \ \ \ \ =\underset{\theta \in NC_{2,6}(8)}{\sum }\,(4)\cdot
\left( 928-2688+1024+768\right) $

\strut

$\ \ \ \ \ \ \ \ \ \ =\underset{\theta \in NC_{2,6}(8)}{\sum }\,(4)\cdot
\left( 32\right) =\underset{\theta \in NC_{2,6}(8)}{\sum }\left( 128\right) $

\strut \strut

$\ \ \ \ \ \ \ \ \ \ =8\cdot (128)=1024.$

\strut

(\textbf{Step 4}) Compute

\strut

$\ \ \underset{\theta \in NC_{4,4}(8)}{\sum }\,(4)\cdot \underset{\pi \in
NC^{(even)}(8),\,\pi \leq \theta }{\sum }\mu _{\pi }^{\theta }\cdot \tau
\left( \Psi _{x}(\pi )\right) $

\strut

$\ \ \ \ \ \ \ \ \ \ \ \ =\underset{\theta \in NC_{4,4}(8)}{\sum }\,(4)\cdot
(\mu _{\theta }^{\theta }\cdot \varphi \left( \Psi _{x}(\theta )\right) +\mu
_{\pi _{1}}^{\theta }\cdot \tau \left( \Psi _{x}(\pi _{1})\right) $

\ \ \ \ \ $\ \ \ \ \ \ \ \ \ \ \ \ \ \ \ \ \ \ \ \ \ \ \ \ \ \ \ \ \ \ \ \ \
+\mu _{\pi _{2}}^{\theta }\cdot \tau \left( \Psi _{x}(\pi _{2})\right) +\mu
_{\pi _{3}}^{\theta }\cdot \tau \left( \Psi _{x}(\pi _{3})\right) $

$\ \ \ \ \ \ \ \ \ \ \ \ \ \ \ \ \ \ \ \ \ \ \ \ \ \ \ \ \ \ \ \ \ \ \ \ \ \
+\mu _{\pi _{4}}^{\theta }\cdot \tau \left( \Psi _{x}(\pi _{4})\right) +\mu
_{\pi _{5}}^{\theta }\cdot \tau \left( \Psi _{x}(\pi _{5})\right) $

$\ \ \ \ \ \ \ \ \ \ \ \ \ \ \ \ \ \ \ \ \ \ \ \ \ \ \ \ \ \ \ \ \ \ \ \ \ \
+\mu _{\pi _{6}}^{\theta }\cdot \tau \left( \Psi _{x}(\pi _{6})\right) +\mu
_{\pi _{7}}^{\theta }\cdot \tau \left( \Psi _{x}(\pi _{6})\right) $

$\ \ \ \ \ \ \ \ \ \ \ \ \ \ \ \ \ \ \ \ \ \ \ \ \ \ \ \ \ \ \ \ \ \ \ \ \ \
+\mu _{\pi _{8}}^{\theta }\cdot \tau \left( \Psi _{x}(\pi _{6})\right) ),$

\strut

since for each $\theta \in NC_{4,4}(8),$ there are proper partitions $\pi
_{1},\pi _{2},\pi _{3},\pi _{4}\in NC_{2,2,4}(8)$ and $\pi _{5},\pi _{6},\pi
_{7},\pi _{8}\in NC_{2,2,2,2}(8).$

\strut

$\ \ \ \ \ \ \ \ \ \ \ \ =\underset{\theta \in NC_{4,4}(8)}{\sum }\,(4)\cdot
(\left( (p_{0}^{4})^{2}+2\right) -p_{0}^{2}p_{0}^{2}p_{0}^{4}$

\ \ \ \ \ $\ \ \ \ \ \ \ \ \ \ \ \ \ \ \ \ \ \ \ \ \ \ \ \ \ \ \ \ \ \ \ \ \
-p_{0}^{2}p_{0}^{2}p_{0}^{4}-p_{0}^{2}p_{0}^{2}p_{0}^{4}$

$\ \ \ \ \ \ \ \ \ \ \ \ \ \ \ \ \ \ \ \ \ \ \ \ \ \ \ \ \ \ \ \ \ \ \ \ \ \
-p_{0}^{2}p_{0}^{2}p_{0}^{4}+(p_{0}^{2})^{4}$

$\ \ \ \ \ \ \ \ \ \ \ \ \ \ \ \ \ \ \ \ \ \ \ \ \ \ \ \ \ \ \ \ \ \ \ \ \ \
+(p_{0}^{2})^{4}+(p_{0}^{2})^{4}$

$\ \ \ \ \ \ \ \ \ \ \ \ \ \ \ \ \ \ \ \ \ \ \ \ \ \ \ \ \ \ \ \ \ \ \ \ \
+(p_{0}^{2})^{4}),$

\strut

$\ \ \ \ \ \ \ \ \ \ \ \ =\underset{\theta \in NC_{4,4}(8)}{\sum }\,(4)\cdot
\left( 786+(4)(-448)+4(256)\right) $

\strut

$\ \ \ \ \ \ \ \ \ \ \ =\underset{\theta \in NC_{4,4}(8)}{\sum }\,(4)\cdot
\left( 786-1792+1024\right) =\underset{\theta \in NC_{4,4}(8)}{\sum }(4)(18)$

$\ \ \ \ \ \ \ \ \ \ \ \ =\underset{\theta \in NC_{4,4}(8)}{\sum }%
(72)=(4)(72)$

$\strut $

$\ \ \ \ \ \ \ \ \ \ \ =288.$

\strut

(\textbf{Step 5}) Compute

\strut

$\ \ 2\cdot \underset{\theta \in NC^{(even)}(8)}{\sum }\mu _{\theta
}^{1_{8}}\cdot \tau \left( \Psi _{x}(\theta )\right) $

\strut

$\ \ \ \ \ \ \ =2\cdot (2\cdot (p_{0}^{2})^{4}\mu (0_{4},1_{4})+8\cdot
(p_{0}^{2})^{4}\mu (0_{2},1_{2})\mu (0_{3},1_{3})$

$\ \ \ \ \ \ \ \ \ \ \ \ \ \ \ \ \ \ \ \ \ \ \ \ \ \ \ \ \ \ \ \ \ \ \ \ \ \
\ \ \ \ \ \ \ \ +4\cdot (p_{0}^{2})^{4}\left( \mu (0_{2},1_{2})\right) ^{3})$

\strut

\ \ \ \ \ $\ \ \ \ \ \ \ +2\cdot (8\cdot (p_{0}^{2})^{2}(p_{0}^{4})\mu
(0_{3},1_{3})+4\cdot (p_{0}^{2})^{2}(p_{0}^{4})\left( \mu
(0_{2},1_{2})\right) ^{2}$

$\ \ \ \ \ \ \ \ \ \ \ \ \ \ \ \ \ \ \ \ \ \ \ \ \ \ \ \ \ \ \ \ \ \ \ \ \ \
\ \ \ \ \ \ \ \ \ \ \ \ \ \ \ \ \ +8\cdot (p_{0}^{2})^{2}(p_{0}^{4})\left(
\mu (0_{2},1_{2})\right) ^{2}$

$\ \ \ \ \ \ \ \ \ \ \ \ \ \ \ \ \ \ \ \ \ \ \ \ \ \ \ \ \ \ \ \ \ \ \ \ \ \
\ \ \ \ \ \ \ \ \ \ \ \ \ \ \ \ +8\cdot (p_{0}^{2})^{2}(p_{0}^{4})\left( \mu
(0_{2},1_{2})\right) ^{2})$

\strut

$\ \ \ \ \ \ \ \ \ \ \ \ +2\cdot \left( 8\cdot p_{0}^{2}p_{0}^{6}\,\mu
(0_{2},1_{2})\right) $

\strut

$\ \ \ \ \ \ \ \ \ \ \ \ +2\cdot \left( \underset{\theta \in NC_{4,4}(8)}{%
\sum }\mu _{\theta }^{1_{8}}\cdot \varphi \left( \Psi _{x}(\theta )\right)
\right) $

\strut

$\ \ \ \ \ \ \ \ \ \ \ \ +2\cdot p_{0}^{8}$

\strut

$\ \ \ \ \ \ \ =2\left( -2560-4096-1024\right) +2\left(
7168+1792+3584+3584\right) $

$\ \ \ \ \ \ \ \ \ \ \ \ \ \ \ \ \ \ \ \ +2\left( -7424\right) +2(\underset{%
\theta \in NC_{4,4}(8)}{\sum }\mu _{\theta }^{1_{8}}\cdot \varphi \left(
\Psi _{x}(\theta )\right) )+2\left( 2092\right) $

\strut

$\ \ \ \ \ \ =-15360+32256-14848$

$\ \ \ \ \ \ \ \ \ \ \ \ \ \ \ \ \ \ \ +2\left( \underset{\theta \in
NC_{4,4}(8)}{\sum }\mu _{\theta }^{1_{8}}\cdot \tau \left( \Psi _{x}(\theta
)\right) \right) +4184$

\strut

$\ \ \ \ \ \ =2\left( \underset{\theta \in NC_{4,4}(8)}{\sum }\mu _{\theta
}^{1_{8}}\cdot \tau \left( \Psi _{x}(\theta )\right) \right) +6232$

\strut

Now, let's compute $\underset{\theta \in NC_{4,4}(8)}{\sum }\mu _{\theta
}^{1_{8}}\cdot \tau \left( \Psi _{x}(\theta )\right) $ ; we have that

\strut

\ \ $\ NC_{4,4}(8)=\{\{(1,2,3,4),(5,6,7,8)\},\{(1,6,7,8),(2,3,4,5)\}$

$\ \ \ \ \ \ \ \ \ \ \ \ \ \ \ \ \ \ \ \ \ \ \ \ \ \ \
\{(1,2,7,8),(3,4,5,6)\},\{(1,2,3,8),(4,5,6,7)\}\}.$

\strut

Hence,

\strut

$\ \Psi _{x}\left( \{(1,2,3,4),(5,6,7,8)\}\right) =\Psi _{x}\left(
1,2,3,4\right) \Psi _{x}(5,6,7,8)$

$\ \ \ \ \ \ \ \ =E(x^{4})\cdot E(x^{4})=\left( (h+h^{-1})+p_{0}^{4}\right)
\left( (h+h^{-1})+p_{0}^{4}\right) $

$\ \ \ \ \ \ \ \ =(h+h^{-1})^{2}+2p_{0}^{4}(h+h^{-1})+(p_{0}^{4})^{2}$

$\ \ \ \ \ \ \ \ =h^{2}+2e+h^{-2}+2p_{0}^{4}(h+h^{-1})+784e$

$\ \ \ \ \ \ \ \ =h^{2}+h^{-2}+2p_{0}^{4}(h+h^{-1})+786e,$

\strut

$\ \Psi _{x}\left( \{(1,6,7,8),(2,3,4,5)\}\right) =\Psi _{x}\left(
1,[4],3\right) $

$\ \ \ \ \ \ \ \ =\left( (h^{2}+h^{-2})+2e\right) +784e$

$\ \ \ \ \ \ \ \ =(h^{2}+h^{-2})+786e,$

\strut

$\ \Psi _{x}\left( \{(1,2,7,8),(3,4,5,6)\}\right) =\Psi _{x}\left(
2,[4],2\right) $

$\ \ \ \ \ \ \ \ =(h^{2}+h^{-2})+786e$

\strut

and

\strut

$\ \Psi _{x}\left( \{(1,2,3,8),(4,5,6,7)\}\right) =\Psi _{x}\left(
3,[4],1\right) $

$\ \ \ \ \ \ \ \ =(h^{2}+h^{-2})+786e.$

\strut

Notice that, for any $\pi \in NC_{4,4}(8),$ we have that

$\strut $

\begin{center}
$\mu _{\pi }^{1_{8}}=\mu (0_{1},1_{1})^{6}\mu (0_{2},1_{2})=-1.$
\end{center}

\strut

Thus,

\strut

$\ \underset{\theta \in NC_{4,4}(8)}{\sum }\mu _{\theta }^{1_{8}}\cdot \tau
\left( \Psi _{x}(\theta )\right) =-\tau \left(
h+h^{-1}+2p_{0}^{4}(h+h^{-1})+786e\right) $

$\ \ \ \ \ \ \ \ \ \ \ \ \ \ \ \ \ \ \ \ \ \ \ \ \ \ \ \ \ \ \ \ \ \ \ \ \ \
\ \ \ \ \ \ \ \ \ \ \ \ \ \ \ -3\cdot \tau \left( (h^{2}+h^{-2})+786e\right) 
$

\strut

$\ \ \ \ \ \ \ \ \ \ \ \ \ \ \ \ \ \ \ \ \ \ \ \ \ \ \ \ \ \ \ \ \ \ \ \ \ \
\ \ \ \ =(-4)(786)=-3144.$

\strut

Therefore,

\strut

$\ \ \ \ \ \ \ \ \ \ 2\cdot \underset{\theta \in NC^{(even)}(8)}{\sum }\mu
_{\theta }^{1_{8}}\cdot \tau \left( \Psi _{x}(\theta )\right) $

$\ \ \ \ \ \ \ \ \ \ \ \ \ \ \ \ \ \ \ \ \ \ \ \ \ \ \strut =2\left( 
\underset{\theta \in NC_{4,4}(8)}{\sum }\mu _{\theta }^{1_{8}}\cdot \tau
\left( \Psi _{x}(\theta )\right) \right) +6232$

$\ \ \ \ \ \ \ \ \ \ \ \ \ \ \ \ \ \ \ \ \ \ \ \ \ =2\left( -3144\right)
+6232$

$\ \ \ \ \ \ \ \ \ \ \ \ \ \ \ \ \ \ \ \ \ \ \ \ \ =-6288+6232=-56$

\strut

\strut

(\textbf{Step 6}) By (\textbf{Step 1}) $\sim $ (\textbf{Step 5}), we can
conclude that

\strut

$\ \ \ \ \ \ \ \ \ \ \ \ \ \ \ \ \ \ \tau \left( (x+y)^{8}\right) =\tau
\left( E\left( (x+y)^{8}\right) \right) $

\strut

$\ \ \ \ \ \ \ \ \ \ \ \ \ \ \ \ \ \ \ \ \ \ \ \ \ \ \ \ \ \ \ \ \ \ \ \ \ \
=57344+(-14336)+1024+288+(-56)$

\strut

$\ \ \ \ \ \ \ \ \ \ \ \ \ \ \ \ \ \ \ \ \ \ \ \ \ \ \ \ \ \ \ \ \ \ \ \ \ \
=44264.$

\strut

\begin{example}
Let $x,y\in L(F_{2})*_{L(F_{1})}L(F_{2})$ be $L(F_{1})$-valued random
variables such that $x=a+b+a^{-1}+b^{-1}$ and $y=c+d+c^{-1}+d^{-1}$ and let $%
E:L(F_{2})*_{L(F_{2})}L(F_{2})\rightarrow L(F_{1})$ and $\tau
:L(F_{2})*_{L(F_{1})}L(F_{2})\rightarrow \Bbb{C}$ be the conditional
expectation (finding $h$-terms) and the canonical trace, respectively. Then

\strut 

\begin{center}
$\tau \left( (x+y)^{8}\right) =\tau \left( E\left( (x+y)^{8}\right) \right)
=44256.$
\end{center}

$\square $
\end{example}

\strut

\begin{remark}
The above result also gotten from the following way ; First, recall that $%
L(F_{2})*_{L(F_{1})}L(F_{2})\simeq L(F_{2}*_{F_{1}}F_{2}).$ Also, we can
regard the group $F_{2}*_{F_{1}}F_{2}$ as a (topological) fundamental group
of torus with genus 2,

\strut 

\begin{center}
$G=\,<a,b,c,d:aba^{-1}b^{-1}d^{-1}c^{-1}dc=e>.$
\end{center}

\strut \strut 

(Remember that $aba^{-1}b^{-1}=cdc^{-1}d^{-1}$ is our $h$ !) We need to
recognize that $aba^{-1}b^{-1}d^{-1}c^{-1}dc$ is a word with length 8,
without considering the relation in the group $G.$ Again, denote $%
a+b+a^{-1}+b^{-1}$ and $c+d+c^{-1}+d^{-1}$ by $x$ and $y,$ respectively.
Now, define the following trace

\strut 

\begin{center}
$\tau _{4}$ $:L(F_{4})\rightarrow \Bbb{C}$
\end{center}

by

\begin{center}
$\tau _{4}\left( \underset{g\in F_{4}}{\sum }\alpha _{g}g\right) =\alpha
_{e_{F_{4}}},$
\end{center}

\strut 

where $F_{4}=\,<a,b,c,d>.$ (Of course, in our case, the word $%
aba^{-1}b^{-1}d^{-1}c^{-1}dc$ is $e,$ but this can be come from making the
words with length 8 !) Now, let's compute $\tau \left( (x+y)^{8}\right) .$
This can be computed by using the method introduced in [35], as follows ;
this method is also used in Section 3.3.

\strut \strut 

\begin{center}
$\ \tau _{4}\left( (x+y)^{8}\right) =44284.$
\end{center}

\strut 

By the relation\strut 

\begin{center}
$aba^{-1}b^{-1}d^{-1}c^{-1}dc=e$
\end{center}

and

\begin{center}
$c^{-1}d^{-1}cdbab^{-1}a^{-1}=e^{-1}=e,$
\end{center}

\strut 

we have to add 16 to $\tau \left( (x+y)^{8}\right) .$ i.e

\strut 

\begin{center}
$
\begin{array}{ll}
\tau \left( (x+y)^{8}\right)  & =\tau \left( E\left( (x+y)^{8}\right)
\right)  \\ 
& =\tau _{4}\left( (x+y)^{8}\right) +16 \\ 
& =44264.
\end{array}
$
\end{center}
\end{remark}

\strut

The above method introduced in the previous remark looks much more easy to
compute the moments of $x+y.$ However, when we deal with the higher degree
computation, it is very hard to find the suitable relation for the identity $%
e.$

\strut

\strut \strut

\strut

\strut \textbf{References}

\strut

\label{REF}

{\small [1] \ A. Nica, R-transform in Free Probability, IHP course note,
available at www.math.uwaterloo.ca/\symbol{126}anica.\strut }

{\small [1] \ A. Nica, R-transforms of Free Joint Distributions and
Non-crossing Partitions, J. of Func. Anal, 135 (1996), 271-296.\strut }

{\small [24] \ A. Nica, R-diagonal Pairs Arising as Free Off-diagonal
Compressions, available at www.math.uwaterloo.ca/\symbol{126}anica.}

{\small [4] \ A. Nica, D. Shlyakhtenko, F. Goodman, Free Probability of Type
B, preprint.\strut \strut }

{\small [5] \ A. Nica, D. Shlyakhtenko and R. Speicher, R-cyclic Families of
Matrices in Free Probability, J. of Funct Anal, 188 (2002), 227-271.\strut }

{\small [6] \ A. Nica, D. Shlyakhtenko and R. Speicher, R-diagonal Elements
and Freeness with Amalgamation, Canad. J. Math. Vol 53, Num 2, (2001)
355-381.\strut }

{\small [7] \ A. Nica, D. Shlyakhtenko, R. Speicher, Operator-Valued
Distributions I. Characterizations of Freeness, preprint.}

{\small [8] \ A. Nica and R. Speicher, R-diagonal Pair-A Common Approach to
Haar Unitaries and Circular Elements, (1995), www.mast.queensu.ca/\symbol{126%
}speicher.\strut }

{\small [9] \ A. Nica and R.Speicher, A ''Fouries Transform'' for
Multiplicative Functions on Noncrossing Patitions, J. of Algebraic
Combinatorics, 6, (1997) 141-160.\strut }

{\small [13] B.Krawczyk and R.Speicher, Combinatorics of Free Cumulants,
available at www.mast.queensu.ca/\symbol{126}speicher.\strut \strut }

{\small [11] D. Shlyakhtenko, Some Applications of Freeness with
Amalgamation, J. Reine Angew. Math, 500 (1998), 191-212.\strut }

{\small [16] D. Shlyakhtenko, A-Valued Semicircular Systems, J. of Funct
Anal, 166 (1999), 1-47.\strut }

{\small [13] D. Voiculescu, Operations on Certain Non-commuting
Operator-Valued Random Variables, Ast\'{e}risque, 232 (1995), 243-275.\strut 
}

{\small [14] D.Voiculescu, K. Dykemma and A. Nica, Free Random Variables,
CRM Monograph Series Vol 1 (1992).\strut }

{\small [15] F. Radulescu, Singularity of the Radial Subalgebra of }$%
L(F_{N}) ${\small \ and the Puk\'{a}nszky Invariant, Pacific J. of Math,
vol. 151, No 2 (1991)\strut , 297-306.\strut \strut }

{\small [16] I. Cho, Amalgamated Boxed Convolution and Amalgamated
R-transform Theory (preprint).\strut }

{\small [17] I. Cho, Compressed Amalgamated R-transform Theory,
preprint.\strut }

{\small [18] I. Cho, Perturbed R-transform Theory, preprint.}

{\small [19] I. Cho, Compatibility of a noncommutative probability space and
an amalgamated noncommutative probability space, preprint}

{\small [20] J.A.Mingo and A.Nica, Annular Noncrossing Permutations and
Partitions and Second-order Asymptotics for Random Variables,
preprint.\strut }

{\small [21] M. Bo\.{z}ejko, M. Leinert and R. Speicher, Convolution and
Limit Theorems for Conditionally Free Random Variables, www.mast.queensu.ca/%
\symbol{126}speicher.\strut }

{\small [22] P.\'{S}niady and R.Speicher, Continous Family of Invariant
Subspaces for R-diagonal Operators, Invent Math, 146, (2001) 329-363.\strut }

{\small [23] R. Speicher, Combinatorial Theory of the Free Product with
Amalgamation and Operator-Valued Free Probability Theory, AMS Mem, Vol 132 ,
Num 627 , (1998).}

{\small [24] R. Speicher, Combinatorics of Free Probability Theory IHP
course note, available at www.mast.queensu.ca/\symbol{126}speicher.\strut }

{\small [25] R. Speicher, A Conceptual Proof of a Basic Result in the
Combinatorial Approach to Freeness, www.mast.queensu.ca/\symbol{126}%
speicher.\strut }

{\small \strut [26] S. Popa, Markov Traces on Universal Jones Algebras and
Subfactors of Finite Index, Invent. Math, 111 (1993), 375-405.\strut }

\strut

\strut

\strut

\strut

\strut

\end{document}